\documentclass{amsart}
\usepackage{amsmath}
\usepackage{amsfonts}
\usepackage{amssymb}
\usepackage{graphicx}
\usepackage{multirow}
\usepackage{stix}
\usepackage{courier}

\numberwithin{equation}{section}
\begin{document}

\title{On the geometry of a triangle \protect\\ in the elliptic and in the extended hyperbolic plane}
\author{Manfred Evers}
\curraddr[Manfred Evers]{Bendenkamp 21, 40880 Ratingen, Germany}
\email[Manfred Evers]{manfred\_evers@yahoo.com}
\date{\today}

\begin{abstract}
We investigate several topics of triangle geometry in the elliptic and in the extended hyperbolic plane. For both planes a uniform metric is used. The concept for this metric was developed by C. V\"or\"os at the beginning of the $20^{\textrm{th}}$ century and subsequently described by Horv\'ath in 2014.
\end{abstract}

\maketitle

\section*{Introduction}
This paper is a continuation of a previous work \cite{Ev}. Whilst the previous paper was restricted to triangles in the elliptic plane, these investigations also deal with triangles in the extended hyperbolic plane. In addition, both papers differ in the selection of centers, central lines, central conics and cubics. \\
\noindent\hspace*{5mm} The first section gives an introduction to the metric used on the projective plane. We assume that the reader is familiar with euclidean triangle geometry, but in order to introduce the terminology and fix notations, we give some basic definitions, rules and theorems. There are many introductory books on euclidean triangle geometry; we refer to presentations of Yiu \cite{Y} and Douillet \cite{D}. For geometry on the sphere / elliptic geometry we recommend the book \cite{TL} by Todhunter and Leathem. Several publications appeared in the last decade dealing with triangle geometry in the hyperbolic plane, to name Ungar \cite{U1, U2}, Wildberger \cite{W1}, Wildberger and Alkhaldi \cite{W2}, Horv\'ath \cite{H1, H2}, Vigara \cite{Vi} and Russell \cite{Ru}.\\
\noindent\hspace*{5mm} The topics of the second section are: centers based on orthogonality, centers related to circumcircles and incircles, radical centers and centers of similitude, orthology, Kiepert perspectors and related objects, Tucker circles, isoptics, and substitutes for the Euler line.\vspace*{0.5 mm}\\
Remark: A detailed proof is not given for every statement, but quite often the statement can be checked by simple computation (supported by a CAS-program). 

\section{Metric geometry in the projective plane \vspace*{0.5 mm}}

\subsection{} \textbf{The projective plane, its points and its lines}\hspace*{\fill} \\
\noindent\hspace*{5mm}Let $V$ be the three dimensional vector space $\mathbb{R}^{3}$, equipped with the canonical dot product $\boldsymbol{p} \cdot \boldsymbol{q} = (p_1,p_2,p_3) \cdot (q_1,q_2,q_3) = p_1q_1+p_2q_2+p_3q_3$, and let
$\mathcal{P}$ denote the projective plane $({V}{-}\{\boldsymbol{0}\})/\mathbb{R}^\times$. The image of a non-zero vector $\boldsymbol{p} = (p_1,p_2,p_3) \in V $ under the canonical projection $\Pi{:\,} V \rightarrow \mathcal{P}$ will be denoted by $(p_1{:}p_2{:}p_3)$ and will be regarded as a point in this plane. \\
\noindent\hspace*{5mm}Given two different points  $P$ and $Q$ in this projective plane, there exists exactly one line incident with these two points. It is called the \textit{join} $P \vee Q$ of $P$ and $Q$.  If $\boldsymbol{p} = (p_1,p_2,p_3)$, $\boldsymbol{q} = (q_1,q_2,q_3)$ are two non zero vectors with $\Pi(\boldsymbol{p}) = P$ and $\Pi(\boldsymbol{q}) = Q$, then the line $P \vee Q$ through $P$ and $Q$ is the set of points $\Pi(s \boldsymbol{p} + t \boldsymbol{q})$ with $s, t \in \mathbb{R}$.
One can find linear forms $l \in V^{*}{-}\{\boldsymbol{0^*}\}$ with $\text{ker}(l) = \text{span}(\boldsymbol{p}, \boldsymbol{q})$. A suitable $l$ is, for example, $l = \,\ast\!(\boldsymbol{p}\times \boldsymbol{q}) = (\boldsymbol{p}\times \boldsymbol{q})^*$, where $\times$ stands for the canonical cross product in $V = \mathbb{R}^{3}$ and $\ast$ for the isomorphism $V \rightarrow V^*, {\ast}(\boldsymbol{r}) = \boldsymbol{r}{\cdot}(.).$ The linear form $l$ is uniquely determined up to a nonzero real factor, so there is a ${1{:}1}$-correspondence between the lines in the projective plane and the elements of $\mathcal{P}^* = (V^{*}{-}\{\boldsymbol{0^*}\})/\mathbb{R}^\times$. We identify the line $l = P \vee Q$ with the element $(p_2q_3 - p_3q_2:p_3q_1 - p_1q_3:p_1q_3 - p_2q_1)^* \in \mathcal{P}^*.$\\
\noindent\hspace*{5mm}In the projective plane, two different lines $k = (k_1{:}k_2{:}k_3)^*, l = (l_1{:}l_2{:}l_3)^*$ always meet at one point $k \wedge l = \Pi((k_1,k_2,k_3)\times (l_1,l_2,l_3))$, the so-called \textit{meet} of these lines. 
\subsection{} \textbf{Collineations, correlations and polarities}\hspace*{\fill} \\
\noindent\hspace*{5mm}The automorphism group Aut$(\mathcal{P})$ of $\mathcal{P}$ can be identified with the projective linear group PGL$(V)$ = GL$(V)/\mathbb{R}^\times$. These automorphisms on $\mathcal{P}$ are called \textit{collineations} because the image of a line under an automorphism is again a line.
A \textit{correlation} is a bijective mapping $c{:}\;\mathcal{P}\cup \mathcal{P}^* \rightarrow \mathcal{P}\cup \mathcal{P}^*$ with the following property:
The restriction $c\vert _\mathcal{P}$ of $c$ to $\mathcal{P}$ is a point-to-line transformation that maps collinear points to concurrent lines, while $c\vert _\mathcal{P^*}$is a line-to-point transformation that maps concurrent lines to collinear points. Thus there are for each correlation $c$ uniquely determined linear mappings $c_1{:\,}\mathcal{P}\rightarrow\mathcal{P}^*,\, c_2{:\,}\mathcal{P}^*\rightarrow\mathcal{P}$ with $c_1 = c\vert _\mathcal{P}$ and $c_2=c\vert _{\mathcal{P}^*}.$
If $c_2$ is the inverse of $c_1$, this correlation is an involution. In this case, $c_1, c_2$ are symmetric mappings since $\mathcal{P}$ and $\mathcal{P}^*$ have even dimensions,  and the correlation $c$ is a \textit{polarity}.\vspace*{0.5 mm} 

\subsection{Elliptic and hyperbolic metric structures on $\mathcal{P}$}\hspace*{\fill}\\
\noindent\hspace*{5mm}We now fix a nonzero real number $\rho$. 
%and consider the mapping $\mathfrak{S} \in \textrm{PGL}(V)$ which is given with respect to the canonical basis of {V} by the matrix $\text{diag}(\rho,1,1)$. 
%\begin{pmatrix} 1 &0 &0\\0 &\sigma &0\\0 &0 &\varepsilon \sigma\end{pmatrix} \in \textrm{PGL}(V)$ 
By $\delta$ we denote the polarity which
maps a point $P = (p_1{:}p_2{:}p_3)$ to its \textit{dual line} $P^\delta = (\rho p_1{:}p_2{:}p_3)^*$ and a line $l = (l_1{:}l_2{:}l_3)^*$ to its \textit{dual point} $l^\delta  = (l_1{:}\rho l_2{:}\rho l_3)$. 
\subsubsection{} We make use of $\delta$ to introduce orthogonality. A line $k = (k_1{:}k_2{:}k_3)^*$ is \textit{orthogonal} (or \textit{perpendicular}) to a line $l = (l_1{:}l_2{:}l_3)^*$  when the dual $k^\delta$ of $k$ is a point on  $l$; this is precisely when  $k_1 l_1 + \rho k_2 l_2 +\rho k_3 l_3 = 0$. Obviously, if $k$ is orthogonal to $l$, then $l$ is orthogonal to $k$. Two points $P = (p_1{:}p_2{:}p_3)$ and $Q = (q_1{:}q_2{:}q_3)$ are \textit{orthogonal} to each other if $\rho p_1 q_1 +  p_2 q_2 + p_3 q_3 = 0$. This is the case exactly when their dual lines $P^\delta$ and $Q^\delta$ are orthogonal. Self-orthogonal points and lines are called \textit{isotropic}.
The isotropic points form the so-called \textit{absolute conic} $\mathcal{C}_{\!abs}$. This is either the empty set, in which case the geometry is called \textit{elliptic}, or it is a proper conic and the geometry is called \textit{hyperbolic}.\hspace*{\fill}\\
\noindent\hspace*{5mm}Each symmetric real 3x3-matrix $\mathfrak{M}$ with $\textrm{det}\, \mathfrak{M}\ne 0$ determines a scalar product (symmetric bilinear form) 
 $\;{\scriptscriptstyle{[\mathfrak{M}]}}{:\,}V{\times}V{\,\rightarrow\,}\mathbb{R}$  by\vspace*{-1 mm}\\
\centerline{$(v_1,v_2,v_3)\,\scriptscriptstyle{[\mathfrak{M}]}$$\,(w_1,w_2,w_3)\, =\, \boldsymbol(v_1,v_2,v_3)\mathfrak{M}\begin{pmatrix} w_1 \\ w_2\\ w_3 \end{pmatrix}.$}\vspace*{1 mm}\\
The orthogonality of points and lines can be expressed with the help of the scalar products ${\scriptscriptstyle{[\mathfrak{S}]}}$ and ${\scriptscriptstyle{[\mathfrak{S}^{-1}]}}\,$ where $\mathfrak{S} = \text{diag}(\rho,1,1)$.
Two points $P = (p_1{:}p_2{:}p_3)$ and $Q = (q_1{:}q_2{:}q_3)$ are orthogonal precisely when $(p_1{,}p_2{,}p_3){\,\scriptscriptstyle{[\mathfrak{S}]}\,}(q_1{,}q_2{,}q_3) = 0$,\, two lines $k = (k_1{:}k_2{:}k_3)^*$ and $l = (l_1{:}l_2{:}l_3)^*$ are orthogonal precisely when $(k_1{,}k_2{,}k_3){\,\scriptscriptstyle{[\mathfrak{S}^{-1}]}\,}(l_1{,}l_2{,}l_3) = 0.$
\subsubsection{}
Consider some point $P = (p_1{:}p_2{:}p_3)$ and some line $l = (l_1{:}l_2{:}l_3)^*$ with $l \ne P^\delta.$ \\
The perpendicular from $P$ to $l$ is the line \vspace*{-1 mm}
\begin{equation*}
\;\;\;\text{perp}(l,P) := P \vee l^\delta = (p_2l_3{-}p_3l_2: \rho(p_3l_1{-}p_1l_3): \rho(p_1l_2{-}p_2l_1))^*, \vspace*{-1 mm}
\end{equation*}                                              
and the  point on $l$  orthogonal to $P$ is \vspace*{-1 mm}
\begin{equation*}
\text{perp}(P,l) := l \wedge P^\delta = (\rho (p_2l_3{-}p_3l_2): p_3l_1{-}p_1l_3: (p_1l_2{-}p_2l_1)).\vspace*{-1 mm}
\end{equation*} 
The line through $P$ \textit{parallel} to $l$ is\vspace*{0.5 mm}\\
\noindent\hspace*{17mm}$\text{par}(l,P) := \text{perp}(\text{perp}(l,P),P)$.\\
The point\vspace*{-1 mm}\\
\noindent\hspace*{17mm}$\text{ped}(P,l) := \text{perp}(\text{perp}(P,l),l)$,\vspace*{0.5 mm}\\ 
is the \textit{orthogonal projection} of $P$ on $l$, also called the \textit{pedal} of $P$ on $l$.
\subsubsection{Line segments and angles.}Define the function $\chi{:}\;V \rightarrow \{-1,0,1\}$ by\vspace*{-1 mm}  
\[\chi(p_0,p_1,p_2) = 
\begin{cases}
 \;\;\,0,\;\;\text{if }\, (p_1,p_2,p_3) = (0,0,0)\;, \\
 \;\;\,1,\;\;\text{if }\, (p_1,p_2,p_3) > (0,0,0)\; \textrm{with respect to the lexicographic order,} \\
-1,\;\;\text{if }\, (p_1,p_2,p_3) < (0,0,0)\; \textrm{with respect to the lexicographic order.} \\
\end{cases}\vspace*{-1 mm}\]
For an anisotropic point $P = (p_1{:}p_2{:}p_3)$ let $P^\circ \in V$ be the vector \vspace*{-2 mm}\\
\[\,P^\circ\!:= \frac{\chi(p_1,p_2,p_3)}{\sqrt{|(p_1,p_2,p_3){\,\scriptscriptstyle{[\mathfrak{S}]}\,}(p_1,p_2,p_3)|}} (p_1,p_2,p_3) =\frac{\chi(p_1,p_2,p_3)}{\sqrt{|\rho p_1^{ 2} + p_2^{ 2}+ p_3^{ 2}|}} (p_1,p_2,p_3).\vspace*{-2 mm}\\
\] 
$P^\circ$ is obviously uniquely determined by $P$. \vspace*{1 mm}\\
For any two different points $P$ and $Q$ we introduce two line segments $[P,Q]_{+}, [P,Q]_{-}$, these are the closures of the two connected components of the set $P \vee Q - \{P, Q\}$. If $R,S$ are two different anisotropic points in $[P,Q]_+$, then $[R,S]_+ = \{\Pi(s R^\circ{+\,}t S^\circ) |\, s,t{\,\in\,}\mathbb{R}, st \geq 0\}$ is a subset of $[P,Q]_+$, and if $R,S$ are two different anisotropic points in $[P,Q]_-$, then $[R,S]_- = \{\Pi(s R^\circ{+\,}t S^\circ) |\, s,t{\,\in\,}\mathbb{R}, st \leq 0\}$ is a subset of $[P,Q]_-$.\vspace*{-1.5 mm}\\
%\\ $[P,Q]_{\pm}$ has the property that $P, Q \in [P,Q]_\pm$ and for any two different anisotropic points $R, S \in [P,Q]_{\pm}$ the point $\Pi(R^\circ \pm S^\circ)$ is also a point in $[P,Q]_\pm$. If $P$ and $Q$ are anisotropic, then 
%$[P,Q]_{+} = \{\Pi(s P^\circ{+\,}t Q^\circ) |\, s,t{\,\in\,}\mathbb{R}, st \geq 0\}$ and $[P,Q]_{-} = \{\Pi(s P^\circ{+\,}t Q^\circ) |  \, s,t{\,\in\,}\mathbb{R},$  $ st\leq 0\}$.
%\\

We define angles as subsets of the pencil of lines through a point, which is the vertex of this angle: Given three noncollinear points $Q, R, S$, put\vspace*{-2 mm}
\begin{equation*}\vspace*{-1 mm}
\angle_{+} QSR := \{S \vee P |\, P \in [Q,R]_+\}\,,\;\angle_{-} QSR := \{S \vee P |\, P \in [Q,R]_-\}\,.
\end{equation*} 
\vspace*{-1 mm}The union of these two angles is the complete pencil of lines through $S$.
\subsubsection{The length of a segment and the measure of an angle.} \label{subsubsec:The length of a segment and the measure of an angle.}
We introduce the measure of segments and angles using concepts that were developed by C. V\"or\"os in Hungary at the beginning of the $20^{\textrm{th}}$ century and more recently (2014)   described by Horv\'ath in two papers \cite{H1, H2}, where he applied these to various configurations in hyperbolic geometry.\\
To each line segment $s$ and to each angle $\phi$ is assigned its measure $\mu(s)$ resp. $\mu(\phi)$ which is a complex number with a real part in $\bar{\mathbb{R}} = \mathbb{R} \cup \{-\infty,\infty\}$ and an imaginary part in the interval $[0,\pi]$. \vspace*{-3 mm}\\

Define for each anisotropic point $P = (p_1{:}p_2{:}p_3)$ a number $\varepsilon_{\!P} \in \{1,-i\}$ by\\
\noindent\hspace*{25mm} \hspace*{20 mm}$\varepsilon_{\!P} := 1/\sqrt{P^\circ {\scriptscriptstyle{[\mathfrak{S}}]} P^\circ}\;$\footnote{$^)$ We follow the convention $\sqrt{r \exp(\phi\,i)} = \sqrt{r}\exp(\frac{1}{2}\phi\,i)$ for $r \in \mathbb{R}^+ \textrm{and}\; 0\leq \phi<2\pi $.}$^).$\vspace*{1 mm}\\
%${\framebox[10pt]{$\scriptscriptstyle{\mathfrak{S}}$}} $
First we describe the measure $\mu([P,Q]_{\pm})$ of line segments with anisotropic endpoints $P,Q,$ $P \ne Q$. (The remaining cases will be treated afterwards.)\\Here, the function $\mu$ satisfies the following rules:\vspace*{1 mm}\\
\noindent\hspace*{5mm}(1)$\;\;\mu([P,Q]_{+})$ and $\mu([P,Q]_{-})$ are finite complex numbers with imaginary parts\\
\noindent\hspace*{11mm}in the interval $[\,0,\pi\,]$ satisfying\\ 
\noindent\hspace*{11mm}$\mu([P,Q]_{+}) + \mu([P,Q]_{-}) = \mu(P \vee Q) = \pi i$.\vspace*{1 mm}\\
\noindent\hspace*{5mm}(2)$\;\;\cosh(\mu[P,Q]_{+}) = \varepsilon_{\!P}\varepsilon_{\!Q} P^\circ{\,\scriptscriptstyle{[\mathfrak{S}]}\,}Q^\circ.$\vspace*{1 mm}\\
\noindent\hspace*{5mm}(3)$\;\;$If $R$ is an anisotropic inner point of $[P,Q]_+$,\\
\noindent\hspace*{11mm}then $\mu([P,R]_+) + \mu([R,Q]_+) = \mu([P,Q]_+).$\vspace*{1 mm}\\
\noindent\hspace*{5mm}(4)$\;\;Q = \text{perp}(P,P{\vee}Q)$\, precisely when \,$\mu([P,Q]_+) = \mu([P,Q]_-) = \frac{\pi}{2} i$.\vspace*{2 mm}\\
Let us look at special cases.\vspace*{1 mm}\\
In the elliptic case, all points $R$ are anisotropic with $\varepsilon_{\!R} = 1$. For two different points $P$ and $Q$ we have
$-1 < P^\circ{\,\scriptscriptstyle{[\mathfrak{S}]}\,}\,Q^\circ < 1$ and $\mu([P,Q]_{\pm}) = s i$ with a real number $s, \,0 < s < \pi.$ \vspace*{1 mm}\\
The hyperbolic case is more complicated. The set of anisotropic points in $\mathcal{P}$ consists of two connected components. One component is the part \textit{outside} the absolute conic; it consists of all nonisotropic points of tangents of $\mathcal{C}_{\!abs}$. The other component contains the points \textit{inside} the conic; each line through a point of this region meets the absolute conic twice. The inner part consists of points $R$ with $\varepsilon_{\!R} = -i$, 
the outer part contains all points $R$ with $\varepsilon_{\!R} = 1$.\vspace*{1 mm}\\ 
If all points on $P \vee Q$ are anisotropic, the situation is similar to the elliptic case. We have $\varepsilon_{\!R} = 1$ for each point $R$ on $P \vee Q$, and we have $-1 < P^\circ{\,\scriptscriptstyle{[\mathfrak{S}]}\,}\,Q^\circ < 1$ and $\mu([P,Q]_{\pm}) = s i$ with a real number $s$,  $0 < s < \pi.$  \vspace*{1 mm}\\
If $P \vee Q$ is isotropic, then it contains exactly one isotropic point, which is the touchpoint of this line with the absolute conic. If this point is outside the segment $[P,Q]_+$, then
$\varepsilon_{\!P}\varepsilon_{\!Q} P^\circ{\,\scriptscriptstyle{[\mathfrak{S}]}\,}\,Q^\circ = 1, \mu([P,Q]_+) = 0$, and $\mu([P,Q]_-) = \pi i$. Otherwise, $\varepsilon_{\!P}\varepsilon_{\!Q} P^\circ{\,\scriptscriptstyle{[\mathfrak{S}]}\,}\,Q^\circ = -1,\, \mu([P,Q]_+) = \pi i$, \,and $\mu([P,Q]_-) = 0$.\vspace*{1 mm}\\
The remaining case: There are two isotropic points on $P \vee Q$. We consider several subcases:\\
Subcase 1: $P$ and $Q$ are inside the absolute conic. Then there is no isotropic point in $[P,Q]_+$, while there are two such points in $[P,Q]_-$. The length $\mu([P,Q]_+)$ has to be a real number in order to satisfy the condition $\cosh(\mu[P,Q]_{+}) = \varepsilon_{\!P}\varepsilon_{\!Q} P^\circ{\,\scriptscriptstyle{[\mathfrak{S}]}\,}\,Q^\circ$, since 
the number on the right side of this equation is a real number $>1$. We have the choice of $\mu([P,Q]_+)$ being positive or negative and decide for the positive value.\vspace*{1 mm}\\
Subcase 2: One of the points $P, Q$ is inside, the other one is outside the absolute conic. Let us assume $P$ is outside; then the point $R := \text{perp}(P,P{\vee}Q)$ lies inside the absolute conic. In this case, \\
\noindent\hspace*{5mm}$\mu([P,Q]_+) = \,\mu([R,Q]_+) + \frac{\pi}{2} i, \; \mu([P,Q]_-) = \,-\mu([R,Q]_+) + \frac{\pi}{2} i$, \;if $R \in [P,Q]_+$, and\\
\noindent\hspace*{5mm}$\mu([P,Q]_+) = -\mu([R,Q]_+) + \frac{\pi}{2} i, \;\mu([P,Q]_-) = \mu([R,Q]_+) + \frac{\pi}{2} i$, \;if $R \in [P,Q]_-$. \vspace*{1 mm}\\
Subcase 3: Both, $P$ and $Q$, are points outside the absolute conic. With $R := \text{perp}(P,P{\vee}Q)$ and $S := \text{perp}(Q,P{\vee}Q)$ we get\vspace*{0.5mm}\\
\noindent\hspace*{15mm}$\cosh(\mu([P,Q]_+) = \varepsilon_{\!P}\varepsilon_{\!Q} P^\circ{{\,\scriptscriptstyle{[\mathfrak{S}]}\,}}\,Q^\circ = \varepsilon_{\!R}\varepsilon_{\!S} R^\circ{{\,\scriptscriptstyle{[\mathfrak{S}]}\,}}\,S^\circ = \cosh(\mu([R,S]_+)$.\vspace*{0.5mm}\\
The rules (1) and (4) require $\mu([P,Q]_+) = -\mu([R,S]_+) < 0$. \vspace*{1 mm}\\
An analysis of the different cases shows that by knowing the two numbers $\varepsilon_{\!P}\varepsilon_{\!Q}$ and $P^\circ{{\,\scriptscriptstyle{[\mathfrak{S}]}\,}}\,Q^\circ$ we can determine $\mu([P,Q]_\pm)$.\vspace*{3 mm}\\
We now set the lengths of the segments $[P,Q]_\pm$ if at least one endpoint is isotropic. In this case, $\mu([P,Q]_{\pm})$ is not an element in $\mathbb{C}$, the real part of $\mu([P,Q]_\pm)$ is either $+\infty$ or $-\infty$. Again we "decline" various cases.\vspace*{1 mm}\\
If $P \vee Q$ is an isotropic line and $P$ is an isotropic point, then $Q$ is  anisotropic and \vspace*{0.5mm}\\
\noindent\hspace*{23mm}$P = \text{perp}(Q,P \vee Q)$,\, $\mu([P,Q]_+) = \mu([P,Q]_-) = \frac{\pi}{2} i$.\vspace*{1 mm}\\ 
If $P \vee Q$ is an anisotropic line and just one of the points $P, Q$, let us say $P$, is isotropic, then: \\
\noindent\hspace*{5mm}$\mu([P,Q]_+) = \infty + \frac{1}{4}\pi i$,\; $\mu([P,Q]_-) = -\infty + \frac{3}{4}\pi i$ in case of $\varepsilon_{\!Q}= -i$.\\
\noindent\hspace*{5mm}$\mu([P,Q]_+) = -\infty + \frac{1}{4}\pi i$,\; $\mu([P,Q]_-) = \infty + \frac{3}{4}\pi i$ \,if $\varepsilon_{\!Q}= 1\,$ and $\,P$ is the only isotropic\\ \noindent\hspace*{5mm}point in $[P,Q]_+$,\\
\noindent\hspace*{5mm}$\mu([P,Q]_+) = \infty + \frac{3}{4}\pi i$,\; $\mu([P,Q]_-) = -\infty + \frac{1}{4}\pi i$ \,if there are two isotropic points\\ \noindent\hspace*{5mm}in $[P,Q]_+$.\\
If $P$ and $Q$ are isotropic, then $\mu([P,Q]_{\pm}) = \pm \infty + \frac{\pi}{2} i$.\vspace*{2 mm}\\
The measure of angles:\; We use the same symbol $\mu$ for the measure of angles as for line segments. The angle $\angle_{+} QSR$ is either identical with $\{P^\delta\,|\,P\in [(S\vee Q)^\delta, (S \vee R)^\delta]_{+}\}$ or with $\{P^\delta\,|\,P\in [(S\vee Q)^\delta, (S \vee R)^\delta]_{-}\}$. In the first case we set $\mu(\angle_{\pm} QSR) = \mu([(S\vee Q)^\delta, (S \vee R)^\delta]_{\pm})$, in the second $\mu(\angle_{\pm} QSR) = \mu([(S\vee Q)^\delta, (S \vee R)^\delta]_{\mp})$. If the points $R, S, Q$ and the lines $S{\vee}Q,\, S{\vee}R$ are anisotropic, then \vspace*{1 mm}\\
%\[
%\begin{split}
%&\cosh(\mu(\angle_{+} QSR)) = -\cosh(\mu(\angle_{-} QSR))\\
%&= \frac{(S^\circ\times Q^\circ)\bullet(S^\circ\times R^\circ)}{\sqrt{(S^\circ\times Q^\circ)\bullet(S^\circ\times Q^\circ)}\,\sqrt{(S^\circ\times &R^\circ)\bullet(S^\circ\times R^\circ)}}\\
%&= \frac{(S^\circ\bullet S^\circ)(Q^\circ\bullet R^\circ) - (S^\circ\bullet Q^\circ)(S^\circ\bullet R^\circ)}{\sqrt{(S^\circ\bullet S^\circ)(Q^\circ\bullet Q^\circ) - (S^\circ\bullet Q^\circ)^2}\,\sqrt{(S^\circ\bullet S^\circ)(R^\circ\bullet R^\circ) - (S^\circ\bullet R^\circ)^2}}\\
%\end{split}
%\]
\centerline{$\cosh(\mu(\angle_{+} QSR)) = -\cosh(\mu(\angle_{-} QSR))$\hspace*{60 mm}}\vspace*{1 mm}\\
\centerline{$\displaystyle = \frac{(S^\circ\times Q^\circ){\,\scriptscriptstyle{[\mathfrak{S}^{-1}]}\,}(S^\circ\times R^\circ)}{\sqrt{(S^\circ\times Q^\circ){\,\scriptscriptstyle{[\mathfrak{S}^{-1}]}\,}(S^\circ\times Q^\circ)}\,\sqrt{(S^\circ\times R^\circ){\,\scriptscriptstyle{[\mathfrak{S}^{-1}]}\,}(S^\circ\times R^\circ)}}$\,\hspace*{29 mm}}\vspace*{1 mm}\\
\centerline{ $\displaystyle= \frac{(S^\circ{\,\scriptscriptstyle{[\mathfrak{S}]}\,} S^\circ)(Q^\circ{\,\scriptscriptstyle{[\mathfrak{S}]}\,} R^\circ) - (S^\circ{\,\scriptscriptstyle{[\mathfrak{S}]}\,} Q^\circ)(S^\circ{\,\scriptscriptstyle{[\mathfrak{S}]}\,} R^\circ)}{\sqrt{(S^\circ{\,\scriptscriptstyle{[\mathfrak{S}]}\,} S^\circ)(Q^\circ{\,\scriptscriptstyle{[\mathfrak{S}]}\,} Q^\circ) - (S^\circ{\,\scriptscriptstyle{[\mathfrak{S}]}\,} Q^\circ)^2}\,\sqrt{(S^\circ{\,\scriptscriptstyle{[\mathfrak{S}]}\,} S^\circ)(R^\circ{\,\scriptscriptstyle{[\mathfrak{S}]}\,} R^\circ) - (S^\circ{\,\scriptscriptstyle{[\mathfrak{S}]}\,} R^\circ)^2}}\;.$}\vspace*{1 mm}
\subsubsection{The distance between points and the (angle) distance between lines.} 
We introduce an order $\prec$ on $\bar{\mathbb{R}} + \mathbb{R} i$\; by\vspace*{-3 mm}\\ 
\centerline{$a_1 + b_1 i \prec a_2 + b_2 i \;\;\;\text{iff} \;\;\;
\begin{cases}
&b_1 < b_2 \vspace*{-2 mm}\\
\text{or}\vspace*{-2 mm}&\\
&b_1 = b_2 \;\;\text{and }\, a_1 < a_2\;,
\end{cases}$}\vspace*{2 mm}\\
and define a function $d{:\,}\mathcal{P}\times \mathcal{P} \rightarrow \{a{+}bi\, |\, a \in \bar{\mathbb{R}}, b \in [0,\frac{\pi}{2}] \}$ by\vspace*{1 mm}\\
\centerline{$d(P,Q) = 
\begin{cases}
0,\;\;\text{if } P = Q\;, \\
\mu([P,Q]_+),\;\text{if }\, P \ne Q \;\,\textrm{and}\;\, \mu([P,Q]_+ \prec  \mu([P,Q]_-\;, \\
\mu([P,Q]_-),\;\text{otherwise}\,.
\end{cases}$}\vspace*{1 mm}\\
We call $d$ a distance function on $\mathcal{P}$, even though it is obviously not a proper metric. This distance is continuous on $\mathcal{P} - \mathcal{C}_{\!abs}$.\vspace*{0.5mm}\\
\noindent\hspace*{5mm} For anisotropic points $P$ and $Q$, the distance $d(P,Q)$ is uniquely determined by the two numbers $\cosh^2(d(P,Q))$ and $\varepsilon_{\!P}\varepsilon_{\!Q}$.
If $P$ and $Q$ are both isotropic, then $\cosh^2(d(P,Q)) = -\infty$ and $d(P,Q) = -\infty + \frac{\pi}{2} i$. If only $P$ is isotropic, then $\cosh^2(d(P,Q)) = \frac{1}{2} -\varepsilon_{\!Q}^{ 2}\infty\,i$ and $d(P,Q) = -\varepsilon_{\!Q}^{2}\infty + \frac{\pi}{4} i$. \vspace*{0.5mm}\\
\noindent\hspace*{5mm} In the same way we define a function measuring the distance between two lines, and we also use the same symbol $d$. Thus, given two lines $k$ and $l$, then $d(k,l) = d(k^\delta,l^\delta)$. 
\subsubsection{Barycentric coordinates of a point on a line.} 
Let $P$ and $Q$ be two different anisotropic points on an anisotropic line and $R_\pm \in [P,Q]_\pm$ be two other anisotropic
%\footnote{$^)$ A closer analysis shows that the restriction of $R_\pm$ being anisotropic can be removed.}$^)\;$$\!$
points on this line, then $R_\pm = \Pi(x P^\circ \pm y{\,}Q^\circ)$ with any two real numbers $x, y$ satisfying\vspace*{1mm}\\
 \centerline{$(\ast)\;\;x\,:\,y = \sinh(d(R_\pm,Q))\,\varepsilon_{\!P}{\,:\,}\sinh(d(R_\pm,P))\,\varepsilon_{\!Q}$ .}\vspace*{1mm}\\
Real numbers $x, y$ satisfying $(\ast)\;$are called \textit{barycentric coordinates} of $P_\pm$ with respect to the points $P,Q$.\vspace*{1 mm}\\
\textit{Proof} of $(\ast):$ Let $R\,$ be an anisotropic point on $P\vee Q$. Assume $R^\circ = x P^\circ {+\,} y\,Q^\circ$, $x, y \in \mathbb{R}$. Then:\vspace*{0 mm}\\

\noindent\hspace*{22mm}$\displaystyle\sinh^2(d(P,R)) = \cosh^2(d(P,R)) - 1 $\vspace*{1mm}\\
\noindent\hspace*{44mm}$= \displaystyle\frac{(R^\circ{\,\scriptscriptstyle{[\mathfrak{S}]}\,}P^\circ)^2 \varepsilon_{\!P}^{\;2}}{R^\circ{\,\scriptscriptstyle{[\mathfrak{S}]}\,} R^\circ} - 1 $ \vspace*{1mm}\\ 
\noindent\hspace*{44mm}$\displaystyle= \frac{((x P^\circ{+\,}y\,Q^\circ){\,\scriptscriptstyle{[\mathfrak{S}]}\,}P^\circ)^2 \varepsilon_{\!P}^{\;2}}{(x P^\circ{+\,}y\,Q^\circ){\,\scriptscriptstyle{[\mathfrak{S}]}}(x P^\circ{+\,}y\,Q^\circ)} - 1$\vspace*{1mm}\\
\noindent\hspace*{44mm}$=\displaystyle\frac{x^2\varepsilon_{\!P}^{\;2} + {2 xy P^\circ{\,\scriptscriptstyle{[\mathfrak{S}]}\,}Q^\circ} + y^2 (P^\circ{\,\scriptscriptstyle{[\mathfrak{S}]}\,}Q^\circ)^2\varepsilon_{\!P}^{\;2}} {x^2\varepsilon_{\!P}^{\;2} + {2 xy P^\circ{\,\scriptscriptstyle{[\mathfrak{S}]}\,}Q^\circ} + y^2\varepsilon_{\!Q}^{\;2}} - 1$\\ 
\noindent\hspace*{44mm}$= \displaystyle\frac{y^2({(P^\circ{\,\scriptscriptstyle{[\mathfrak{S}]}\,}Q^\circ})^2 \varepsilon_{\!P}^{\;2}-\varepsilon_{\!Q}^{\;2})}{{R^\circ{\,\scriptscriptstyle{[\mathfrak{S}]}\,} R^\circ}}.$ \\

\noindent\hspace*{0mm}In the same way we get  
$\hspace*{1mm}\displaystyle\sinh^2(d(Q,R)) = \frac{x^2({(P^\circ{\,\scriptscriptstyle{[\mathfrak{S}]}\,}Q^\circ})^2 \varepsilon_{\!Q}^{\;2}-\varepsilon_{\!P}^{\;2})}{{R^\circ{\,\scriptscriptstyle{[\mathfrak{S}]}\,} R^\circ}}\;.$\hspace*{0mm}\\

\noindent\hspace*{0mm}Thus:
$\displaystyle\frac{\varepsilon_{\!Q}^{\;2}\sinh^2(d(P,R))}{\varepsilon_{\!P}^{\;2}\sinh^2(d(Q,R))} = \frac{y^2({(P^\circ{\,\scriptscriptstyle{[\mathfrak{S}]}\,}Q^\circ})^2 \varepsilon_{\!P}^{\;2}\varepsilon_{\!Q}^{\;2}-1)}{x^2({(P^\circ{\,\scriptscriptstyle{[\mathfrak{S}]}\,}Q^\circ})^2 \varepsilon_{\!P}^{\;2}\varepsilon_{\!Q}^{\;2}-1)} =\frac{y^2(\cosh^2(d(P,Q))-1)}{x^2(\cosh^2(d(P,Q))-1)}  = \frac{y^2}{x^2}\;.$\vspace*{2mm}\\

\noindent Now we look at different cases.\\
If all points of $P \vee Q$ are anisotropic, then $\varepsilon_{\!P}=\varepsilon_{\!Q}=\varepsilon_{\!R}=1$ and $d(P,R), d(Q,R) \in [0,\frac 12 \pi i]$. Therefore, $(x,y) \in \mathbb{R} (-i\varepsilon_{\!P}\sinh(d(Q,R)), -i\varepsilon_{\!Q}\sinh(d(P,R))).$\hspace*{0.3mm}\\
If $P\vee Q$ is a line passing through two isotropic points, we consider four subcases:\\
$(1)\,: \varepsilon_{\!P}=\varepsilon_{\!Q}=\varepsilon_{\!R}$. Then $d(P,R), d(Q,R) \in \mathbb{R}$ and \\ 
\hspace*{8 mm}$(x,y) \in \mathbb{R} (\varepsilon_{\!P}\sinh(d(Q,R)), \varepsilon_{\!Q}\sinh(d(P,R)))$ if $\varepsilon_{\!P} = 1$,\\
\hspace*{8 mm}$(x,y) \in \mathbb{R} (-i\varepsilon_{\!P}\sinh(d(Q,R)), -i\varepsilon_{\!Q}\sinh(d(P,R)))$ if $\varepsilon_{\!P} = -i$.\\
$(2):$ $\varepsilon_{\!P}=\varepsilon_{\!Q}\ne\varepsilon_{\!R}$, hence $d(P,R), d(Q,R) \in \frac12 \pi i +\mathbb{R}$.\\
$\hspace*{9mm}$If $\varepsilon_{\!P}{\,=\,}\varepsilon_{\!Q}{\,=\,}1$, then  $(x,y) \in $ $\mathbb{R} (-i\varepsilon_{\!P}\sinh(d(Q,R)), -i\varepsilon_{\!Q}\sinh(d(P,R)))$.\\
$\hspace*{9mm}$If $\varepsilon_{\!P}=\varepsilon_{\!Q} = -i$, then $(x,y) \in \mathbb{R}(\varepsilon_{\!P}\sinh(d(Q,R)),\varepsilon_{\!Q}\sinh(d(P,R))).$\\
$(3):$ If $\varepsilon_{\!P}=\varepsilon_{\!R}\ne\varepsilon_{\!Q}$, then  $d(P,R) \in \mathbb{R}\;$ and $\;d(Q,R) \in \frac12 \pi i +\mathbb{R}$.\\
$\hspace*{9mm}$In this case, $(x,y) \in $ $\mathbb{R} (-i\varepsilon_{\!P}\sinh(d(Q,R)), -i\varepsilon_{\!Q}\sinh(d(P,R)))$ if $\varepsilon_{\!P} = 1$,\\
$\hspace*{9mm}$and $(x,y) \in \mathbb{R}(\varepsilon_{\!P}\sinh(d(Q,R)),\varepsilon_{\!Q}\sinh(d(P,R)))$ if $\varepsilon_{\!P} = -i$.\\
$(4):$ If $\varepsilon_{\!P}\ne\varepsilon_{\!R}=\varepsilon_{\!Q}$, then  $d(P,R) \in \frac12 \pi i +\mathbb{R}\;$ and $\;d(Q,R) \in \mathbb{R}$.\\
$\hspace*{9mm}$In this case, $(x,y) \in \mathbb{R}(\varepsilon_{\!P}\sinh(d(Q,R)),\varepsilon_{\!Q}\sinh(d(P,R)))$ if $\varepsilon_{\!P} = 1$,\\
$\hspace*{9mm}$and $(x,y) \in $ $\mathbb{R} (-i\varepsilon_{\!P}\sinh(d(Q,R)), -i\varepsilon_{\!Q}\sinh(d(P,R)))$ if $\varepsilon_{\!P} = -i. \;\;\;\;\;\;\;\;\Box\vspace*{0 mm}$\vspace*{0.5mm}\\
\subsubsection{Semi-midpoints.} 
The points $R_\pm = P \pm Q$ $:=\Pi(P^\circ \pm Q^\circ)$ we call \textit{semi-midpoints} of $P$ and $Q$. These points were introduced and investigated by Wildberger and Alkhaldi \cite{W2} under the name \textit{smydpoints}.\\
Properties of semi-midpoints:\\
If $\varepsilon_{\!P}=\varepsilon_{\!Q}$, then $P + Q$ and $P - Q$ are (proper) midpoints: $d(P \pm Q, P) = d(P \pm Q, Q)$.\\ 
If $\varepsilon_{\!P}\ne\varepsilon_{\!Q}$, then $d(P \pm Q, P) \ne d(P \pm Q, Q)$. So, $P + Q$ and $P - Q$ are not proper midpoints of $P$ and $Q$. We will call them \textit{pseudo-midpoints}; in \cite{W2} they are called \textit{sydpoints}. If $P \ne \text{perp}(Q, P\vee Q)$, then $P + Q$ and $P - Q$ are anisotropic and 
$d(P \pm Q, P)\varepsilon_{P \pm Q}\varepsilon_{\!P} = d(P \mp Q, Q)\varepsilon_{P \mp Q}\varepsilon_{\!Q}$. If $P = \text{perp}(Q, P\vee Q)$, then $P + Q$ and $P - Q$ are isotropic and $d(P\pm Q, P) = -\varepsilon_{\!P}^{2}\infty + \frac{\pi}{4} i\,,\;d(P\pm Q, Q) = -\varepsilon_{\!Q}^{2}\infty + \frac{\pi}{4} i$.\\
$P, Q, P{\,+\,}Q, P{\,-\,}Q$ always form a harmonic range.\newpage
\subsection{The use of barycentric coordinates.}\vspace*{-2 mm}
\subsubsection{Barycentric coordinates with respect to a re\-ference triple $\Delta$ of points.} We fix a re\-ference triple $\Delta = ABC$ of non-collinear, anisotropic points $A, B, C \in \mathcal{P}$. For every point $P$ we can find a triple $(p_1,p_2,p_3)$ of real numbers such that $P = p_1A + p_2B + p_3C := \Pi(p_1A^\circ + p_2B^\circ + p_3C^\circ)$. Such a triple of real numbers will be called (triple of) \textit{barycentric coordinates} of $P$ with respect to $\Delta$. The triple $(p_1,p_2,p_3)$ is not uniquely determined by $P$, but every other triple of barycentric coordinates is a real multiple of $(p_1,p_2,p_3)$. The point $P$ is determined by the homogenous triple $(p_1{:}p_2{:}p_3)$ and $\Delta$, and we write $P = [p_1{:}p_2{:}p_3]_\Delta$. The triple can be calculated by $p_1{:}p_2{:}p_3 = P^\circ{\cdot}\,(B^\circ\,{\times}\,C^\circ):P^\circ{\cdot}\,(C^\circ{\times}\,A^\circ):\!P^\circ{\cdot}\,(A^\circ{\times}\,B^\circ)$.\\
\textit{Remark}: As barycentric coordinates we also accept a triple of complex numbers as long as it is a complex multiple of a real triple. 
\subsubsection{Lines.}Let $P=[p_1\,{:}\,p_2\,{:}\,p_3]_\Delta$ and $Q=[p_1\,{:}\,p_2\,{:}\,p_3]_\Delta$ be two different points, then 
\noindent\hspace*{0mm}the line $P\,{\vee}\,Q$ consists of all points $X=[x_1\,{:}\,x_2:x_3]_\Delta$ satisfying the equation \vspace*{0.5 mm}\\
\centerline{$((p_1,p_2,p_3) \times (p_1,p_2,p_3)) \cdot (x_1,x_2,x_3) = 0.$}\vspace*{2 mm}\\
If $R \vee S$ is a line through $R = [r_1\!:\!r_2\!:\!r_3]_\Delta$ and $S = [s_1\!:\!s_2\!:\!s_3]_\Delta$, different from $P \vee Q$, then both lines meet at a point 
$T = [t_1\!:\!t_2\!:\!t_3]_\Delta$ with \vspace*{1 mm} \\
\centerline{$
(t_1,t_2,t_3) = ((p_1,p_2,p_3)\times ((q_1,q_2,q_3))\times((r_1,r_2,r_3)\times (s_1,s_2,s_3)).$
}
\subsubsection{Triangles.} The reference triple $\Delta = ABC$ determines four \textit{triangles} that share the same vertices $A, B, C$ and the same sidelines $a:= B \vee C$, $b:= C \vee A$, $c:= A \vee B$. These triangles are the closures in $\mathcal{P}$ of the four connected components of $\mathcal{P} - \{a,b,c\}$:\vspace*{1 mm}\\
\noindent\hspace*{5mm}$\Delta_0 := \{[p_1:p_2:p_3]_\Delta |\, p_1,p_2,p_3 > 0\} \cup [B,C]_+ \cup [C,A]_+ \cup [A,B]_+$,\\
\noindent\hspace*{5mm}$\Delta_1 := \{[p_1:p_2:p_3]_\Delta |\,{-}p_1,p_2,p_3 > 0\} \cup [B,C]_+ \cup [C,A]_- \cup [A,B]_-$,\\
\noindent\hspace*{5mm}$\Delta_2 := \{[p_1:p_2:p_3]_\Delta |\, p_1,-p_2,p_3 > 0\} \cup [B,C]_- \cup [C,A]_+ \cup [A,B]_-$,\\
\noindent\hspace*{5mm}$\Delta_3 := \{[p_1:p_2:p_3]_\Delta |\, p_1,p_2,-p_3 > 0\} \cup [B,C]_- \cup [C,A]_- \cup [A,B]_+$.\vspace*{2 mm}\\
Remark: We do not consider the closures of the sets $\mathcal{P} - \Delta_k ,\, k = 0,1,2,3,\;$ as triangles.\vspace*{0.5 mm}
\subsubsection{Note:}\label{subsubsec:Note1}In the following, we always assume that not only the vertices but also the sidelines of $\Delta$ are anisotropic.\vspace*{0 mm}
\subsubsection{The characteristic matrix of $\Delta$.}\label{subsubsec:characteristic matrix} We introduce the symmetric matrix \vspace*{1.5 mm}\\
\centerline{$
\mathfrak{C} = (\mathfrak{c}_{ij})_{i,j=1,2,3}:= \begin{pmatrix} A^\circ{\,\scriptscriptstyle{[\mathfrak{S}]}\,} A^\circ &A^\circ\!{\,\scriptscriptstyle{[\mathfrak{S}]}\,} B^\circ &A^\circ\!{\,\scriptscriptstyle{[\mathfrak{S}]}\,} C^\circ\\B^\circ\!{\,\scriptscriptstyle{[\mathfrak{S}]}\,} A^\circ &B^\circ\!{\,\scriptscriptstyle{[\mathfrak{S}]}\,} B^\circ &B^\circ \!{\,\scriptscriptstyle{[\mathfrak{S}]}\,} C^\circ \\C^\circ\!{\,\scriptscriptstyle{[\mathfrak{S}]}\,} A^\circ &C^\circ\!{\,\scriptscriptstyle{[\mathfrak{S}]}\,} B^\circ &C^\circ\!{\,\scriptscriptstyle{[\mathfrak{S}]}\,} C^\circ \end{pmatrix}.\vspace*{1.5 mm}
$}
$\mathfrak{C}$ is a regular matrix with $\det\, \mathfrak{C} = \rho\, (A^\circ{\cdot}\,(B^\circ\times C^\circ))^2 \ne 0$ and inverse \vspace*{1.5 mm}\\
\centerline{$\displaystyle
\mathfrak{D} = (\mathfrak{d}_{ij})_{i,j=1,2,3} = \frac{1}{\det\, \mathfrak{C}} \begin{pmatrix} \mathfrak{c}_{22}\mathfrak{c}_{33}-\mathfrak{c}_{23}^2&\mathfrak{c}_{23}\mathfrak{c}_{31}-\mathfrak{c}_{12}\mathfrak{c}_{33}&\mathfrak{c}_{12}\mathfrak{c}_{23}-\mathfrak{c}_{31}\mathfrak{c}_{22}\\\mathfrak{c}_{23}\mathfrak{c}_{31}-\mathfrak{c}_{12}\mathfrak{c}_{33}&\mathfrak{c}_{33}\mathfrak{c}_{11}-\mathfrak{c}_{31}^2&\mathfrak{c}_{31}\mathfrak{c}_{12}-\mathfrak{c}_{23}\mathfrak{c}_{11}\\\mathfrak{c}_{12}\mathfrak{c}_{23}-\mathfrak{c}_{31}\mathfrak{c}_{22}&\mathfrak{c}_{31}\mathfrak{c}_{12}-\mathfrak{c}_{23}\mathfrak{c}_{11}&\mathfrak{c}_{11}\mathfrak{c}_{22}-\mathfrak{c}_{12}^2\end{pmatrix}.
$}\vspace*{1.5 mm}\\
$\mathfrak{C}$ is the Gram matrix of the basis $A^\circ,B^\circ,C^\circ$ of the vectorspace $V$ with the scalar product ${{\,\scriptscriptstyle{[\mathfrak{S}]}\,}}$. We call $\mathfrak{C}$ the characteristic matrix of $\Delta$, since from this matrix, as we will see, one can read off the lengths of the sides and the measures of the angles of the triangles $\Delta_k, \, k = 0,1,2,3$.\\\vspace*{-2 mm}

Let us denote the lengths of the sides of $\Delta_0$ by\vspace*{1 mm}\\ 
\centerline{$\mathscr{a}{:=}\mu([B,C]_+),\; \mathscr{b}{:=}\mu([C,A]_+),\; \mathscr{c}{:=}\mu([A,B]_+)$,}\\
then\\
\centerline{$\displaystyle \cosh(\mathscr{a}) = \varepsilon_{\!B}\varepsilon_{\!C} B^\circ{{\,\scriptscriptstyle{[\mathfrak{S}]}\,}}\,C^\circ = \frac{1}{\sqrt{\mathfrak{c}_{22}}}\frac{1}{\sqrt{\mathfrak{c}_{33}}}\mathfrak{c}_{23}$,\hspace*{12 mm}\;}\vspace*{1 mm}\\
\centerline{$\cosh(\mathscr{b}) \,= \varepsilon_{\!C}\varepsilon_{\!A} C^\circ{{\,\scriptscriptstyle{[\mathfrak{S}]}\,}}\,A^\circ,\,\;\; \cosh(\mathscr{c}) = \varepsilon_{\!A}\varepsilon_{\!B} A^\circ{{\,\scriptscriptstyle{[\mathfrak{S}]}\,}}\,B^\circ.$}\vspace*{1 mm}\\
By knowing $\mathfrak{C}$, we also know $\varepsilon_{\!A}, \varepsilon_{\!B}, \varepsilon_{\!C}, \cosh(\mathscr{a}), \cosh(\mathscr{b}), \cosh(\mathscr{c})$. With the help of these numbers we can calculate the lengths of the sides of $\Delta_0$, but also the lengths of the sides of the triangles $\Delta_1, \Delta_2,\Delta_3$. For example, the sidelengths of $\Delta_1$ are $\mathscr{a}, \pi i - \mathscr{b}, \pi i - \mathscr{c}.$\vspace*{-2 mm}\\

Denote the measures of the (inner) angles of $\Delta_0$ by \vspace*{1 mm}\\
\centerline{$\alpha := \mu(\angle_{+} BAC),\, \beta := \mu(\angle_{+} CBA),\, \gamma := \mu(\angle_{+} ACB)$.}\vspace*{1 mm}\\

We can calculate $\cosh(\alpha), \cosh(\beta), \cosh(\gamma)$ from $\mathfrak{C}$ and get for example,\vspace*{1.5 mm}\\
\centerline{$\displaystyle
\cosh(\alpha) = - \frac{\mathfrak{d}_{23}}{\sqrt{\mathfrak{d}_{22}}\,\sqrt{\mathfrak{d}_{33}}} = - \frac{\mathfrak{c}_{31}\mathfrak{c}_{12}-\mathfrak{c}_{23}\mathfrak{c}_{11}}{\sqrt{\mathfrak{c}_{11}\mathfrak{c}_{22}-\mathfrak{c}_{12}^2}\sqrt{\mathfrak{c}_{33}\mathfrak{c}_{11}-\mathfrak{c}_{31}^2}}. 
$}\vspace*{3 mm}
%The denominators of the fractions are different from zero because of the assumption in \ref{subsubsec:Note1}.\vspace*{-1 mm}\\
Of course, when we know $\alpha, \beta$ and $\gamma$, we also know the measures of the angles of the other triangles $\Delta_k, k=1,2,3$.
\subsubsection{Some useful trigonometric formulae.\vspace*{0 mm}} 
\begin{equation*}
\begin{split}
\!\cosh^2(\alpha)\; &= \frac{(\varepsilon_{\!\!A}^{\,2}\varepsilon_{\!B}^{\,}\varepsilon_{C}^{\,}(\mathfrak{c}_{31}\mathfrak{c}_{12}-\mathfrak{c}_{23}\mathfrak{c}_{11}))^2}{({\varepsilon_{\!\!A}^{\;2}\varepsilon_{\!B}^{\;2}(\mathfrak{c}_{12}^{\;2}-\mathfrak{c}_{11}^{\,}\mathfrak{c}_{22}^{\,})})\,({\varepsilon_{\!\!A}^{\;2}\varepsilon_{\!C}^{\;2}(\mathfrak{c}_{13}^{\;2}-\mathfrak{c}_{11}^{\,}\mathfrak{c}_{33}^{\,})})}
=\left(\frac{\cosh(\mathscr{b})\,\cosh(\mathscr{c})-\cosh(\mathscr{a})}{\sinh(\mathscr{b}) \sinh(\mathscr{c})}\right)^{\!2}\\
\cosh^2(\mathscr{a})\, &= \cosh^2(\pi i - \mathscr{a}) = \left(\frac{\cosh(\pi i - \beta)\,\cosh(\pi i - \gamma)-\cosh(\pi i - \alpha)}{\sinh(\pi i - \beta) \sinh(\pi i - \gamma)}\right)^2\\
&= \left(\frac{\cosh(\beta)\,\cosh(\gamma)+\cosh(\alpha)}{\sinh(\beta) \sinh(\gamma)}\right)^{\!2}\\
\cosh^2(\mathscr{a})\, &= \left(1 + \frac{\cosh(\beta)\,\cosh(\gamma)+\cosh(\alpha)-\sinh(\beta) \sinh(\gamma)}{\sinh(\beta) \sinh(\gamma)}\right)^{\!2}\\
 &= \left(1 + \frac{\cosh(\alpha)+\cosh(\beta+\gamma)}{\sinh(\beta) \sinh(\gamma)}\right)^{\!2}\\
 &= \left(1 + \frac{2\cosh(\sigma)\cosh(\sigma{-}\alpha)}{\sinh(\beta) \sinh(\gamma)}\right)^{\!2}\; \textrm{with }\sigma = \frac 12 (\alpha+\beta+\gamma)\\
\frac{\sinh^2(\alpha)}{\sinh^2(\mathscr{a})} \;&= \frac{\cosh^2(\alpha)-1}{\sinh^2(\mathscr{a})} = \frac{1}{\sinh^2(\mathscr{a})}\left(\left(\frac{\cosh(\mathscr{b})\,\cosh(\mathscr{c})-\cosh(\mathscr{a})}{\sinh(\mathscr{b}) \sinh(\mathscr{c})}\right)^2 - 1\right)\\
&=\frac{-2 \cosh(\mathscr{a})\cosh(\mathscr{b})\cosh(\mathscr{c})+\cosh^2(\mathscr{a})+\cosh^2(\mathscr{b})+\cosh^2(\mathscr{c})-1}{\sinh^2(\mathscr{a})\sinh^2(\mathscr{b})\sinh^2(\mathscr{c})} \\
&= \frac{\sinh^2(\beta)}{\sinh^2(\mathscr{b})} = \frac{\sinh^2(\gamma)}{\sinh^2(\mathscr{c})}\\
\sinh^2(\alpha)\;\; &= \left(\frac{\cosh(\mathscr{b})\,\cosh(\mathscr{c})-\cosh(\mathscr{a})}{\sinh(\mathscr{b}) \sinh(\mathscr{c})} -1\right)\left(\frac{\cosh(\mathscr{b})\,\cosh(\mathscr{c})-\cosh(\mathscr{a})}{\sinh(\mathscr{b}) \sinh(\mathscr{c})} +1\right)\\
&=\frac{(\cosh(\mathscr{b}+\mathscr{c})-\cosh(\mathscr{a}))(\cosh(\mathscr{b}-\mathscr{c})-\cosh(\mathscr{a}))}{\sinh^2(\mathscr{b}) \sinh^2(\mathscr{c})}\\
\end{split}
\end{equation*}
\begin{equation*}
%\begin{split}
=-\frac{4\sinh(s)\sinh(s-\mathscr{a})\sinh(s-\mathscr{b})\sinh(s-\mathscr{c})}{\sinh^2(\mathscr{b}) \sinh^2(\mathscr{c})}\;\, \textrm{with}\;s=\frac{1}{2}(\mathscr{a}+\mathscr{b}+\mathscr{c})\\
%\end{split}
\end{equation*}
\vspace*{0.2 mm}
\subsubsection{The distance between two points which are given by barycentric coordinates.} \label{subsubsec:The length of a segment in barycentric coordinates.}$\;\;\;\;\;\;$ A point $S = [s_1{:}s_2{:}s_3]_\Delta$ is isotropic precisely when $(s_1,s_2,s_3){\,\scriptscriptstyle{[\mathfrak{C}]}\,}(s_1,s_2,s_3) = 0$.\\ 
Given two different anisotropic points $P=[p_1\,{:}\,p_2\,{:}\,p_3]_\Delta$ and $Q=[p_1\,{:}\,p_2\,{:}\,p_3]_\Delta$, then \vspace*{-1 mm}\\
\begin{equation*}
\cosh^2(d(P,Q)) = \frac{((p_1,p_2,p_3){\,\scriptscriptstyle{[\mathfrak{C}]}\,}(q_1,q_2,q_3))^2}{((p_1,p_2,p_3){\,\scriptscriptstyle{[\mathfrak{C}]}\,}(p_1,p_2,p_3))\,((q_1,q_2,q_3){\,\scriptscriptstyle{[\mathfrak{C}]}\,}(q_1,q_2,q_3))}
\end{equation*}
\begin{equation*}
\varepsilon_{\!P}^{\;2} = \frac{(p_1,p_2,p_3){\,\scriptscriptstyle{[\mathfrak{C}]}\,}(p_1,p_2,p_3)}{|(p_1,p_2,p_3){\,\scriptscriptstyle{[\mathfrak{C}]}\,}(p_1,p_2,p_3)|} \;,\; \varepsilon_{\!Q}^{\;2} = \frac{(q_1,q_2,q_3){\,\scriptscriptstyle{[\mathfrak{C}]}\,}(q_1,q_2,q_3)}{|(q_1,q_2,q_3){\,\scriptscriptstyle{[\mathfrak{C}]}\,}(q_1,q_2,q_3)|} .\vspace*{2 mm}
\end{equation*}
The three numbers $\cosh^2(d(P,Q)), \varepsilon_{\!P}^{\;2}, \varepsilon_{\!Q}^{\;2}$ determine $d(P,Q)$.\vspace*{2 mm}

\subsubsection{The dual of $\Delta$.}\label{subsubsec:The dual triple}
Put $A':= a^\delta, B':= b^\delta, C':= c^\delta$. The triple $\Delta' = A'B'C'$ is called the dual of $\Delta$. The representation of $A'$ by barycentric coordinates with respect to $\Delta$ is \\$[(B^\circ{\times}\,C^\circ){\,\scriptscriptstyle{[\mathfrak{S}^{-1}]}\,}(B^\circ{\times}\,C^\circ):(B^\circ{\times}\,C^\circ){\,\scriptscriptstyle{[\mathfrak{S}^{-1}]}\,}(C^\circ{\times}\,A^\circ):(B^\circ{\times}\,C^\circ){\,\scriptscriptstyle{[\mathfrak{S}^{-1}]}\,}(A^\circ{\times}\,B^\circ)]_\Delta$, with $\mathfrak{S} = \text{diag}(\rho,1,1)$. It can be easily checked (see also \ref{subsubsec:conics}) that \\
\centerline{$A' = [\mathfrak{d}_{11}{:}\mathfrak{d}_{12}{:}\mathfrak{d}_{13}]_\Delta = [\mathfrak{c}_{22}\mathfrak{c}_{33}-\mathfrak{c}_{23}^2 :\mathfrak{c}_{23}\mathfrak{c}_{31}-\mathfrak{c}_{12}\mathfrak{c}_{33} :\mathfrak{c}_{12}\mathfrak{c}_{23}-\mathfrak{c}_{31}\mathfrak{c}_{22}]_\Delta\,.$}\vspace*{1 mm}\\
In the same way the barycentric coordinates of $B'$ and $C'$ can be calculated, getting\vspace*{1 mm}\\
\centerline{$B' = [\mathfrak{d}_{21}{:}\mathfrak{d}_{22}{:}\mathfrak{d}_{23}]_\Delta, 
C' = [\mathfrak{d}_{31}{:}\mathfrak{d}_{32}{:}\mathfrak{d}_{33}]_\Delta, $}\\
and we conclude:\vspace*{0.5 mm}\\
\centerline{$A = [\mathfrak{c}_{11}{:}\mathfrak{c}_{12}{:}\mathfrak{c}_{13}]_{\Delta'}\,,\; B = [\mathfrak{c}_{21}{:}\mathfrak{c}_{22}{:}\mathfrak{c}_{23}]_{\Delta'}\,,\; C = [\mathfrak{c}_{31}{:}\mathfrak{c}_{32}{:}\mathfrak{c}_{33}]_{\Delta'}$\,.}\vspace*{2 mm}

%\noindent\textit{Remark}: If we take $A^\circ, B^\circ, C^\circ$ as a covariant basis of $(V,{\,\scriptscriptstyle{[\mathfrak{S}]}\,})$, then the corresponding contravariant basis consists of vectors $\boldsymbol{a'},\boldsymbol{b'},\boldsymbol{c'}$ with $A'=\Pi(\boldsymbol{a'}), B'=\Pi(\boldsymbol{b'}), C'=\Pi(\boldsymbol{c'})$.
\subsubsection{The dual triangles.}\label{subsubsec:The dual triangles.} 
Just like $\Delta$, the triple $\Delta'$ determines four triangles $\Delta'_{\,k}, k = 0,1,2,3$. But, in general, triangle $\Delta'_{\,k}$ is not the dual of $\Delta_k$. In fact, the dual of $\Delta_k$ is   \centerline{$\Delta^{k} := \mathcal{P} - \bigcup_{P\, \text{inner point of}\, \Delta_k} P^\delta $.}\vspace*{1 mm}\\
We put $a':= B'{\vee}C',\, b':= C'{\vee}A',\, c':= A'{\vee}B'$ and denote by $\mathscr{a}', \mathscr{b}', \mathscr{c}'$ the lengths of the sides and by $\alpha', \beta', \gamma'$ the measures of the angles of triangle $\Delta^0$. Then (see \ref{subsubsec:The length of a segment and the measure of an angle.}),\vspace*{-3 mm}\\
\begin{equation*}
\begin{split}
\mathscr{a}' &= \pi i - \alpha,\; \mathscr{b}' = \pi i - \beta,\, \mathscr{c}' = \pi i - \gamma\,,\\
\alpha' &= \pi i - \mathscr{a},\,\beta' = \pi i - \mathscr{b},\, \gamma' = \pi i - \mathscr{c}\,.
\end{split}
\end{equation*}
\subsubsection{Reflections.} 
For each anisotropic point $M = [m_1{:}m_2{:}m_3]_\Delta$ we define a mapping $\sigma_M{:\,}\mathcal{P}-\mathcal{C}_{\!abs}\rightarrow\mathcal{P}$ as follows: The image of an anisotropic point $P = [p_1{:}p_2{:}p_3]_\Delta$ under $\sigma_M$ is the point $Q = [q_1{:}q_2{:}q_3]_\Delta$ with
\begin{equation*}
%\begin{split}
(q_1,q_2,q_3) = (p_1,p_2,p_3) - 2\frac{(m_1,m_2,m_3){\,\scriptscriptstyle{[\mathfrak{C}]}\,}(p_1,p_2,p_3)}{(m_1,m_2,m_3){\,\scriptscriptstyle{[\mathfrak{C}]}\,}(m_1,m_2,m_3)} (m_1,m_2,m_3).
%\end{split}
\end{equation*}
The points $P, M, Q$ are obviously collinear, and it can be easily checked that $P$ and $Q$ have the same distance from $M$. In particular, $Q$ is also an anisotropic point, lying together with $P$ in the same connected component of $\mathcal{P}-\mathcal{C}_{\!abs}$. It can also be verified that $\sigma_M(Q) = P$, and that $P = Q$ only if either $P = M$ or $d(P,M) = \frac{\pi}{2} i.$ \\
We extend $\sigma_M$ to a mapping which is continuous on $\mathcal{P}$:
If $P$ is the only isotropic point on the (isotropic) line $P\vee M$, then $\sigma_M(P) = P$; if the line $P\vee M$ contains still another isotropic point $Q$ besides $P$, then $\sigma_M(P) = Q$.\\ This mapping $\sigma_M{:\,}\mathcal{P}\rightarrow\mathcal{P}$ is an involution.
We call it the reflection in the point $M$ and call $Q$ the mirror image of $P$ under $\sigma_M$.
It is not difficult to show that this reflection preserves the distance between points and the (angle) distance between lines.\vspace*{1 mm}\\
\noindent \textit{Remark:} A reflection $\sigma_M$ in an anisotropic point $M$ can also be interpreted as a reflection $\sigma_\ell$ in the line $\ell = M^\delta.$ 
\vspace*{-2 mm}\\

If a nonempty subset $\mathcal{S}$ of $\mathcal{P}$ is invariant under the reflection in an anisotropic point $M$, then $M$ is called a \textit{symmetry point} of $\mathcal{S}$; the line $M^\delta$ is a \textit{symmetry axis} of $\mathcal{S}$. A special case: If $P$ and $Q$ are two different anisotropic points, the points $P{+\,}Q$ and $P{-\,}Q$ are symmetry points of $\{P,Q\}$ precisely when they are (proper) midpoints. 

\subsubsection{The pedals and antipedals of a point.} We calculate the coordinates of the pedals of a point $P=[p_1:p_2:p_3]_\Delta$ on the sidelines of $\Delta$:
\begin{equation*}
\begin{split}
&A_{[P]} := \text{ped}(P,a) = (P\vee C') \wedge a = [0:\mathfrak{d}_{11}p_2-\mathfrak{d}_{12}p_1:\mathfrak{d}_{11}p_3-\mathfrak{d}_{31}p_1]_\Delta,\\
&B_{[P]} = [\mathfrak{d}_{22}p_1-\mathfrak{d}_{12}p_2:0:\mathfrak{d}_{22}p_3-\mathfrak{d}_{23}p_2]_\Delta,\\
&C_{[P]} = [\mathfrak{d}_{33}p_1-\mathfrak{d}_{31}p_3:\mathfrak{d}_{33}p_2-\mathfrak{d}_{23}p_3:0]_\Delta.
\end{split}
\end{equation*}
The \textit{antipedal points} of $P$ are\vspace*{-2 mm}\\
\begin{equation*}
\begin{split}
&A^{[P]} := \text{perp}(B\vee P,B)\wedge \text{perp}(C\vee P,C) = [-1:\frac{\mathfrak{d}_{22}p_1-\mathfrak{d}_{12}p_2}{\mathfrak{d}_{11}p_2-\mathfrak{d}_{12}p_1}:\frac{\mathfrak{d}_{33}p_1-\mathfrak{d}_{31}p_3}{\mathfrak{d}_{11}p_3-\mathfrak{d}_{31}p_1}]_\Delta,\\
&B^{[P]} = [\frac{\mathfrak{d}_{11}p_2-\mathfrak{d}_{12}p_1}{\mathfrak{d}_{22}p_1-\mathfrak{d}_{12}p_2}:-1:\frac{\mathfrak{d}_{33}p_2-\mathfrak{d}_{23}p_3}{\mathfrak{d}_{22}p_3-\mathfrak{d}_{23}p_2}]_\Delta,\\
&C^{[P]} = [\frac{\mathfrak{d}_{11}p_3-\mathfrak{d}_{31}p_1}{\mathfrak{d}_{33}p_1-\mathfrak{d}_{31}p_3}:\frac{\mathfrak{d}_{22}p_3-\mathfrak{d}_{23}p_2}{\mathfrak{d}_{33}p_2-\mathfrak{d}_{23}p_3}:-1]_\Delta,\\
\end{split}
\end{equation*}

\subsubsection{The cevian and the anticevian triple of a point.}
If $P = [p_1{:}p_2{:}p_3]_\Delta$ is a point different from $A,B,C$, then the lines $P\vee A,$ $P\vee B,$ $P\vee C$ are called the \textit{cevian lines} of $P$ with respect to $\Delta$. The cevian lines meet the sidelines $a, b, c$ at the points $A_{P} =[0{:}p_2{:}p_3]_\Delta, B_{P} = [p_1{:}0{:}p_3]_\Delta, C_{P} = [p_1{:}p_2{:}0]_\Delta$, respectively. These points are called the \textit{cevian points} or the \textit{traces}, and the triple $A_{P}B_{P}C_{P}$ is called the \textit{cevian triple} of $P$ with respect to $\Delta$. The \textit{anticevian triple} $A^{P}B^{P}C^{P}$ consists of the \textit{harmonic associates} $A^{P} = [-p_1{:}p_2{:}p_3]_\Delta$, $B^{P} = [p_1{:}-p_2{:}p_3]_\Delta$, $C^{P} = [p_1{:}p_2{:}-p_3]_\Delta$ of $P$.
\subsubsection{Tripolar and tripole}
Given a point $P = [p_1{:}p_2{:}p_3]_\Delta$, then the point $[0{:}-p_2{:}p_3]_\Delta$ is the harmonic conjugate of $A_P$ with respect to $\{B,C\}$. Correspondingly, the harmonic conjugates of the traces of $P$ on the other sidelines are $[-p_1{:}0{:}p_3]_\Delta$  and $[p_1{:}-p_2{:}0]_\Delta$. These three harmonic conjugates are collinear; the equation of the line $l$ is $p_2p_3 x_1 + p_3p_1 x_2 + p_1p_2x_3 = 0.$ This line is called the \textit{tripolar line} or the \textit{tripolar} of $P$ and we denote it by $P^\tau$. $P$ is the \textit{tripole} of $l$ and we write $P = l^\tau$.
%Theorem: 

%G\"ulicher's theorem: $A_P B_P C_P$ is the cevian triple of a point $P = [p_1{:}p_2{:}p_3]_\Delta$  with respect to $\Delta$. Now, let $Q_1, Q_2, Q_3$ be the cevian triple of a point $Q = [q1{:}q2{:}q3]_\Delta$ with respect to $A_P B_P C_P$.  Then the lines $Q_1 \vee A, Q_2 \vee B, Q_3 \vee C$ are concurrent at a point $R$. \\
%We calculate the coordinates of $R$: $R = [p_1/(q_1p_2p_3-q_2p_3p_1-q_3p_1p_2):\cdots:\cdots]_\Delta.$
\subsubsection{Isoconjugation} Let $P = [p_1{:}p_2{:}p_3]_\Delta$ be a point not on a sideline of $\Delta$ and let the point $Q = [q_1{:}q_2{:}q_3]_\Delta$ be not a vertex of $\Delta$, then the point $R =[p_1q_2q_3{:}p_2q_3q_1{:}p_3q_1q_2]_\Delta$ is called 
%the \textit{isoconjugate} of $Q$ with respect to $P$ or shorter 
the $P$-isoconjugate of $Q$ with respect to $\Delta$. \\
Special cases: If $P = G = [1{:}1{:}1]_\Delta$ (the centroid of $\Delta_0$), $P$-isoconjugation is called \textit{isotomic} conjugation, and for $P = K = [\mathfrak{d}_{11}{:}\mathfrak{d}_{22}{:}\mathfrak{d}_{33}]_\Delta$ (the  
symmedian of $\Delta_0$) it is called  \textit{isogonal} conjugation. 
The isotomic (resp. isogonal) conjugate of a point $P$ with respect to $\Delta$ agrees with the isogonal (resp. isotomic) conjugate of $P$ with respect to $\Delta'$.

\subsubsection{The area of a triangle.} For a triangle $\mathcal{T}$ \footnote{$^)$ Implicitly it is always assumed that the vertices and the sidelines of a triangle are anisotropic.}$^)$ we define its \textit{area} (also called its \textit{excess}) by $\textrm{area}(\mathcal{T}) = \alpha + \beta + \gamma - \pi i$, where $\alpha, \beta, \gamma$ are the measures of the inner angles of $\mathcal{T}$. This function is additive: If we dissect $\mathcal{T}$ in finitely many triangles, then the sum of the areas of these parts equals the area of $\mathcal{T}$. Adding up all the areas of the triangles $\Delta_k, k = 1,2,3,4$, we get $2\pi i$ for the area of the whole plane $\mathcal{P}$.

\subsubsection{Conics.}\label{subsubsec:conics} Let $\mathfrak{M} = (\mathfrak{m}_{ij})_{i,j=1,2,3}$ be a symmetric matrix, then we denote by $\mathcal{Q}(\mathfrak{M})$ the quadratic form \vspace*{1 mm}\\
\centerline{$q(x_1,x_2,x_3) = \mathfrak{m}_{11}x_1^2 + \mathfrak{m}_{22}x_2^2 + \mathfrak{m}_{33}x_3^2 + 2\mathfrak{m}_{23}x_2 x_3 + 2\mathfrak{m}_{31}x_3 x_1 + 2\mathfrak{m}_{12}x_1 x_2.$}\vspace*{1 mm}\\
If $\mathfrak{M}$ is indefinite, the set of points $[x_1{:}x_2{:}x_3]_\Delta$ which satisfy 
the quadratic equation $q(x_1,x_2,x_3)=0$
is a \textit{real conic} which we denote by $\mathcal{C}(\mathfrak{M})$. This conic is called \textit{nondegenerate} if det$(\mathfrak{M}) \ne 0$. \vspace*{-2 mm}\\

The polar of a point $P = [p_1{:}p_2{:}p_3]_\Delta $ with respect to $\mathcal{Q}(\mathfrak{M})$ is the line with the equation $(p_1,p_2,p_3){\,\scriptscriptstyle{[\mathfrak{M}]\,}}(x_1,x_2,x_3) = 0$, and	
	%\[ (p_1,p_2,p_3)\,\mathfrak{M}\!\begin{pmatrix} x_1 \\ x_2\\ x_3 \end{pmatrix} = 0.\]
the pole of the line $\ell\!\!: \ell_{\!1} x_1 + \ell_{\!2} x_2 + \ell_{\!3} x_3 = 0$ with respect to $\mathcal{Q}(\mathfrak{M})$ is the point $P = [p_1{:}p_2{:}p_3]_\Delta$ with
	$(p_1,p_2,p_3) = (\ell_{\!1},\ell_{\!2}, \ell_{\!3})\mathfrak{M}^{\#}$, 
where\vspace*{0.5mm} \\ 
\centerline{$\mathfrak{M}^\# = \begin{pmatrix} \mathfrak{m}_{22}\mathfrak{m}_{33}-\mathfrak{m}_{23}^2&\mathfrak{m}_{23}\mathfrak{m}_{31}-\mathfrak{m}_{12}\mathfrak{m}_{33}&{\mathfrak{m}_{12}\mathfrak{m}_{23}-\mathfrak{m}_{31}\mathfrak{m}_{22}}\\ 
\mathfrak{m}_{23}\mathfrak{m}_{31}-\mathfrak{m}_{12}\mathfrak{m}_{33}&\mathfrak{m}_{33}\mathfrak{m}_{11}-\mathfrak{m}_{31}^2&\mathfrak{m}_{31}\mathfrak{m}_{12}-\mathfrak{m}_{23}\mathfrak{m}_{11} \\
\mathfrak{m}_{12} \mathfrak{m}_{23}-\mathfrak{m}_{31}\mathfrak{m}_{22}&\mathfrak{m}_{31}\mathfrak{m}_{12}-\mathfrak{m}_{23}\mathfrak{m}_{11}&\mathfrak{m}_{11}\mathfrak{m}_{22}-\mathfrak{m}_{12}^2 \end{pmatrix}$}\vspace*{1.5mm} \\
\noindent is the adjoint of $\mathfrak{M}.$ \vspace*{1.5 mm}

\noindent A special case: In the hyperbolic plane ($\rho<0$),  the conic $\mathcal{C}(\mathfrak{C})$ is the \textit{absolute conic} which consists of all the isotropic points. In the elliptic plane, $\mathfrak{C}$ is definite and no real isotropic points exist.\vspace*{-1mm}\\

A point $P$ is a symmetry point of a nondegenerate real conic $\mathcal{C}(\mathfrak{M}), \mathfrak{M}\ne\mathfrak{C}, $ precisely when the polar of $P$ equals $P^\delta$.\vspace*{1 mm}\\
\noindent \textit{Proof}: Take any line through $P$ that meets $\mathcal{C}(\mathfrak{M})$ in two points $Q$ and $R$. This line meets the polar of $P$ at the harmonic conjugate  $S$ of $P$ with respect to $Q$ and $R$. Precisely when $P$ is a proper midpoint of $\{Q, R\}$, the point $S$ is also a proper midpoint of $\{Q, R\}$ and lies on $P^\delta.\;\;\Box$ \\

If $\mathfrak{M}$ is a diagonal matrix, then the polar lines of $A, B, C$ with respect to $\mathcal{Q}(\mathfrak{M})$ are the lines $B\vee C$, $C\vee A$, $A\vee B$, respectively. If $\mathfrak{M}$ is not diagonal,
then the poles of $B\vee C$, $C\vee A$, $A\vee B$ with respect to $\mathcal{Q}(\mathfrak{M})$ form a triple, perspective to $\Delta$ at perspector \vspace*{-1 mm}
	\[[\frac 1 {\mathfrak{m}_{11} \mathfrak{m}_{23}- \mathfrak{m}_{31} \mathfrak{m}_{12}}: \frac 1 {\mathfrak{m}_{22} \mathfrak{m}_{31}- \mathfrak{m}_{12} \mathfrak{m}_{23}}:\frac 1 {\mathfrak{m}_{33} \mathfrak{m}_{12}- \mathfrak{m}_{23} \mathfrak{m}_{31}}]_\Delta.\]
This perspector is called the \textit{perspector of} $\mathcal{Q}(\mathfrak{M})$ with respect to $\Delta$ and, if $\mathcal{C}(\mathfrak{M})$ is a real conic, the \textit{perspector of} $\mathcal{C}(\mathfrak{M})$ with respect to $\Delta$.\vspace*{-2 mm}\\

Let $\mathcal{C}(\mathfrak{M})$ be a nondegenerate real conic. The set $\{p^\delta |\, p\;\text{is a tangent of}\; \mathcal{C}(\mathfrak{M})\}$ is the conic $\mathcal{C}(\mathfrak{C M}^{\sharp}\mathfrak{C})$; we call it the \textit{dual} of $\mathcal{C}(\mathfrak{M})$. Both conics, $\mathcal{C}(\mathfrak{M})$ and its dual, share the same symmetry points and symmetry axes. \vspace*{-2 mm}\\

There are other conics named "dual of $\mathcal{C}(\mathfrak{M})$". One, for example, is $\mathcal{C}(\mathfrak{M^\sharp})$. To avoid a name collision, we call this one the \textit{adjoint conic}. 
\subsubsection{Circles} We assume that $\mathcal{C}(\mathfrak{M}), \mathfrak{M}\ne\mathfrak{C},$ is a nondegenerate real conic. Then $\mathcal{C}(\mathfrak{M})$ is a \textit{circle} if there exists a line $\ell$ consisting of symmetry points. In this case, the point $\ell^\delta$ is also a symmetry point of $\mathcal{C}(\mathfrak{M})$ and is called the \textit{center} of the circle.\\
\textit{Remarks}:\\
$\bullet$ A circle may have (up to two) isotropic points.\\
$\bullet$ In the elliptic plane, the center of a circle always lies inside the circle. In the hyperbolic plane, the center of a circle lies inside the circle precisely when this center also lies inside the absolute conic.\vspace*{-2 mm}\\

If $P$ and $Q$ are two distinct anisotropic points on a circle with center $M$, then $d(P,M) = d(Q,M)$, because: 
Let $R$ be the intersection of the lines $P\vee Q$ and $M^\delta$, then the triple $PMQ$ is mapped onto the triple $QMP$ by the (distance preserving) reflection $\sigma_{R}$. \vspace*{0.5 mm}\\
So we can define the \textit{radius} of a circle as the common distance between its center and any of its anisotropic points. \vspace*{-2 mm}\\

Given two anisotropic points $M$ and $P$, then there exists a unique nondegenerate circle $\mathscr{C}(M,P)$ with center $M$ through the point $P$ precisely when $0 \prec d(M,P) \prec \frac{\pi}{2} i$. \\
Besides the nondegenerate circles, one can regard the following degenerate conics as circles:  It is common to look at a  double line $l$ (a line with multiplicity 2) as a circle with center $l^\delta$ and radius $\frac{\pi}{2} i$. A double point can be considered as a circle with radius $0$ around this point. And also a degenerate conic consisting of two different isotropic lines can be seen as a circle; its center is the meet of these lines and its radius is $0$. 

\section{Special triangle centers, central lines, conics and cubics}\vspace*{2 mm}

\subsection{Triangle centers and central triangles based on orthogonality}
\subsubsection{The common orthocenter of $\Delta$ and $\Delta'$.} The perspector of $\mathcal{Q}(\mathfrak{C})$ is the common \textit{orthocenter} $H$ of $\Delta$ and $\Delta'$,\vspace*{1 mm}\\
\centerline{$H=[\displaystyle\frac{1}{\mathfrak{c}_{31}\mathfrak{c}_{12}-\mathfrak{c}_{23}\mathfrak{c}_{11}}:\frac{1}{\mathfrak{c}_{12}\mathfrak{c}_{23}-\mathfrak{c}_{31}\mathfrak{c}_{22}}:\frac{1}{\mathfrak{c}_{23}\mathfrak{c}_{31}-\mathfrak{c}_{12}\mathfrak{c}_{33}}]_\Delta=[{\mathfrak{d}_{31}}{\mathfrak{d}_{12}}{\,:\,}{\mathfrak{d}_{12}}{\mathfrak{d}_{23}}{\,:\,}{\mathfrak{d}_{23}}{\mathfrak{d}_{31}}]_\Delta$.}\\

Since $A{\vee}H$ is orthogonal to $B{\vee}C$, $B{\vee}H$ orthogonal to $C{\vee}A$ and $C{\vee}H$ orthogonal to $A{\vee}B$, the points $A, B, C, H$ together with the
lines $a, b, c, A{\vee}H,$ $B{\vee}H, C{\vee}H$ form an \textit{orthocentric system}; each of the points is the orthocenter of the other three. The same applies to $A', B', C', H$, together with the lines through any two of them. If we combine these two systems and add the points $a{\wedge}a', b{\wedge}b', c{\wedge}c'$ and the line $H^\delta$, we get a system which - if $\Delta_0$ has neither a right angle nor a right side - consists of 10 points and 10 lines such that the dual of each of its points is one of its lines and the dual of each of its lines is one of its points, and each point is incident with 3 lines and each line incident with 3 points.

\subsubsection{Another triple perspective to $\Delta$ and $\Delta'$ at perspector $H$.}
Define $A_b := b\wedge a'$, $B_a := a\wedge b'$, and the points $A_c, B_c, C_a, C_b$ accordingly. 
We calculate the coordinates of these points: $A_b = [\mathfrak{c}_{31}{\,:\,}0{\,:}-\mathfrak{c}_{11}]_\Delta, B_a = [0{\,:\,}\mathfrak{c}_{23}{\,:}-\mathfrak{c}_{22}]_\Delta$, etc. 
The point-triple $\tilde{A}\tilde{B}\tilde{C}$, \vspace*{-3 mm}\\
\begin{equation*}
\begin{split}
\tilde{A}&:= (A_b\vee B_a)\wedge(A_c\vee C_a) = [\frac{\mathfrak{d}_{11}}{\mathfrak{c}_{11}}\frac{\mathfrak{c}_{31} \mathfrak{c}_{12}}{\mathfrak{c}_{23}}:\mathfrak{d}_{12}:\mathfrak{d}_{13}]_\Delta,\\
\tilde{B}&:= (B_c\vee C_b)\wedge(B_a\vee A_b) = [\mathfrak{d}_{21}:\frac{\mathfrak{d}_{22}}{\mathfrak{c}_{22}}\frac{\mathfrak{c}_{12} \mathfrak{c}_{23}}{\mathfrak{c}_{31}}:\mathfrak{d}_{23}]_\Delta,\\
\tilde{C}&:= (C_a\vee A_c)\wedge(C_b\vee B_c) = [\mathfrak{d}_{31}:\mathfrak{d}_{32}:\frac{\mathfrak{d}_{33}}{\mathfrak{c}_{33}}\frac{\mathfrak{c}_{23} \mathfrak{c}_{31}}{\mathfrak{c}_{12}}]_\Delta,\vspace*{-1 mm}\\
\end{split}
\end{equation*}
is perspective to $\Delta$ and $\Delta'$ at perspector $H$, cf. Vigara \cite{Vi} section 8.2. \vspace*{1 mm}\\

\begin{figure}[!htbp]
\includegraphics[height=9cm]{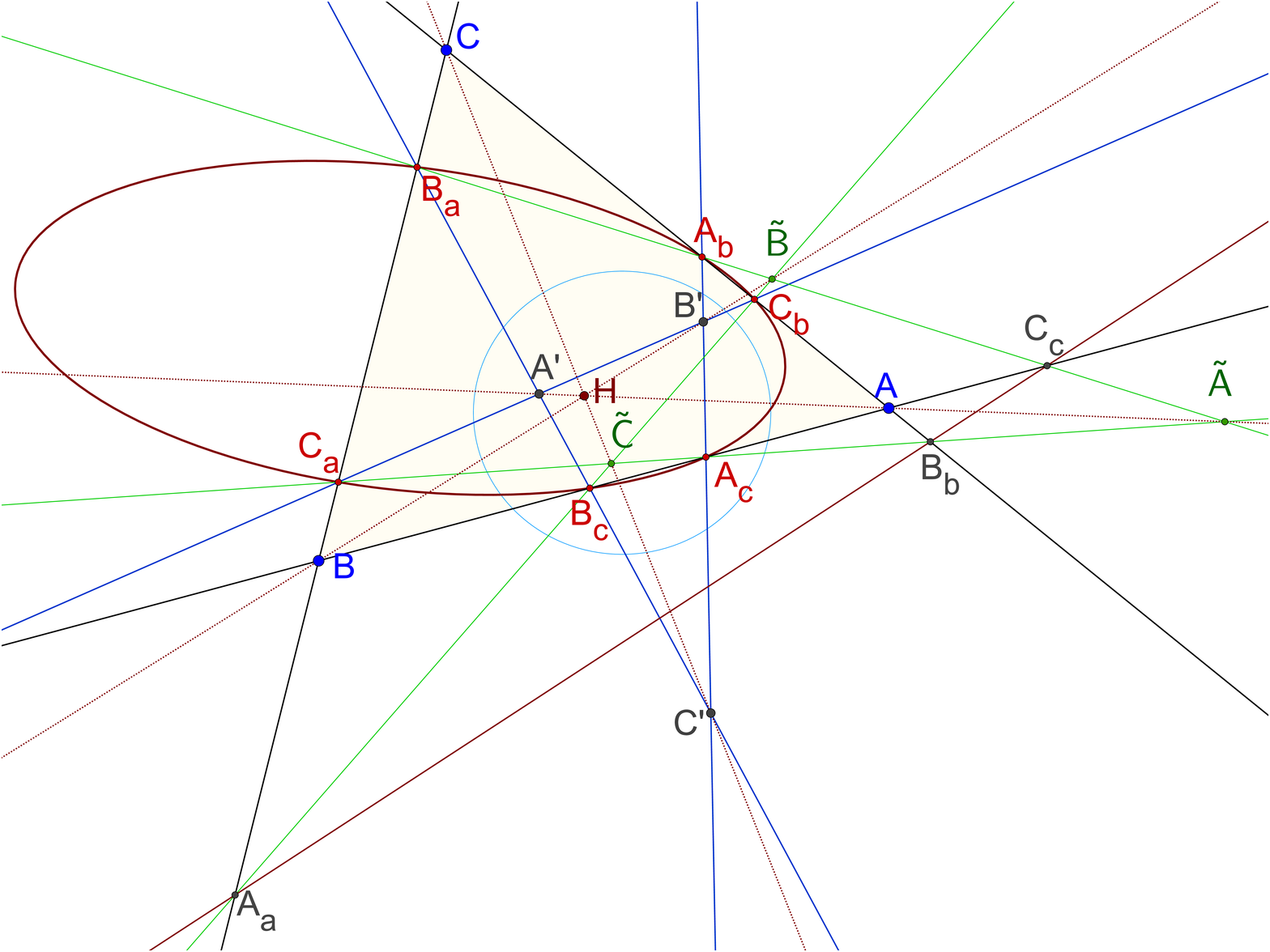}
\caption{The triples $\Delta, \Delta'$ and $\tilde{A}\tilde{B}\tilde{C}$. The light blue circle is the absolute conic in the hyperbolic plane.
\newline
All figures were created with the software program GeoGebra \cite{GG}.}
\end{figure}

The points $B_a, C_a, C_b, A_b, A_c, B_c$ lie on the conic\vspace*{-2 mm}
\begin{equation*}
\{[x_1{:}x_2{:}  x_3]_\Delta\, |\, \sum\limits_{j=1}^3 \mathfrak{c}_{j,j}^{\;2}x_j^{\;2} +\,(\mathfrak{c}_{j+1,j+2} + \frac{\mathfrak{c}_{j+1,j+1} \mathfrak{c}_{j+2,j+2}}{\mathfrak{c}_{j+1,j+2}})x_{j+1}x_{j+2}  = 0 \}\;\;\;(\text{indices mod}\,3)\vspace*{-2 mm}\end{equation*}

with perspector $[\mathfrak{c}_{23}(1+\mathfrak{c}_{12}^{\;2})(1+\mathfrak{c}_{31}^{\;2})-2\mathfrak{c}_{12}\mathfrak{c}_{31}(1+\mathfrak{c}_{23}^{\;2}):\cdots:\cdots]_\Delta$.\vspace*{1 mm}\\
The three points $A_a := (\tilde{B}\vee\tilde{C})\wedge a, B_b := (\tilde{C}\vee\tilde{A})\wedge b, C_c := (\tilde{A}\vee\tilde{B})\wedge c$ lie on the line \vspace*{-1 mm}
\begin{equation*}
\{[x_1{:}x_2{:}x_3]_\Delta\, |\, \mathfrak{c}_{11} \mathfrak{c}_{23} x_1 + \mathfrak{c}_{22} \mathfrak{c}_{31} x_2 + \mathfrak{c}_{33} \mathfrak{c}_{12} x_3 =0.\}
\end{equation*}
By a routine calculation we can determine the coordinates of the orthocenter of the triple $\tilde{\Delta}=\tilde{A}\tilde{B}\tilde{C}$. The expression for it is rather complicated, so we do not present it here. But it might be worth mentioning that the euclidean limit ($\rho \rightarrow \pm \infty$) of this point is the {de{\,}Longchamps point} $X_{20}$ \footnote{$^)$ The designation of the triangle centers corresponds to ETC \cite{ETC}. }$^)$ of $\Delta$.\\ The dual $\tilde{\Delta}'$ of $\tilde{\Delta}$ is perspective to $\Delta$. The perspector is $O^+ =[\mathfrak{d}_{11}\mathfrak{d}_{23}{:}\mathfrak{d}_{22}\mathfrak{d}_{31}{:}\mathfrak{d}_{33}\mathfrak{d}_{12}]_\Delta. $\\The euclidean limit of $\tilde{\Delta}'$ is the mirror image of $\Delta$ in the circumcenter of $\Delta_0$.\vspace*{1 mm}

The six points $(B{\vee} A')\wedge(A{\vee} C')$, $(C{\vee} A')\wedge(A{\vee} B')$, $(C{\vee} B')\wedge(B{\vee} A')$, $(A{\vee} B')\wedge(B{\vee} A')$, $(A{\vee} C')\wedge(C{\vee} B')$, $(B{\vee} C')\wedge(C{\vee} A')$ lie on a conic. The perspector of this conic with respect to $\Delta$ is the point 
$\displaystyle[\frac{\mathfrak{d}_{11}^{\,}}{\mathfrak{d}_{11} \mathfrak{d}_{23} - \mathfrak{d}_{31} \mathfrak{d}_{12}}:\frac{\mathfrak{d}_{22}}{\mathfrak{d}_{22} \mathfrak{d}_{31} - \mathfrak{d}_{12} \mathfrak{d}_{23_{{\,}^{\;}}}}:\frac{\mathfrak{d}_{33}}{\mathfrak{d}_{33_{\,}^{\,}} \mathfrak{d}_{12} - \mathfrak{d}_{23} \mathfrak{d}_{31}} ]_{\Delta}$,
it is the isogonal 
 conjugate of the {de}$\!$ {Longchamps} point $L$, which will be introduced in the next subsection. The center of the euclidean limit of the conic is the circumcenter of $\Delta_0$.
\subsubsection{The double triangle.} We assume that the orthocenter $H$ of $\Delta$ is not a vertex. Wildberger \cite{W1} showed that the antipedals $A^{[H]}, B^{[H]}, C^{[H]}$ of $H$ are the harmonic associates of the point\vspace*{1 mm}\\ 
\centerline{$G^+ ={[\mathfrak{c}_{23}{:}\mathfrak{c}_{31}{:}\mathfrak{c}_{12}]_\Delta} = [\mathfrak{d}_{11} \mathfrak{d}_{23} - \mathfrak{d}_{12} \mathfrak{d}_{31} : \mathfrak{d}_{22} \mathfrak{d}_{31} - \mathfrak{d}_{12} \mathfrak{d}_{23}  : \mathfrak{d}_{33} \mathfrak{d}_{12}  - \mathfrak{d}_{23} \mathfrak{d}_{31}]_\Delta$.}\vspace*{1 mm}\\
It can be easily checked that $G^+$ is the tripole of the dual line of $H$.
The points $A^{[H]}, B^{[H]}$, $C^{[H]}$ are the vertices of a triangle which contains $G^+$. In \cite{W1}, this triangle is called the \textit{double triangle} of $\Delta$, and $G^+$ is called the \textit{double point}.\\ Furthermore, Wildberger proved that the points $A,B,C$ are proper midpoints of $\{\!B^{[H]}{,}C^{[H]}\!\}$, $\{C^{[H]},A^{[H]}\}$, \;$\{A^{[H]},B^{[H]}\}$, respectively. Thus, the point $H$ is the center of a circle through the points $A^{[H]}, B^{[H]}, C^{[H]}$. The perspector of this circle with respect to $\Delta$ is $H$, the perspector with respect to the triple $A^{[H]}B^{[H]}C^{[H]}$ (the Lemoine point of the double triangle) is the point $[\mathfrak{d}_{23}^{\; } \mathfrak{c}_{23}^{\;2}{:}\mathfrak{d}_{31}^{\; } \mathfrak{c}_{31}^{\;2}{:}\mathfrak{d}_{12}^{\; } \mathfrak{c}_{12}^{\;2}]_\Delta$.\vspace*{1 mm}\vspace*{1 mm} \\
The point $O^+$ is identical with a \textit{pseudo-circumcenter}, introduced in \cite{Vi} as the meet of \textit{perpendicular pseudo-bisectors} $A_{G^+}\!\vee A', B_{G^+}\!\vee B', C_{G^+}\!\vee C'$. And it is also identical with the \textit{basecenter} introduced in \cite{W1}.\vspace*{1 mm}\\
The triple $\Delta^{[H]} = A^{[H]}B^{[H]}C^{[H]}$ is perspective to $\Delta'$ at perspector \vspace*{1 mm}\\
\centerline{$L = [\mathfrak{d}_{11} \mathfrak{d}_{23}\! + \mathfrak{d}_{31} \mathfrak{d}_{12} : \mathfrak{d}_{22} \mathfrak{d}_{31}\! + \mathfrak{d}_{12} \mathfrak{d}_{23} : \mathfrak{d}_{33} \mathfrak{d}_{12}\!  + \mathfrak{d}_{23} \mathfrak{d}_{31}]_\Delta$.}\vspace*{1 mm}\\
This point $L$ we will call the \textit{de}$\!$ \textit{Longchamps} \textit{point} of $\Delta$.\\

\subsubsection{The orthic triangle.} We assume again that the orthocenter $H$ of $\Delta$ is not a vertex. The \textit{orthic triangle} has the vertices $A_{[H]}, B_{[H]}, C_{[H]}$, it contains $H$, and it is the dual of the double triangle of $\Delta'$. The point $H$ is the incenter and the points $A, B, C$ are the excenters of the orthic triangle. The pedals of $H$ on the sidelines of $\Delta_{[H]}$ are the traces with respect to $\Delta_{[H]}$ of the point \vspace*{-1 mm}
\[ 
[\mathfrak{d}_{31}^{\;}\mathfrak{d}_{12}^{\;} (\mathfrak{d}_{31}^{\;2}\mathfrak{d}_{22}^{\;}-2 \mathfrak{d}_{31}^{\;}\mathfrak{d}_{12}^{\;}\mathfrak{d}_{23}^{\;}+\mathfrak{d}_{12}^{\;2}\mathfrak{d}_{33}^{\;})(2\mathfrak{d}_{11}^{\;}\mathfrak{d}_{23}^{\;2}+\mathfrak{d}_{31}^{\;2}\mathfrak{d}_{22}^{\;}-4 \mathfrak{d}_{31}^{\;}\mathfrak{d}_{12}^{\;}\mathfrak{d}_{23}^{\;}+\mathfrak{d}_{12}^{\;2}\mathfrak{d}_{33}^{\;}):\cdots:\cdots]_\Delta. \vspace*{0 mm}
\]
This point has euclidean limit $X_{52}$.\vspace*{-2 mm}\\
%                  2      2                              2      2                              2 
%d12 d31 (2 d11 d23  + d12  d33 - 4 d12 d23 d31 + d22 d31 ) (d12  d33 - 2 d12 d23 d31 + d22 d31 )
\subsubsection{The antipedal triple of $O^+$} The points $A^{[O^+]} = [-\mathfrak{d}_{11}/\mathfrak{c}_{23}, \mathfrak{d}_{22}/\mathfrak{c}_{31}, \mathfrak{d}_{33}/\mathfrak{c}_{12}]_\Delta$, $B^{[O^+]}$, $C^{[O^+]}$ form a triple $\Delta^{[O^+]}$, which is perspective to $\Delta$ at the isogonal conjugate of $G^+$ and perspective to the orthic triple $\Delta_{[H]}$ at 
\[[\frac{\mathfrak{d}_{23}^{\;}(\mathfrak{d}_{11}^{\;}(\mathfrak{d}_{22}^{\;}\mathfrak{d}_{33}^{\;}-\mathfrak{d}_{23}^{\;2})+\mathfrak{d}_{22}^{\;}\mathfrak{d}_{31}^{\;2}+\mathfrak{d}_{33}^{\;}\mathfrak{d}_{12}^{\;2}-2\mathfrak{d}_{22}^{\;}\mathfrak{d}_{33}^{\;}\mathfrak{d}_{31}^{\;}\mathfrak{d}_{12}^{\;}}{\mathfrak{d}_{23}\mathfrak{c}_{23}}:\cdots:\cdots]_\Delta.
\]
\subsubsection{The star of a point with respect to the triple $\Delta$.} Define the mapping \vspace*{-1 mm}
%$\star: \mathcal{P} - \{A, B, C\} \rightarrow \mathcal{P}$ by\vspace*{-3 mm}\\ 
\begin{equation*}
\begin{split}&\star = \star_\Delta: \mathcal{P} - \{A, B, C\} \rightarrow \mathcal{P} \; \text{by}\\
&P^\star=\star([p_1{:}p_2{:}p_3]_\Delta)=[q_1{:}q_2{:}q_3]_\Delta\,,\; (q_1,q_2,q_3)=(p_2p_3,p_3p_1,p_1p_2)\mathfrak{D}.\vspace*{1 mm}
\end{split}
\end{equation*}
$P^\star$ is the dual of the tripolar of $P$, so the preimage of a point $Q \in \mathcal{P} - \{A', B', C'\}$ under $\star$ is the point $Q^{\tau \delta}$.\vspace*{1 mm} 
%The inverse of $\star$ is 
%\[
%\begin{split} 
%\star^{-1}: \mathcal{P} - \{A', B', C'\},  
%\star^{-1}([p_1{:}p_2{:}p_3]_\Delta)=([q_1{:}q_2{:}q_3]_\Delta\,)^{\tau \delta}}.\vspace*{1 mm}
%\]
%The euclidean limit of $P^\star$ is either the point $O$ or a point on the line at infinity.
If $\Delta_0$ is not right-angled, then we get for the special case $P = H$:\vspace*{-3 mm}\\
\begin{equation*}
H^\star=[\mathfrak{d}_{11} \mathfrak{d}_{23}{\,+\,}2 \mathfrak{d}_{31} \mathfrak{d}_{12}{\,:\,}\mathfrak{d}_{22} \mathfrak{d}_{31}{\,+\,}2 \mathfrak{d}_{12} \mathfrak{d}_{23}{\,:\,}\mathfrak{d}_{33} \mathfrak{d}_{12}{\,+\,}2 \mathfrak{d}_{23} \mathfrak{d}_{31}]_\Delta.\vspace*{1 mm}
\end{equation*}
The orthic axis and its dual are defined even if $\Delta_0$ has one or two right angles, so we extend $\star$ in this case. \\The point $H^\star$ is introduced in \cite{W1} under the name \textit{orthostar}; we adopt the name.\\ 
%The star of $H$ with respect to $\Delta'$, $\star_{\Delta'}(H)$, is the point $G^+$, and $G^+$ is also the preimage of $H$ under $\star$.

\begin{figure}[!htbp]
\includegraphics[height=9cm]{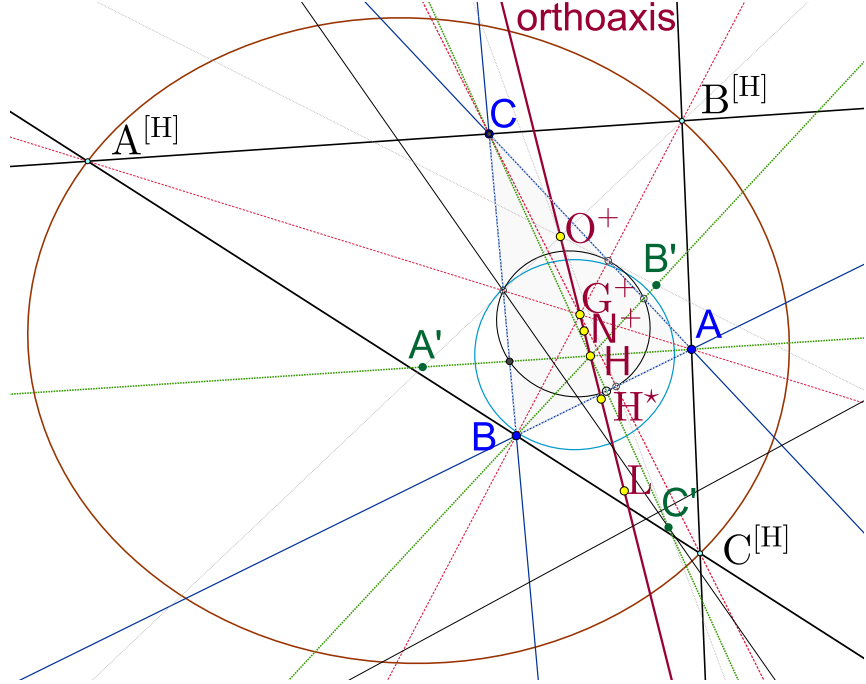}
\caption{Triangle centers on the orthoaxis. The three conics are the absolute conic (light blue), the circumcircle of $\Delta^{[H]}$ with center $H$ (brown) and the bicevian conic of the points $H$ and $G^+$ (black).
%\newline
%All figures were created with the software program GeoGebra [].
}
\end{figure}

\subsubsection{Centers on the orthoaxis.}\label{subsubsec:orthoaxis} Wildberger \cite{W1} showed that the line $H\vee H^\star$,
\[
\mathfrak{d}_{23}^{\;}(\mathfrak{d}_{22}^{\;} \mathfrak{d}_{31}^{\;2} - \mathfrak{d}_{33}^{\;} \mathfrak{d}_{12}^{\;2}) x_1 + \mathfrak{d}_{31}^{\;}(\mathfrak{d}_{33}^{\;} \mathfrak{d}_{12}^{\;2} - \mathfrak{d}_{11}^{\;} \mathfrak{d}_{23}^{\;2}) x_2 + \mathfrak{d}_{12}^{\;}(\mathfrak{d}_{11}^{\;} \mathfrak{d}_{23}^{\;2} - \mathfrak{d}_{22}^{\;} \mathfrak{d}_{31}^{\;2}) x_3 =0,
\]
 which he named \textit{orthoaxis}, is incident with the points $L, O^+, G^+$. Furthermore, he proved that $L,O^+\!,H,H^\star$ as well as $L,O^+\!,G^+\!,H$ form harmonic ranges. The euclidean limits of the points $L,O^+\!,$ $G^+\!,H,H^\star$ are $L=X_{20},O=X_3,G=X_2,H=X_4$, the Euler infinity point $X_{30}$, respectively.\vspace*{-2 mm}\\ 

The orthoaxis is a symmetry axis of the bicevian conic through the traces of $H$ and $G^+$. The points $H + O^+$ and $H - O^+$ are, in addition to the dual point of the orthoaxis, symmetry points of this conic precisely when these two points are proper midpoints of $H$ and $O^+$.\vspace*{1 mm}\\
\textit{Proof} of the last two statements: For two points $P = [p_1{:}p_2{:}p_3]_\Delta$ and $Q = [q_1{:}q_2{:}q_3]_\Delta$, not on any of the lines of $\Delta_0$, the conic which passes through the traces of $P$ and $Q$ has the matrix (cf. \cite{Y})\vspace*{-2 mm} 
	\[\mathfrak{M} = (\mathfrak{m}_{ij})_{i,j=1,2,3}\; \text{with}\; \mathfrak{m}_{ii} = -\frac {2}{p_i q_i}, \mathfrak{m}_{ij} = \frac {1}{p_i q_i}+\frac {1}{p_j q_j}\;\text{for}\; i,j = 1,2,3, i\ne j.
\]

(\textit{Remark}: For $P = Q$, this is the matrix of an inconic with perspector $P$.)\\

\noindent It can now be checked that for the special case $P = H$ and $Q = G^+$, the pole of the orthoaxis with respect to the bicevian conic is the point \vspace*{-1 mm}
\[ 
[\mathfrak{c}_{23}^{\;} (\mathfrak{d}_{22}^{\;} \mathfrak{d}_{31}^{\;2} - \mathfrak{d}_{12}^{\;2} \mathfrak{d}_{33}^{\;}), \mathfrak{c}_{31}^{\;} (\mathfrak{d}_{12}^{\;2} \mathfrak{d}_{33}^{\;} - \mathfrak{d}_{11}^{\;} \mathfrak{d}_{23}^{\;2}), \mathfrak{c}_{12}^{\;} (\mathfrak{d}_{11}^{\;} \mathfrak{d}_{23}^{\;2} - \mathfrak{d}_{22}^{\;} \mathfrak{d}_{31}^{\;2})]_\Delta.
\]
But this is also the dual point of the orthoaxis.\\
We show that $H\pm O^+$ is a point on the polar of $H\mp O^+$ with respect to the bicevian conic $\mathcal{C}(\mathfrak{M})$ of $H$ and $G^+$ by proving the correctness of the equation\\ $(H^\circ+O^{+\circ}){{\,\scriptscriptstyle{[\mathfrak{M}]}\,}}(H^\circ-O^{+\circ}) = 0$. We transform this equation equivalently:\vspace*{-3mm}\\
 \[
\begin{split}		
&\;\;\;\;\;\;\;\;\;\;\;\;\;\;\;\;\;\;(H^\circ + O^{+\circ}){{\,\scriptscriptstyle{[\mathfrak{M}]}\,}}(H^\circ - O^{+\circ}) =0\\
\Leftrightarrow\;\;&\;\;\;\;\;\;\;\;\;\;\;\;\;\;\;\;\;\;\;\;\;H^\circ{{\,\scriptscriptstyle{[\mathfrak{M}]}\,}}\, H^\circ - O^{+\circ}{{\,\scriptscriptstyle{[\mathfrak{M}]}\,}}\, O^{+\circ} =0\\
\Leftrightarrow\;\;&(\boldsymbol{h}\,{{\,\scriptscriptstyle{[\mathfrak{M}]}\,}}\boldsymbol{h})(\boldsymbol{o^+_{\,}}{{\,\scriptscriptstyle{[\mathfrak{C}]}\,}} \boldsymbol{o^+}) - (\boldsymbol{h}\,{{\,\scriptscriptstyle{[\mathfrak{C}]}\,}}\,\boldsymbol{h})(\boldsymbol{o^+_{\,}}{{\;\scriptscriptstyle{[\mathfrak{M}]}\,}} \boldsymbol{o^+}) =0\;\;\\ 
&\;\text{with}\; \boldsymbol{h} = ({\mathfrak{d}_{31}}{\mathfrak{d}_{12}}{\,,\,}{\mathfrak{d}_{12}}{\mathfrak{d}_{23}}{\,,\,}{\mathfrak{d}_{23}}{\mathfrak{d}_{31}})  \;\textrm{and}\;\boldsymbol{o^+} = ({\mathfrak{d}_{11}}{\mathfrak{d}_{23}}{\,,\,}{\mathfrak{d}_{22}}{\mathfrak{d}_{31}}{\,,\,}{\mathfrak{d}_{33}}{\mathfrak{d}_{12}}) .
\vspace*{-2mm}
\end{split}		
\]
The last equation holds, as can be checked with the help of a CAS.
Now it follows: Necessary and sufficient for $H + O^+$ and $H - O^+$ to be symmetry points of $\mathcal{C}(\mathfrak{M})$ is that their distance is $\frac{\pi}{2} i$. But this holds precisely when $H + O^+$ and $H - O^+$ are proper midpoints of $\{H, O^+\}$.
$\;\Box$\vspace*{-1mm}\\

The orthoaxes of the triples $ABC, AH\!B, BHC, {CH}\!A$ meet at the point \vspace*{1mm}\\
\centerline{$N^+:=[\mathfrak{d}_{11} \mathfrak{d}_{23} - 2 \mathfrak{d}_{12} \mathfrak{d}_{31}, \mathfrak{d}_{22} \mathfrak{d}_{31} - 2 \mathfrak{d}_{12} \mathfrak{d}_{23}, \mathfrak{d}_{33} \mathfrak{d}_{12} - 2 \mathfrak{d}_{23} \mathfrak{d}_{31}]_\Delta.$}\vspace*{-1mm}\\

The coordinates of $H, N^+, O^+,H^\star$ indicate that these four points form (in this order) a harmonic range. The euclidean limit of $N^+$ is $N = X_5$. \\

The orthic triple of the orthic triple $\Delta_{[H]}$ is perspective to $\Delta$ at perspector \vspace*{-1mm}
\[ P = [\mathfrak{d}_{31}^{\;}\mathfrak{d}_{12}^{\;}(\mathfrak{d}_{11}^{\;}\mathfrak{d}_{23}^{\;2}-\mathfrak{d}_{22}^{\;}\mathfrak{d}_{31}^{\;2}-\mathfrak{d}_{33}^{\;}\mathfrak{d}_{12}^{\;2}+2\mathfrak{d}_{12}^{\;}\mathfrak{d}_{23}^{\;}\mathfrak{d}_{31}^{\;}):\cdots:\cdots]_\Delta.
\]
This is a point on the orthoaxis with euclidean limit $X_{24}$.\vspace*{-1mm}\\

The de\! Longchamps point, the double point, the pseudo-circumcenter and the orthostar of $\Delta'$ are $O^+, H^\star, L$ and $G^+$, respectively.\vspace*{-1mm}\\

In case of $\varepsilon_{\!A}=\varepsilon_{\!B}=\varepsilon_{\!C}$ and $\varepsilon_{\!A'}=\varepsilon_{\!B'}=\varepsilon_{\!C'}$, the coordinates of the above given 
points on the orthoaxis can be written
\[
\begin{split}	
G^{+} &= [\sinh(\alpha)\big(\cosh(\alpha)+\cosh(\beta)\cosh(\gamma)\big):\cdots:\cdots]_\Delta\\
O^+ &= [\sinh(2\alpha):\cdots:\cdots]_\Delta,\\
H \;\;& = [\tanh(\alpha):\cdots:\cdots]_\Delta,\\  
N^+ &= [\sinh(\alpha)\big(\cosh(\alpha)+2\cosh(\beta)\cosh(\gamma)\big):\cdots:\cdots]_\Delta,\\ 
H^{\star} &= [\sinh(\alpha)\big(\cosh(\alpha)-2\cosh(\beta)\cosh(\gamma)\big):\cdots:\cdots]_\Delta,\\
L \;\;\,& = [\sinh(\alpha)\big(\cosh(\alpha)-\cosh(\beta)\cosh(\gamma)\big):\cdots:\cdots]_\Delta,\\
P \;\;\,&=[\tanh(\alpha)
\big(\cosh^2(\alpha)-\cosh^2(\beta)-\cosh^2(\gamma)-2\cosh(\alpha)\cosh(\beta)\cosh(\gamma)\big):\cdots:\cdots]_\Delta.
\end{split}	
\]
\subsubsection{The Taylor conic.} By projecting the pedals $A_{[H]}, B_{[H]}, C_{[H]}$, lying on $a\,$ resp. $\!b\,$ resp. $c$, onto the other sidelines of $\Delta$, we get altogether six points \\$B_{[A_{[H]}]} = [\mathfrak{d}_{12}^{\;2}{:} 0{:} \mathfrak{d}_{12}^{\;} \mathfrak{d}_{23}^{\;}{-}\,\mathfrak{d}_{22}^{\;} \mathfrak{d}_{31}^{\;}]_\Delta$, $C_{[A_{[H]}]}= [\mathfrak{d}_{13}^{\;2}{:}\mathfrak{d}_{31}^{\;} \mathfrak{d}_{23}^{\;}{-}\,\mathfrak{d}_{33}^{\;} \mathfrak{d}_{12}^{\;}{:}0]_\Delta$, $C_{[B_{[H]}]}$, $A_{[B_{[H]}]}$, $A_{[C_{[H]}]},B_{[C_{[H]}]}$ which all lie on a conic with the equation\vspace*{-1mm}

\[\sum\limits_{i=1}^3 \mathfrak{d}_{i+1,i+2}^{\;} x_i^{\;2} + \frac{(\mathfrak{d}_{i,i} \mathfrak{d}_{i+1,i+2}-\mathfrak{d}_{i+2,i} \mathfrak{d}_{i,i+1})^2 + \mathfrak{d}_{i+2,i} \mathfrak{d}_{i,i+1}}{\mathfrak{d}_{i,i} \mathfrak{d}_{i+1,i+2}-\mathfrak{d}_{i+2,i} \mathfrak{d}_{i,i+1}}x_{i+1} x_{i+2} = 0\;\;(\text{indices mod}\, 3).\vspace*{0mm}
\]
We like to call this conic \textit{Taylor conic}, since its euclidean limit is the Taylor circle. \vspace*{-2mm}\\

Let $H_A, H_B, H_C$ be the orthocenters of $AB_{[H]}C_{[H]}, BC_{[H]}A_{[H]}, CA_{[H]}B_{[H]}$, respectively. The triple $H_AH_BH_C$ is perspective to $A_{[H]}B_{[H]}C_{[H]}$; the perspector is \vspace*{-1mm}
\[[\mathfrak{d}_{11}^{\;} (\mathfrak{d}_{23}^{\;}(\mathfrak{d}_{22}^{\;} \mathfrak{d}_{31}^{\;2}+\mathfrak{d}_{33}^{\;} \mathfrak{d}_{12}^{\;2})- \mathfrak{d}_{22}^{\;}\mathfrak{d}_{33}^{\;}\mathfrak{d}_{31}^{\;}\mathfrak{d}_{12}^{\;}):\cdots:\cdots]_\Delta.\vspace*{0mm}
\] 
The euclidean limit of this point is $X_{389}$, the center of the Taylor circle. But, in general, this point is not a symmetry point of the Taylor conic.
\vspace*{0.5mm}
\subsubsection{Orthotransversal and orthocorrespondent of a point}$\!$For any point $P =$ $[p_1{:}p_2{:}p_3]_\Delta$ which is not a vertex of $\Delta$, the three points $\text{perp}(A\vee P,P)\wedge a ,  \text{perp}(B\vee P,P)\wedge b$ and $\text{perp}(C\vee P,P)\wedge c$ 
are collinear on the line\vspace*{-1mm}  
\[\sum\limits_{i=1,2,3} \,\frac{1}{p_i(\mathfrak{d}_{i+1,i+2}p_{i} - \mathfrak{d}_{i+2,i}p_{i+1}- \mathfrak{d}_{i,i+1}p_{i+2})+ \mathfrak{d}_{i,i}p_{i+1}p_{i+2}}\,x_i = 0\;\;\;\;\;\;(\text{indices mod}\, 3).\vspace*{0mm}
\] 
The line is called the \textit{orthotransversal} (line) of $P$, and its tripole is called the \textit{orthocorrespondent} of $P$.
As a special case, the orthocorrespondent of $H$ is $G^+$.\vspace*{-2mm} \\

The points $\text{perp}(A\vee P,P)\wedge a' ,  \text{perp}(B\vee P,P)\wedge b'$ and $\text{perp}(C\vee P,P)\wedge c'$ are also collinear; they lie on the dual of $P$.\vspace*{1mm} \\
We dualize these statements:\vspace*{-2mm} \\

The three lines $\text{perp}(a\wedge P^\tau,P^\tau)\vee A ,  \text{perp}(b\wedge P^\tau,P^\tau)\vee B$ and $\text{perp}(c\wedge P^\tau,P^\tau)\vee C$ meet at the point $Q=[\displaystyle\frac{p_1}{p_1(\mathfrak{c}_{11} p_1 - \mathfrak{c}_{12} p_2 - \mathfrak{c}_{31} p_3)+\mathfrak{c}_{23} p_2 p_3}:\cdots:\cdots]_\Delta$.\vspace*{1mm}\\
A special case: If $P = G^+$, then $Q = H$.\vspace*{-1mm} \\

If $P$ is not a vertex of $\Delta'$, then the points $\text{par}(a,P)\wedge a, \text{par}(b,P)\wedge b, \text{par}(c,P)\wedge c$ lie on the line $P^\delta$.
 
\subsection{Center functions / triangle centers and their mates}\hspace*{\fill}\vspace*{1mm}\\
$\hspace*{3mm}$Given a triangle center $T$ of $\Delta_0$, then there exists a homogenous function ${f\!\!:}\; \mathbb{R}^{6}\to \mathbb{R}\,$ such that\\
$\bullet\;  f(x_1,x_2,x_3,x_4,x_5,x_6) = f(x_1,x_3,x_2,x_4,x_6,x_5)$ \hspace*{2 mm}and\\
$\bullet\; T\! = [f(\mathfrak{d}_{11},\!\mathfrak{d}_{22},\!\mathfrak{d}_{33},\!\mathfrak{d}_{23},\!\mathfrak{d}_{31},\!\mathfrak{d}_{12}){:}f(\mathfrak{d}_{22},\!\mathfrak{d}_{33},\!\mathfrak{d}_{11},\!\mathfrak{d}_{31},\!\mathfrak{d}_{12},\!\mathfrak{d}_{23}){:}f(\mathfrak{d}_{33},\!\mathfrak{d}_{11},\!\mathfrak{d}_{22},\!\mathfrak{d}_{12},\!\mathfrak{d}_{23},\!\mathfrak{d}_{31})]_\Delta$.\vspace*{1 mm}\\
A function satisfying these conditions is called a center function of $T$ with respect to the triangle $\Delta_0$ and the matrix $\mathfrak{D}$. \\
We give two examples: Center functions of $H$ and $L$ are $f_{\!H}(x_1,x_2,x_3,x_4,x_5,x_6) = x_5 x_6$ and $f_{\!L}(x_1,x_2,x_3,x_4,x_5,x_6) = x_1 x_4 + x_5 x_6.$\\
Obviously, by knowing a center function of a triangle center of $\Delta_0$, one also knows the barycentric coordinates of that center.\\
The triangle center $T = T_0$ of $\Delta_0$ is accompanied by its three mates $T_1,T_2,T_3$; these are the corresponding triangle centers of $\Delta_1,\Delta_2,\Delta_3$, respectively. Their representation by barycentric coordinates are\vspace*{0.5mm}\\

\noindent $T_1 = [-f(\mathfrak{d}_{11},\!\mathfrak{d}_{22},\!\mathfrak{d}_{33},\!\mathfrak{d}_{23},\!-\mathfrak{d}_{31},-\mathfrak{d}_{12})$\\
$\hspace*{29mm}{:}f(\mathfrak{d}_{22},\!\mathfrak{d}_{33},\!\mathfrak{d}_{11},-\mathfrak{d}_{31},-\mathfrak{d}_{12},\!\mathfrak{d}_{23}){:}f(\mathfrak{d}_{33},\!\mathfrak{d}_{11},\!\mathfrak{d}_{22},-\mathfrak{d}_{12},\!\mathfrak{d}_{23},-\mathfrak{d}_{31})]_\Delta$,\\
$T_2 = [f(\mathfrak{d}_{11},\!\mathfrak{d}_{22},\!\mathfrak{d}_{33},-\mathfrak{d}_{23},\!\mathfrak{d}_{31},-\mathfrak{d}_{12})$\\
$\hspace*{29mm}{:}-f(\mathfrak{d}_{22},\!\mathfrak{d}_{33},\!\mathfrak{d}_{11},\!\mathfrak{d}_{31},-\mathfrak{d}_{12},-\mathfrak{d}_{23}){:}f(\mathfrak{d}_{33},\!\mathfrak{d}_{11},\!\mathfrak{d}_{22},-\mathfrak{d}_{12},\!\mathfrak{d}_{23},-\mathfrak{d}_{31})]_\Delta$,\\
$T_3 = [f(\mathfrak{d}_{11},\!\mathfrak{d}_{22},\!\mathfrak{d}_{33},-\mathfrak{d}_{23},-\mathfrak{d}_{31},\!\mathfrak{d}_{12})$\\
$\hspace*{29mm}{:}f(\mathfrak{d}_{22},\!\mathfrak{d}_{33},\!\mathfrak{d}_{11},-\mathfrak{d}_{31},\!\mathfrak{d}_{12},-\mathfrak{d}_{23}){:}-f(\mathfrak{d}_{33},\!\mathfrak{d}_{11},\!\mathfrak{d}_{22},\!\mathfrak{d}_{12},-\mathfrak{d}_{23},-\mathfrak{d}_{31})]_\Delta$.\vspace*{1.5mm}\\
All triangle centers in the last subsection are \textit{absolute triangle centers}, their mates do not differ from the main center.\vspace*{1mm}\\
\textit{Remark}: Without giving a formal introduction, we also use (and used already) center functions which depend on the matrix $\mathfrak{C}$, or on the sidelengths of $\Delta_0$, or on its inner angles.\\

\subsection{Circumcircles, incircles and related triangle centers}
\subsubsection{Twin circles, circumcircles and incircles}
\begin{figure}[!htbp]
\includegraphics[height=9cm]{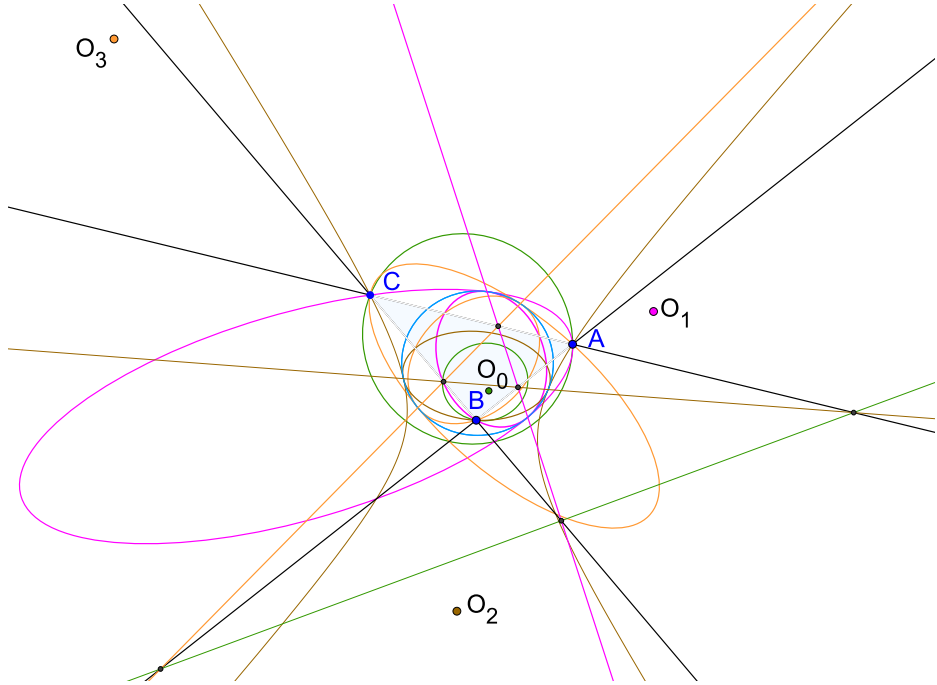}
\caption{The four pairs of twin circumcircles (green, purple, brown, orange) of $ABC$, their centers and the dual lines of these centers. The blue circle is the absolute conic.}
\end{figure}
The points $G\;{=}\;G_0\;{=}\;[1,1,1]_\Delta,\, G_1\;{=}\;\hspace*{0.2mm}$ $[-1,1,1]_\Delta, G_2\,{=}\,[1,-1,1]_\Delta$, $G_3\,{=}\,[1,1,-1]_\Delta$ are the \textit{centroids} of triangle $\Delta_0,\Delta_1,\Delta_2,\Delta_3$, respectively. The point $O_{k}{\,=\,}G_{k}^{\star}$ is called \textit{circumcenter} of $\Delta_k\;\;(k = 0,{\cdots},3)$, cf. \cite{W2}. In the elliptic plane, there always exists a circle around $O_k$, passing through the vertices $A, B, C$. The situation is different in the hyperbolic plane. Unless $\varepsilon_{\!A}{\,=\,}\varepsilon_{\!B}{\,=\,}\varepsilon_{\!C}$, there does not exist any circle which passes through all three vertices $A, B, C$. If for example $\varepsilon_{\!A}{\,=\,}\varepsilon_{\!B}{\,\ne\,}\varepsilon_{C}$, then a circle with center $O_k$ through $A$ will also pass through {B}, but it differs from a circle centered at $O_k$ passing through $C$. Wildberger and Alkhaldi discovered a special relation between these two circles: They form a couple of \textit{twin circles}. Wildberger and Alkhaldi \cite{W2} proved that, given any circle $\mathcal{C}$ with center $P$ and radius $r \in \mathbb{C}$ in the hyperbolic plane, there exists a uniquely determined twin circle of $\mathcal{C}$. This is a circle with center $P$ and radius $r'$, satisfying the equation $\cosh^2(r') = - \cosh^2(r)$. Thus, in the hyperbolic plane, one can always find exactly one pair of twin circles with center $O_k$ such that their union contains the vertices $A,B,C$.\\The coordinates of the circumcenters of $\Delta$ are \vspace*{0mm}
 \[
\begin{split}	&O{\,=\,}O_0{\,=\,}[\mathfrak{d}_{11}{+\,}\mathfrak{d}_{12}{+\,}\mathfrak{d}_{31}, \mathfrak{d}_{12}{+\,}\mathfrak{d}_{22}{+\,}\mathfrak{d}_{23}, \mathfrak{d}_{31}{+\,}\mathfrak{d}_{23}{+\,}\mathfrak{d}_{33}]_\Delta,\\  
&O_1{\,=\,}[{-\,}\mathfrak{d}_{11}{+\,}\mathfrak{d}_{12}{+\,}\mathfrak{d}_{31}, {-\,}\mathfrak{d}_{12}{+\,}\mathfrak{d}_{22}{+\,}\mathfrak{d}_{23}, -\mathfrak{d}_{31}{+\,}\mathfrak{d}_{23}{+\,}\mathfrak{d}_{33}]_\Delta,\; \text{etc}.
\end{split}
\]
For the distances $d(O,A)$ and $d(O_1,A)$ we get:
\[
\begin{split}	
&\cosh^2(d(O,A)) = {1}/{(\mathfrak{c}_{11}(\mathfrak{d}_{11}{+\,}\mathfrak{d}_{22}{+\,}\mathfrak{d}_{33}+2(\mathfrak{d}_{23}{+\,}\mathfrak{d}_{31}{+\,}\mathfrak{d}_{12})))},\\ 
&\cosh^2(d(O_1,A)) = {1}/{(\mathfrak{c}_{11}(\mathfrak{d}_{11}{+\,}\mathfrak{d}_{22}{+\,}\mathfrak{d}_{33}+2(\mathfrak{d}_{23}{-\,}\mathfrak{d}_{31}{-\,}\mathfrak{d}_{12})))}.\vspace*{-1mm} 
\end{split}\]
The circumcenters of $\Delta'$ are the incenters of $\Delta$. To be more precise, the circumcenter of $\Delta^k$ is the incenter $I_k$ of $\Delta_k, k = 0,\cdots,3$.
We calculate the coordinates:
\[
I = I_0 = [\sqrt{|\mathfrak{d}_{11}|} ,\sqrt{|\mathfrak{d}_{22}|},\sqrt{|\mathfrak{d}_{33}|}\,]_\Delta,\;\; I_1 = [-\sqrt{|\mathfrak{d}_{11}|},\sqrt{|\mathfrak{d}_{22}|},\sqrt{|\mathfrak{d}_{33}|}\,]_\Delta, \cdots.
\]
In the elliptic plane, the incenters are the centers of incircles of $\Delta$, and this is also true in the hyperbolic plane if $\varepsilon_{\!A'}{\,=\,}\varepsilon_{\!B'}{\,=\,}\varepsilon_{C'}$. But in general, twin circles with center $I_k$ are needed to touch all sidelines $a, b, c$.
We calculate the distances $d(I,a) = d(I, \text{ped}(I,a)) = \pi i- d(I,A')$  and $d(I_1,a)\!:\vspace*{-2mm}$
\[
\begin{split}
&\cosh^2(d(I,a)) = \textrm{sgn}(\mathfrak{d}_{11})/\big(\!\!\!\sum\limits_{\,i=1,2,3}\!(\mathfrak{c}_{ii}|\mathfrak{d}_{ii}|+2 \mathfrak{c}_{i,i+1}\sqrt{|\mathfrak{d}_{ii}\mathfrak{d}_{i+1,i+1}|}\,)\big) ,\\
&\cosh^2(d(I_1,a)) = \textrm{sgn}(\mathfrak{d}_{11})/\big(\!\!\!\sum\limits_{\,i=1,2,3}\!(\mathfrak{c}_{ii} |\mathfrak{d}_{ii}|) + 2 (\mathfrak{c}_{23}\sqrt{|\mathfrak{d}_{22}\mathfrak{d}_{33}|}{-}\,\mathfrak{c}_{31}\!\sqrt{|\mathfrak{d}_{33}\mathfrak{d}_{11}|}{-}\,\mathfrak{c}_{12}\!\sqrt{|\mathfrak{d}_{11}\mathfrak{d}_{22}|})\big).
\end{split}	
\]

The triple $O_1O_2O_3$ is perspective to $\Delta'$ at perspector $O_0$, the triple $I_1I_2I_3$ is perspective to $\Delta$ at perspector $I_0$.\vspace*{-1mm} \\ 

We now assume $\varepsilon_{\!A}{\,=\,}\varepsilon_{\!B}{\,=\,}\varepsilon_{C}$ and put $\mathfrak{c}\!:=\!\mathfrak{c}_{11}{\,=\,}\mathfrak{c}_{22}{\,=\,}\mathfrak{c}_{33} \in \{-1, 1\}\,$, then :\vspace*{1mm}\\
\centerline{$ \mathfrak{c}\mathfrak{c}_{23} = \cosh(\mathscr{a}), \; \mathfrak{c}\mathfrak{c}_{23} + 1 = 2 \cosh^2(\frac{\mathscr{a}}{2}), \; \mathfrak{c}\mathfrak{c}_{23} - 1 = 2 \sinh^2(\frac{\mathscr{a}}{2}),\, \textrm{etc.}$\vspace*{-1mm}}
and\vspace*{1mm}\\
\centerline{\;\;$O_0 = [(\mathfrak{c}_{23}{-}\mathfrak{c})(\mathfrak{c}_{23}{-}\mathfrak{c}_{31}{-}\mathfrak{c}_{12}{+}\mathfrak{c}){\,:\,}(\mathfrak{c}_{31}{-}\mathfrak{c})({-}\mathfrak{c}_{23}{+}\mathfrak{c}_{31}{-}\mathfrak{c}_{12}{+}\mathfrak{c}){\,:\,}(\mathfrak{c}_{12}{-}\mathfrak{c})({-}\mathfrak{c}_{23}{-}\mathfrak{c}_{31}{+}\mathfrak{c}_{12}{+}\mathfrak{c})]_\Delta\;\;$ 
}
\centerline{$ \;\;\,\, = [(\sinh^2(\frac{\mathscr{a}}{2}))({-}\sinh^2(\frac{\mathscr{a}}{2})+\sinh^2(\frac{\mathscr{b}}{2})+\sinh^2(\frac{\mathscr{c}}{2})):\cdots:\cdots]_\Delta,\;\;\;\;\;\;\;\;\;\;\;\;\;\;\;\;\;\;\;\;\;\;\;\;\;\;\;\;\;\;\;\vspace*{-2.5mm}$}\\
\centerline{$O_1 = [(\mathfrak{c}_{23}{-}\mathfrak{c})(\mathfrak{c}_{23}{+}\mathfrak{c}_{31}{+}\mathfrak{c}_{12}{+}\mathfrak{c}){\,:\,}(\mathfrak{c}_{31}{+}\mathfrak{c})({-}\mathfrak{c}_{23}{-}\mathfrak{c}_{31}{+}\mathfrak{c}_{12}{+}\mathfrak{c}){\,:\,}(\mathfrak{c}_{12}{+}\mathfrak{c})({-}\mathfrak{c}_{23}{+}\mathfrak{c}_{31}{-}\mathfrak{c}_{12}{+}\mathfrak{c})]_\Delta,$}\\
\centerline{$ \;\;= [(\sinh^2(\frac{\mathscr{a}}{2}))(\sinh^2(\frac{\mathscr{a}}{2}){+}\cosh^2(\frac{\mathscr{b}}{2}){+}\cosh^2(\frac{\mathscr{c}}{2})):(\cosh^2(\frac{\mathscr{b}}{2}))({-}\sinh^2(\frac{\mathscr{a}}{2}){-}\cosh^2(\frac{\mathscr{b}}{2})$}\\
\centerline{$\hspace*{31mm}{+}\cosh^2(\frac{\mathscr{c}}{2})):(\cosh^2(\frac{\mathscr{c}}{2}))({-}\sinh^2(\frac{\mathscr{a}}{2}){+}\cosh^2(\frac{\mathscr{b}}{2}){-}\cosh^2(\frac{\mathscr{c}}{2}))]_\Delta,
\cdots.$}\vspace*{-1.5mm}\\

The four points $O_0,O_1,O_2,O_3$ form an orthocentric system.\vspace*{-1.5mm}\\

The triple $O_1O_2O_3$ is perspective to $\Delta$ at the isogonal conjugate of $O = O_0$, and as a consequence of this, the triples ${O_0} {O_3} {O_2}$, ${O_3} {O_0} {O_1}$, ${O_2} {O_1} {O_0}$ are  perspective to $\Delta$ at the isogonal conjugates of ${O}_1, {O}_2, {O}_3$, respectively.\vspace*{-1mm}\\

Let $\mathscr{C}_k$ denote the circumcircle with center $O_k, k = 0,\cdots,3$. All inner points of $\Delta_0$ lie inside $\mathscr{C}_0$, which is not true for any other circumcircle $\mathscr{C}_k, k\ne 0.$ Therefore, we  call $\mathscr{C}_0$ the \textit{main circumcircle} of $\Delta_0$, while the others will be regarded as its \textit{mates}.\vspace*{-1mm}\\

The perspector of $\mathscr{C}_k$ is called the \textit{Lemoine point} $\tilde{K}_k$ of $\Delta_k:\vspace*{1mm}$\\
%\begin{split}
\centerline{$\tilde{K}_0=[\mathfrak{c}_{23}-\mathfrak{c}:\mathfrak{c}_{31}-\mathfrak{c}:\mathfrak{c}_{12}-\mathfrak{c}]_\Delta = [\sinh^2(\frac{\mathscr{a}}{2}):\sinh^2(\frac{\mathscr{b}}{2})\,:\sinh^2(\frac{\mathscr{c}}{2})]_\Delta$,\;\;\;\,\;\;\;\;\;}
\centerline{$\tilde{K}_1=[\mathfrak{c}_{23}-\mathfrak{c}:\mathfrak{c}_{31}+\mathfrak{c}:\mathfrak{c}_{12}+\mathfrak{c}]_\Delta = [-\sinh^2(\frac{\mathscr{a}}{2}):\cosh^2(\frac{\mathscr{b}}{2}):\cosh^2(\frac{\mathscr{c}}{2})]_\Delta,\;....$}\vspace*{1mm} \\
The isogonal conjugate of $\tilde{K}_k$ is  $\tilde{G}_k$ :\vspace*{1mm}\\
\centerline{$\tilde{G}_0=[\mathfrak{c}_{23}+\mathfrak{c}:\mathfrak{c}_{31}+\mathfrak{c}:\mathfrak{c}_{12}+\mathfrak{c}]_\Delta = [\cosh^2(\frac{\mathscr{a}}{2}):\cosh^2(\frac{\mathscr{b}}{2}):\cosh^2(\frac{\mathscr{c}}{2})]_\Delta$\;,\;\;\;\;\;\;\;\;}\\
\centerline{$\;\; \tilde{G}_1=[\mathfrak{c}_{23}+\mathfrak{c}:\mathfrak{c}_{31}-\mathfrak{c}:\mathfrak{c}_{12}-\mathfrak{c}]_\Delta =[-\cosh^2(\frac{\mathscr{a}}{2}):\sinh^2(\frac{\mathscr{b}}{2}):\,\sinh^2(\frac{\mathscr{c}}{2})]_\Delta$\;,\;....\;\;\;}\vspace*{1mm}\\ 
The triple $\tilde{K_1} \tilde{K_2} \tilde{K_3}$ is perspective to $\Delta$ at $\tilde{G}_0$ and, as a consequence, 
the triple $\tilde{G_1} \tilde{G_2} \tilde{G_3}$ is perspective to $\Delta$ at $\tilde{K}_0$. The triple $\tilde{G_2} \tilde{G_3} \tilde{G_0}$ is perspective to $\Delta$ at $\tilde{K}_1$, etc. The four lines $\tilde{K}_j\vee \tilde{G}_j , 0\le j\le 3$, meet at the point $G = [1,1,1]_\Delta.$\vspace*{-0.5mm} \\

The euclidean limits of $\tilde{K}=\tilde{K}_0$ and $\tilde{G}=\tilde{G}_0$ are the symmedian $K = X_6$ and the centroid $G = X_2$ of $\Delta_0$, respectively. The euclidean limit of the circle $\mathscr{C}_1$ is the union of two lines, one is the sideline $a = B\vee C$, the other its parallel through $A$.\vspace*{3mm} \\
\textit{Remark:} In euclidean geometry, the perspector of the circumcircle $\mathscr{C}_0$ coincides with iso\-gonal conjugate of $G_0$. This is not the case in elliptic and in hyperbolic geometry; so we have to find different names for the two points. As proposed by Horv\'ath, the perspector of the circumcircle will be called \textit{Lemoine point}, the isogonal conjugate $K$ of $G$  \textit{symmedian}. In euclidean geometry, the centroid of the pedal triangle $\Delta_{[K],0}$ of $K$ agrees with $K$. Neither of the two points $K, \tilde{K}$ has this property in elliptic and in hyperbolic geometry.\vspace*{1mm}\\
The Lemoine point $\tilde{K}$ we can get by the following construction, cf. \cite{Gr} for the euclidean case: Let $l_A$ be the line that connects the point $A_G$ with the harmonic conjugate of $A'$ with respect to $A,A_H$,  and define the lines $l_B,l_C$ accordingly. The lines $l_A,l_B,l_C$ meet at $\tilde{K}$.\vspace*{-0.5mm}\\
%\noindent\hspace*{0mm}The triple $\tilde{K}_{i+1} \tilde{K}_{i+2} \tilde{K}_{i+3}$ is perspective to $\Delta$ at perspector $\tilde{G}_k$, and
%the triple $\tilde{G}_{i+1} \tilde{G}_{i+2} \tilde{G}_{i+3}$ is perspective to $\Delta$ at perspector $\tilde{K}_i, i=0,1,2,3$, indices mod 3.\vspace*{1mm}\\

If  $\mathcal{C}$ is any conic and $P$ a point on $\mathcal{C}$, then let $T_{\!P}\mathcal{C}$ denote the tangent of $\mathcal{C}$ at $P$. The points $A, T_B\mathscr{C}_2\wedge T_C\mathscr{C}_3, T_B\mathscr{C}_1\wedge T_C\mathscr{C}_1$ are collinear.\vspace*{-1mm}\\

Let $Q_{0j}$ be the non trivial intersection point of $\mathscr{C}_0$ and $\mathscr{C}_{\!j}, j{=}1,2,3.$ The triple $Q_{01}Q_{02}Q_{03}$ is perspective to $\Delta$ at $G_0$.\vspace*{1mm}\\

If we assume $\varepsilon_{\!A'} = \varepsilon_{\!B'} = \varepsilon_{C'}$, we get:\vspace*{1mm}\\
\centerline{$\textrm{sgn}(\mathfrak{d}_{11})= \textrm{sgn}(\mathfrak{d}_{22})= \textrm{sgn}(\mathfrak{d}_{33})$}\vspace*{-1mm}\\
and\vspace*{0mm}\\
%\begin{split}
\centerline{$I\, = I_0 = [\sqrt{\mathfrak{d}_{11}} :\sqrt{\mathfrak{d}_{22}}:\sqrt{\mathfrak{d}_{33}}\,]_\Delta = [\sinh(\mathscr{a}):\sinh(\mathscr{b}):\sinh(\mathscr{c})]_\Delta,$\;\;}\\
\centerline{$\;\;\;I_1 \!= [-\sqrt{\mathfrak{d}_{11}}:\sqrt{\mathfrak{d}_{22}}:\sqrt{\mathfrak{d}_{33}}\,]_\Delta = [-\sinh(\mathscr{a}):\sinh(\mathscr{b}):\sinh(\mathscr{c})]_\Delta,...$\;\;\;\;}\vspace*{-3mm}\\
%&I_1 = [(\mathfrak{c}_{23}{-}c)({-}\mathfrak{c}_{23}{+}\mathfrak{c}_{31}{+}\mathfrak{c}_{12}{+}c){\,:\,}(\mathfrak{c}_{31}{+}c)({-}\mathfrak{c}_{23}{-}\mathfrak{c}_{31}{+}\mathfrak{c}_{12}{+}c){\,:\,}(\mathfrak{c}_{12}{+}c)({-}\mathfrak{c}_{23}{+}\mathfrak{c}_{31}{-}\mathfrak{c}_{12}{+}c)]_\Delta.\\
%\end{split}

$I_0,I_1,I_2,I_3$ form an orthocentric system.\vspace*{-1.5mm}\\

\begin{figure}[!htbp] 
\includegraphics[height=8.5cm]{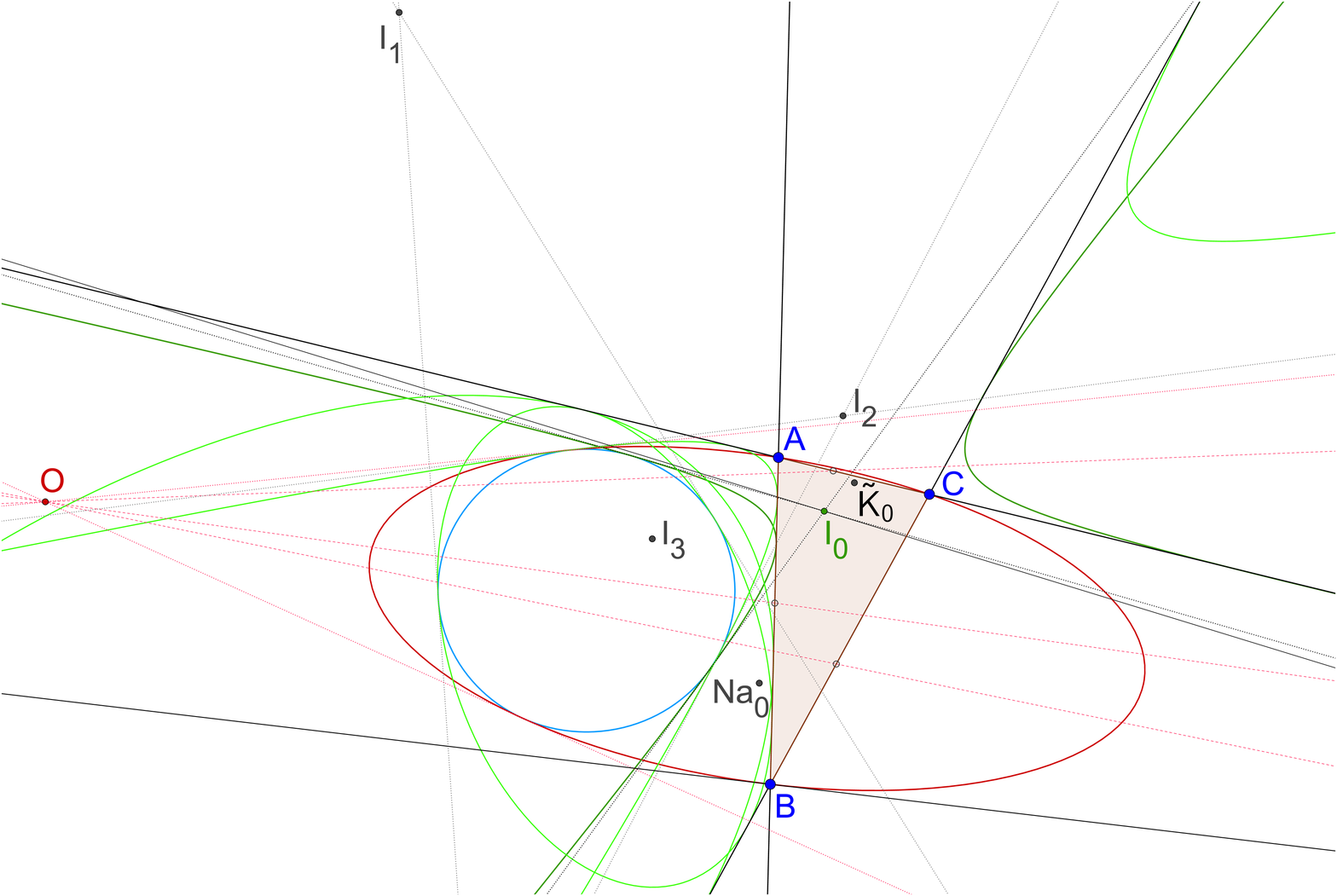} 
\caption{An anti-$\!$ de$\!$ Sitter triangle together with its circumcircle $\mathscr{C}_0$ (red), its incircle $\mathscr{I}_0$ (green) and excircles (light-green).}\label{fig:de Sitter triangle}
\end{figure}

Let $\mathscr{I}_k$ denote the incircle with center $I_k, k = 0,\cdots,3$. The center $I_k$ of $\mathscr{I}_k$ is always a point inside of $\Delta_k$. We call $\mathscr{I}_k$  the \textit{proper incircle} of $\Delta_k$, while the others will be called the \textit{excircles} of $\Delta_k$.  Caution: The inner points of $\mathscr{I}_k$ can completely lie outside of $\Delta_k$, cf. Figure \ref{fig:de Sitter triangle}.\vspace*{-2mm}\\

The perspector of $\mathscr{I}_k$ is called the Gergonne point $Ge_k$ of $\Delta_k$:\vspace*{1.5mm}\\
$\hspace*{4.5mm}{Ge_0}=[\sqrt{\mathfrak{d}_{22}}\sqrt{\mathfrak{d}_{33}} + \mathfrak{d}_{23}:\sqrt{\mathfrak{d}_{33}}\sqrt{ \mathfrak{d}_{11}} + \mathfrak{d}_{31}:\sqrt{\mathfrak{d}_{11}}\sqrt{\mathfrak{d}_{22}} + \mathfrak{d}_{12}]_\Delta$\vspace*{1.3mm}\\
$\hspace*{11mm}=[(\sqrt{\mathfrak{d}_{33}}\sqrt{ \mathfrak{d}_{11}} - \mathfrak{d}_{31})(\sqrt{\mathfrak{d}_{11}}\sqrt{\mathfrak{d}_{22}} - \mathfrak{d}_{12}):(\sqrt{\mathfrak{d}_{11}}\sqrt{\mathfrak{d}_{22}} - \mathfrak{d}_{12})(\sqrt{\mathfrak{d}_{22}}\sqrt{\mathfrak{d}_{33}} - \mathfrak{d}_{23})$\vspace*{1.3mm}\\
         $\hspace*{14.5mm}:(\sqrt{\mathfrak{d}_{33}}\sqrt{ \mathfrak{d}_{11}} - \mathfrak{d}_{31}))(\sqrt{\mathfrak{d}_{22}}\sqrt{\mathfrak{d}_{33}} - \mathfrak{d}_{23})]_\Delta$\vspace*{0.5mm}\\
$\hspace*{11mm}=[\displaystyle\sqrt{\mathfrak{c}_{22}}\sqrt{\mathfrak{c}_{33}}\, \frac{\cosh(\alpha){-}1}{\sinh(\mathscr{a})}:\sqrt{\mathfrak{c}_{33}}\sqrt{ \mathfrak{c}_{11}}\, \frac{\cosh(\beta){-}1}{\sinh(\mathscr{b})}:\sqrt{\mathfrak{c}_{11}}\sqrt{\mathfrak{c}_{22}}\, \frac{\cosh(\gamma){-}1}{\sinh(\mathscr{c})}]_\Delta,$\vspace*{1mm}\\
$\hspace*{4.5mm}{Ge_1}=[-(\sqrt{\mathfrak{d}_{22}}\sqrt{\mathfrak{d}_{33}} + \mathfrak{d}_{23}):\sqrt{\mathfrak{d}_{33}}\sqrt{ \mathfrak{d}_{11}} - \mathfrak{d}_{31}:\sqrt{\mathfrak{d}_{11}}\sqrt{\mathfrak{d}_{22}} - \mathfrak{d}_{12}]_\Delta$\vspace*{-2mm}\\

We introduce the Nagel point $N\!a_k$ of $\Delta_k, k=0,1,2,3$ , by \vspace*{1mm}\\
$\hspace*{4.5mm}{N\!a_0}=[\sqrt{\mathfrak{d}_{22}}\sqrt{\mathfrak{d}_{33}} - d_{23}:\sqrt{\mathfrak{d}_{33}}\sqrt{\mathfrak{d}_{11}} - d_{31}:\sqrt{\mathfrak{d}_{11}}\sqrt{\mathfrak{d}_{22}} - d_{12}]_\Delta$\vspace*{1mm}\\
$\hspace*{8.5mm}=[(\sqrt{\mathfrak{d}_{33}}\sqrt{\mathfrak{d}_{11}} + d_{31})(\sqrt{\mathfrak{d}_{11}}\sqrt{\mathfrak{d}_{22}} + d_{12}):(\sqrt{\mathfrak{d}_{11}}\sqrt{\mathfrak{d}_{22}} + d_{12})(\sqrt{\mathfrak{d}_{22}}\sqrt{\mathfrak{d}_{33}} + d_{23})$\vspace*{1mm}\\
$\hspace*{11mm} :(\sqrt{\mathfrak{d}_{33}}\sqrt{\mathfrak{d}_{11}} + d_{31}))(\sqrt{\mathfrak{d}_{22}}\sqrt{\mathfrak{d}_{33}} + d_{23})]_\Delta,\vspace*{1mm}$\\
$\hspace*{4.5mm}{N\!a_1}=[-(\sqrt{\mathfrak{d}_{22}}\sqrt{\mathfrak{d}_{33}} - \mathfrak{d}_{23}):\sqrt{\mathfrak{d}_{33}}\sqrt{ \mathfrak{d}_{11}} + \mathfrak{d}_{31}:\sqrt{\mathfrak{d}_{11}}\sqrt{\mathfrak{d}_{22}} + \mathfrak{d}_{12}]_\Delta$.\vspace*{-2mm}\\

\noindent\vspace*{0.5mm}$N\!a_k$ and $Ge_k$  are isotomic conjugates for $ 0\le\! k\! \le 3$. The triples $Ge_{1}Ge_{2}Ge_{3}$ and $N\!a_{1}N\!a_{2}N\!a_{3}$ are perspective to $\Delta$ at $N\!a_0$ and $Ge_0$, respectively.\vspace*{-1mm}\\

The lines $I_0 \vee Ge_0, I_1 \vee Ge_1,I_2 \vee Ge_2,I_3 \vee Ge_3$ concur at the de\! Longchamps point $L$,  the lines $I_0 \vee N\!a_0, I_1 \vee N\!a_1,I_2 \vee N\!a_2,I_3 \vee N\!a_3$ at the point $G^+$, and the lines $Ge_0 \vee N\!a_0,$  $Ge_1 \vee N\!a_1,Ge_2 \vee N\!a_2,Ge_3 \vee N\!a_3$ meet at the isotomic conjugate of the orthocenter $H$.

\subsubsection{Note:} \underline{From now on, we always assume  $\hspace*{0mm}\varepsilon_{\!A}=\varepsilon_{\!B}=\varepsilon_{C}$ and $\varepsilon_{\!A'}=\varepsilon_{\!B'}=\varepsilon_{C'}$.}\vspace*{1mm}\\ This is the "classical case". The points $A,B,C$ lie either in the elliptic plane or in a special part of the extended hyperbolic plane.\vspace*{3mm}

\centerline{The classical Cayley-Klein cases in the extended hyperbolic plane:}\vspace*{1mm}
\begin{center}
\noindent
\begin{tabular}{ c| c}\hline\hline
$\varepsilon_{\!A}=\varepsilon_{\!B}=\varepsilon_{C}=-i$ and $\varepsilon_{\!A'}=\varepsilon_{\!B'}=\varepsilon_{C'}=1$ & proper hyperbolic / Lobachevsky\\ \hline
$\varepsilon_{\!A}=\varepsilon_{\!B}=\varepsilon_{C}=1$ \;\,and $\varepsilon_{\!A'}=\varepsilon_{\!B'}=\varepsilon_{C'}=1$ & double-hyperbolic / de\! Sitter \\ \hline
$\varepsilon_{\!A}=\varepsilon_{\!B}=\varepsilon_{C}=1$ and $\varepsilon_{\!A'}=\varepsilon_{\!B'}=\varepsilon_{C'}=-i$ & dual-hyperbolic / anti-de\! Sitter\\\hline 
 \end{tabular}\noindent\\
\end{center}
\vspace*{3mm}
In this case, the following trigonometric formulae apply:\vspace*{1.5mm}\\
\noindent $\hspace*{4mm}\cosh(\alpha) \,= \displaystyle\frac{\cosh(\mathscr{b})\,\cosh(\mathscr{c})-\cosh(\mathscr{a})}{\sinh(\mathscr{b}) \sinh(\mathscr{c})}$\vspace*{1mm}\\
\noindent $\hspace*{4mm}\cosh(\mathscr{a}) = \displaystyle\frac{\cosh(\beta)\,\cosh(\gamma)+\cosh(\alpha)}{\sinh(\beta) \sinh(\gamma)} = 1 + \frac{2\sinh(\epsilon)\cosh(\epsilon{-}\alpha)}{\sinh(\beta) \sinh(\gamma)},\;\ \epsilon = \frac12 \textrm{area}(\Delta_0) .$\vspace*{1mm}\\
\noindent $\hspace*{4mm}\displaystyle\frac{\sinh(\alpha)}{\sinh(\mathscr{a})} = \frac{\sinh(\beta)}{\sinh(\mathscr{b})} = \frac{\sinh(\gamma)}{\sinh(\mathscr{c})}.$\\

The coordinates of $I, Ge, N\!a, O, \tilde{G}, \tilde{K}$ can now be written as functions of the angles:\vspace*{-0.5mm}
\[
\begin{split}
I &= [\sinh(\alpha):\cdots:\cdots]_\Delta,\\
Ge &= [\tanh(\displaystyle\frac{\alpha}{2}):\cdots:\cdots]_\Delta,\;\;\;\;\;\;\;\;\;\;\;\;\;\;\;\hspace*{0mm}  N\!a = [\coth(\displaystyle\frac{\alpha}{2}):\cdots:\cdots]_\Delta,\\
O &= [\sinh(\alpha) \cosh(\alpha\textrm{-}\epsilon):\cdots:\cdots]_\Delta,\\ 
\tilde{K} &= [\sinh(\alpha) \sinh(\alpha\textrm{-}\epsilon):\cdots:\cdots]_\Delta,\;\;\;\;\hspace*{0mm} \tilde{G} = [\sinh(\alpha)/\!\sinh(\alpha\textrm{-}\epsilon):\cdots:\cdots]_\Delta.\\
\end{split}
\]\hspace*{-1mm}
The isotomic conjugate of $I$ is a point on the line $Ge\vee N\!a$.
The isogonal conjugate of $O$ is the point $H^- = [\sinh(\alpha)/\!\cosh(\alpha\textrm{-}\epsilon):\cdots:\cdots]_\Delta$. The cevian line $A\vee H^-$ is perpendicular to the sideline $B_G \vee C_G$ of the medial triangle.

 \vspace*{-0.5mm}
\subsubsection{Other triangle centers related to the circumcenters and the incenters}\label{subsubsec:Spieker}\hspace*{\fill}\vspace*{0.5mm} \\
\noindent$\bullet\;\;$The triple $I_1I_2I_3$ is perspective to $\Delta'$ at the \textit{Bevan point} \vspace*{1.5mm}\\ 
$\hspace*{7.5mm}Be = Be_0 = [\sqrt{\mathfrak{d}_{11}}(\sqrt{\mathfrak{d}_{11}}\sqrt{\mathfrak{d}_{22}}\sqrt{\mathfrak{d}_{33}}-\sqrt{\mathfrak{d}_{11}}\mathfrak{d}_{23}+
\sqrt{\mathfrak{d}_{22}}\mathfrak{d}_{31}+\sqrt{\mathfrak{d}_{33}}\mathfrak{d}_{12}):\cdots:\cdots]_\Delta$\vspace*{1mm}\\
$\hspace*{22.5mm}  = [\sinh(\alpha)(1+\cosh(\alpha)-\cosh(\beta)-\cosh(\gamma)):\cdots:\cdots]_\Delta$ \vspace*{1mm}\\
$\hspace*{22.5mm}    = [\sinh(\alpha)(\sinh^2(\frac12\alpha)-\sinh^2(\frac12\beta)-\sinh^2(\frac12\gamma)):\cdots:\cdots]_\Delta.$\vspace*{-2.5mm}\\

\noindent \hspace*{0mm}The euclidean limit of $Be$ is the point  $X_{40}$. $Be$ is the incenter of the extangents triangle of $\Delta_0$ (see \ref{subsubsec:extangents triangle}).  But, different from its euclidean limit, it is in general not the circumcenter of the excentral triangle $(I_1I_2I_3)_0$, neither a point on $I\vee O$.\\ The four lines $I_k\vee Be_k, 0\le k \le3$, meet at the point $O^+$.\vspace*{-2.5mm}\\

$\bullet\;\;$The four lines $O_k\vee I_k, k = 0,\cdots,3$, meet at the point\vspace*{-1mm}
\[ 
\begin{split}
&[\sqrt{\mathfrak{d}_{11}}(\sqrt{\mathfrak{d}_{22}}(\mathfrak{d}_{11}\mathfrak{d}_{33}+\mathfrak{d}_{31}(\mathfrak{d}_{23}-\mathfrak{d}_{33}))+\sqrt{\mathfrak{d}_{33}}(\mathfrak{d}_{11}\mathfrak{d}_{22}+\mathfrak{d}_{12}(\mathfrak{d}_{23}-\mathfrak{d}_{22}))) \\
& + \sqrt{\mathfrak{d}_{22}}\sqrt{\mathfrak{d}_{33}}\,\mathfrak{d}_{11}(\mathfrak{d}_{11}-\mathfrak{d}_{12}-\mathfrak{d}_{31}) + \mathfrak{d}_{11}(\mathfrak{d}_{22}\mathfrak{d}_{33}-\mathfrak{d}_{23}^{\;2}):\cdots:\cdots]_\Delta.
\end{split}
\]
This is a second point with euclidean limit $X_{40}$ and, in general, it also differs from all the circumcenters of excentral triangles of $\Delta$.\vspace*{-1.5mm}\\ 

\noindent$\bullet\;\;$The triples $\Delta_{[O]} = \Delta_G$ and $I_1 I_2  I_3$ are perspective at the \textit{Mittenpunkt} of $\Delta_0$,\vspace*{0.5mm}\\

\noindent $\hspace*{7.5mm}M\hspace*{-0.3mm}i = M\hspace*{-0.3mm}i_0 =
[\sqrt{\mathfrak{d}_{11}}({-}\sqrt{\mathfrak{d}_{11}}+\sqrt{\mathfrak{d}_{22}}+\sqrt{\mathfrak{d}_{33}}):\cdots:\cdots]_\Delta$\vspace*{0.5mm}\\
$\hspace*{23.5mm}= [\sinh(\alpha)({-}\sinh(\alpha)+\sinh(\beta)+\sinh(\gamma)):\cdots:\cdots]_\Delta,$\vspace*{0.5mm}\\
a point on the line $I\hspace*{-0.3mm}\vee K$. It has euclidean limit $X_9$.\vspace*{-2mm}\\

\noindent$\bullet\;\;$Define the points $P_a$ and $Q_a$ by {$P_a:=(I_0\vee M\!i_0)\wedge(I_1\vee M\!i_1)$} and 
$Q_a:={(I_2\vee M\!i_2)}\wedge{(I_3\vee M\!i_3)}$ and the points $P_b, Q_b, P_c, Q_c$ accordingly. The points $B,P_a,C,Q_a$ form a harmonic range. The triple $P_aP_bP_c$ is perspective to $\Delta$ at the perspector $[\displaystyle\frac{\sinh(\alpha)}{ \sinh(\beta){-}\sinh(\gamma)}:\cdots:\cdots]_\Delta$, a point with euclidean limit $X_{\!100}$.\vspace*{-1mm}\\

\noindent$\bullet\;\;$Vigara \cite{Vi} proved that the triples $\Delta_{[O]}$ and $\Delta'_{\;\,[I]}$ are perspective; the perspector he named \textit{pseudo- Spieker center}. But this point is in fact a good candidate for the Spieker center in elliptic and in hyperbolic geometry, as it is one of four radical centers of the three incircles $\mathscr{I}_k, k = 1,2,3$; and it is the one that lies inside $\Delta^0$. Its coordinates are\vspace*{0.5mm}\\
\noindent\hspace*{1mm}$S\!p = S\!p_0 = [(s_1{+}s_2{+}s_3)(s_1 c_2 c_3{+}s_2 c_3 c_1{+}s_3 c_1 c_2{+}s_1 s_1 s_3){\,+\,}(s_2 c_3{+}s_3 c_2)(s_2{+}s_3{{-}}2s_1)$\\
 $\hspace*{82mm}{\,+\,}s_1(2s_2{+}2s_3{{-}}s_1){:}\cdots:\cdots]_\Delta,$\vspace*{0.5mm} \\
with $s_1\!:=\sinh(\mathscr{a}),\, c_1\!:=\cosh(\mathscr{a}),\, s_2:=\sinh(\mathscr{b}),....$\vspace*{-2mm}\\

\noindent\hspace*{5mm}More general, each triple $\Delta_{[O_j]}$
is perspective to each triple $\Delta'_{\;\,[I_k]},\; j,k \in \{0,1,2,3\}$. By this, we get altogether 16 perspective centers. 
Let $P_{kj}$ be the perspector for $\Delta_{[0_k]}$ and $\Delta'_{\;\,[I_j]}$. GeoGebra-constructions indicate:\\
- The six points $P_{12}, P_{21}, P_{23}, P_{32}, P_{31}, P_{13}$ lie on a singular conic (union of two lines).\\
- Put $Q_1 = (P_{12}{\vee}P_{21}){\wedge}(P_{31}{\vee}P_{13}), Q_2 = (P_{23}{\vee}P_{32}){\wedge}(P_{12}{\vee}P_{21}), Q_3 = (P_{31}{\vee}P_{13}){\wedge}(P_{23}{\vee}P_{32});$\\
\noindent\hspace*{1.7mm}the triple $Q_1Q_2Q_3$ is perspective to $\Delta$ and to $\Delta'$ at points which lie on the line $O\vee I$.\vspace*{-1 mm}\\

\noindent$\bullet\;\;$The points $H^-, M\!i, Be, S\!p$ are collinear.\vspace*{-1 mm}\\

\noindent$\bullet\;\;$Define the point $S_1$ by $S_1 := (I_1\vee O^+)\wedge a$ and the points $S_2, S_3$ likewise. The triple $S_1S_2S_3$ is perspective to $\Delta$ at the point \vspace*{1 mm}\\
\centerline{$\displaystyle [\frac{\sinh(\alpha)}{\cosh(\beta)+\cosh(\gamma)}:\frac{\sinh(\beta)}{\cosh(\gamma)+\cosh(\alpha)}:\frac{\sinh(\gamma)}{\cosh(\alpha)+\cosh(\beta)}]_\Delta.$}\vspace*{1.5 mm}\\
We call this point \textit{pseudo- Schiffler point}; the euclidean limit of this point is the Schiffler point $X_{21}$, cf. \cite{EE}.\vspace*{1mm}\\

\noindent$\bullet\;\;$The antipedal triple $\Delta^{[O_0]}$ of $O_0$ is perspective to $O_1 O_2  O_3$ and to $\Delta'$ at $O_0$.\vspace*{-0.5mm}\\

Let $P$ be the point\vspace*{-1mm}\\
\centerline{$\displaystyle[\frac 1{(\sqrt{\mathfrak{d}_{22}}(\mathfrak{c}-\mathfrak{c}_{12})-\sqrt{\mathfrak{d}_{33}}(\mathfrak{c}-\mathfrak{c}_{31})}:\cdots:\cdots]_\Delta =
[\displaystyle\frac{\sinh(\displaystyle\frac\alpha 2)}{\sinh(\displaystyle\frac{\beta -\gamma} 2)}:\cdots:\cdots]_\Delta, \vspace*{0mm}$}\\
then:\hspace*{2mm} {$ (B^{[O_0]}\vee I_2) \wedge (C^{[O_0]}\vee I_3) = A_P,\; (C^{[O_0]}\vee I_3) \wedge (A^{[O_0]}\vee I_1)= B_P,$\;}  \vspace*{0.5mm}\\
\hspace*{52.3mm}and$\;\;\;(A^{[O_0]}\vee I_1) \wedge B^{[O_0]}\vee I_3) = C_P.$\vspace*{1mm}\\
$P$ is a point on the circumcircle $\mathscr{C}_0$, its euclidean limit is $X_{100}$.\vspace*{1.5mm}\\

\noindent$\bullet\;\;$The pedal triple $\Delta_{[I_0]}$ of $I_0$ is perspective to $I_1 I_2  I_3$ at\vspace*{1mm}\\

\noindent \hspace*{11.5mm}$[\sqrt{\mathfrak{d}_{11}}(-\,\mathfrak{d}_{23}\sqrt{\mathfrak{d}_{11}}+\mathfrak{d}_{31}\sqrt{\mathfrak{d}_{22}}+\mathfrak{d}_{12}\sqrt{\mathfrak{d}_{33}}-\sqrt{\mathfrak{d}_{11}}\sqrt{\mathfrak{d}_{22}}\sqrt{\mathfrak{d}_{33}}):\cdots:\cdots]_\Delta$ \\ \vspace*{1mm}
\noindent \hspace*{6.5mm}$=\, [\displaystyle\sinh(\alpha)\big(-\cosh(\frac\alpha 2)+\cosh(\frac\beta 2)+\cosh(\frac\gamma 2)\big):\cdots:\cdots]_\Delta.\vspace*{1mm}$\\
This point is also the orthocorrespondent of $I$. 
Its euclidean limit is $X_{57}.$\vspace*{0.5mm}\\
%The pedal triple $\Delta_{[I_0]}$ of $I_0$ is perspective to $\Delta'$ at $I_0$.\\
%This point has euclidean limit $X_{57}$.
%[w1 ((- d23 w1 + d31 w2 + d12 w3) - w1 w2 w3), w2 ((d23 w1 - d31 w2 + d12 w3) - w1 w2 w3), w3 ((d23 w1 + d31 w2 - d12 w3) - w1 w2 w3)]

\noindent$\bullet\;\;$The tripole $[1/(\mathfrak{c}_{31}{-}\mathfrak{c}_{12}):1/(\mathfrak{c}_{12}{-}\mathfrak{c}_{23}):1/(\mathfrak{c}_{23}{-}\mathfrak{c}_{31})]$ of the line $\tilde{K}\vee G^+$ is a point on the circumcircle and has euclidean limit $X_{99}$.\vspace*{0.5mm}\\

\noindent$\bullet\;\;$The circumcenters of the triangles $(O_0BC)_0, (O_0CA)_0, (O_0AB)_0$ form (in this order) a triple which is perspective to $\Delta$ at the \textit{Kosnita point} \vspace*{-1.5 mm}
\[
\begin{split}
&[1/\big((\mathfrak{d}_{12}+\mathfrak{d}_{22})(\mathfrak{d}_{31}+\mathfrak{d}_{33})- \mathfrak{d}_{23}(\mathfrak{d}_{11}+\mathfrak{d}_{23}+\mathfrak{d}_{31}+\mathfrak{d}_{12}-c\lambda)\big) :\cdots:\cdots]_\Delta\\
&\text{with\;\,} \lambda = \sqrt{|(1,1,1){\scriptscriptstyle{[\mathfrak{D}]}}(1,1,1)|}=\sqrt{|\mathfrak{d}_{11}+\mathfrak{d}_{22}+\mathfrak{d}_{33}+2(\mathfrak{d}_{12}+\mathfrak{d}_{23}+\mathfrak{d}_{31})|}.
\end{split}
\]\vspace*{-3 mm}
The euclidean limit of this point is $X_{54}$.\vspace*{1mm}\\

\noindent$\bullet\;\;$Put $P_1 := (O_2\vee C)\wedge (O_3\vee B), P_2 := (O_3\vee A)\wedge (O_1\vee C), P_3 := (O_1\vee B)\wedge (O_2\vee A)$. The triple $P_1P_2P_3$ is perspective to $\Delta$ at $O$ and to $\Delta'$ at \vspace*{1 mm}\\
$[(\mathfrak{c}_{23} + \mathfrak{c})(\mathfrak{c}_{23}^{\;2}(\mathfrak{c}_{23}+\mathfrak{c}_{31}+\mathfrak{c}_{12}+\mathfrak{c})-\mathfrak{c}_{23}(\mathfrak{c}_{31}+\mathfrak{c}_{12}+\mathfrak{c})^2 - (\mathfrak{c}_{31}+\mathfrak{c}_{12})((\mathfrak{c}_{31}-\mathfrak{c}_{12})^2+\mathfrak{c}_{31}+\mathfrak{c}_{12}+\mathfrak{c}) \,-\, \mathfrak{c}){:}\cdots{:}\cdots]_\Delta.$\vspace*{1 mm}\\
The euclidean limit of this point is the de{\!} Longchamps point $X_{20}$.\vspace*{-0.5mm}\\

\noindent$\bullet\;\;$The incenters of the triangles $(I_0BC)_0, (I_0CA)_0, (I_0AB)_0$ form (in this order) a triple perspective to $\Delta$. The perspector is the ${1}^{st}$ \textit{de\! Villiers point}\vspace*{2 mm}\\ 
$\displaystyle \hspace*{22mm}\frac{1}{\sqrt{2 \sqrt{\mathfrak{d}_{22}^{\;}}\sqrt{\mathfrak{d}_{33}^{\;}}(\sqrt{\mathfrak{d}_{22}^{\;}}\sqrt{\mathfrak{d}_{33}^{\;}}-\mathfrak{d}_{23}^{\;})} + \mathfrak{c}\sqrt{\mathfrak{d}_{22}^{\;}}\sqrt{\mathfrak{d}_{33}^{\;}}}:\cdots:\cdots]_\Delta$\\
$\displaystyle \hspace*{20mm}=\;\, \frac{\sinh(\alpha)}{2\cosh(\frac {1}{2}\alpha)+1}:\frac{\sinh(\beta)}{2\cosh(\frac {1}{2}\beta)+1}:\frac{\sinh(\gamma)}{2\cosh(\frac {1}{2}\gamma)+1}]_\Delta .$ \vspace*{1 mm}\\
It has euclidean limit $X_{1127}$.\vspace*{0.5 mm}\\
Experimental constructions using GeoGebra indicate that there also exist elliptic and hyperbolic analogues of the $2^{nd}$ de\! Villiers point $X_{1128}$ and of the three Stevanovic points $X_{1130}$, $X_{1488}$  and $X_{1489}$.\vspace*{-1mm}\\
\subsubsection{Kimberlings "Four-Triangle Problem"} Let $T, \tilde{T}$ be triangle centers of $\Delta_0$, and let $\tilde{T}_1,\tilde{T}_2,\tilde{T}_3$ be the $\tilde{T}$-centers of the triangulation triangles $(AT\!B)_0,(BTC)_0,(CT\!A)_0$, respectively. If the triple $\tilde{T}_1\tilde{T}_2\tilde{T}_3$ is perspective to $\Delta_0$, we will say that $T\# \tilde{T}$ is well-defined, and $T\# \tilde{T}$ will stand, in this case, for the perspector. There is a problem, which are the centers $T$ and $\tilde{T}$ such that $T\# \tilde{T}$ is well-defined. Kimberling posed this problem (the "Four-Triangle Problem") in \cite{Ki} for the special case $T = \tilde{T}$, and as far as I know, this problem is still open.\\
It was shown above that $O\# O$ is the Kosnita point and that $I\# I$ is the de\! Villiers point. Experiments with GeoGebra indicate that $T\# T$ is well-defined for the absolute centers $T = H, G^+, O^+, N^+,$ $L, H^{\star}$ and that $L\# L = H$, but for the absolute center $P$ on the orthoaxis (see \ref{subsubsec:orthoaxis}) it is not well-defined. If $T{\,\in\,}\{\tilde{K}, \tilde{G}\}$, then $T\# T$ is well-defined, whereas for $T{\,\in\,}\{Ge, N\!a\}$ it is not.\vspace*{1mm}\\ \newpage

\subsection{Circles, radical centers and centers of similitude}\hspace*{\fill} \\ 

\noindent The following two theorems together with their proofs were presented by Ungar \cite{U2} for the proper hyperbolic case.\hspace*{-2.5mm}
\subsubsection{The Inscribed Angle Theorem} 
Let $\mu$ be the measure of the angle $\angle_+(BOC)$, then\vspace*{-1mm}
\[ 
\begin{split}
&\sinh(\alpha) = \frac{\sqrt{|(1,1,1){\,\scriptscriptstyle{[\mathfrak{D}]}\,}(1,1,1)|}}{\cosh(\frac{\mathscr{b}}{2})\cosh(\frac{\mathscr{c}}{2})} \sinh(\mu/2) \;\;\textrm{and}\;\\
&\sinh(\alpha - \frac{1}{2} \textrm{area}(\Delta_0)) = \sinh(\frac{1}{2}(\mu - \textrm{area}((BOC)_0))).\vspace*{-1mm}
\end{split}
\] 
A special case: If $B, C, O$ are collinear, then 
$\alpha - \frac{1}{2} \textrm{area}(\Delta_0) = \frac{1}{2}\pi \; \textrm{and} \; \alpha = \beta + \gamma.$
\subsubsection{Tangent-Secant Theorem}The tangent at $A$ of the circumcircle $\mathscr{C}_0$ of $\Delta_0$ meets the line $a=B\vee C$ at the point $P = [0:\mathfrak{c}-\mathfrak{c}_{13}:\mathfrak{c}_{12}-\mathfrak{c}]_\Delta = [0:\sinh^2(\frac{\mathscr{b}}{2}):-\sinh^2(\frac{\mathscr{c}}{2})]_\Delta$,  the harmonic conjugate of $A_{\tilde{K}}$ with respect to $\{B,C\}$, and \vspace*{-1mm}
\[\sinh^2(\mu([P,A]_+)) \cosh^2(\frac{\mathscr{a}}{2}) = \sinh(\mu([P,B]_+))\sinh(\mu([P,C]_+)).\]
\begin{figure}[h]
\includegraphics[height=9cm]{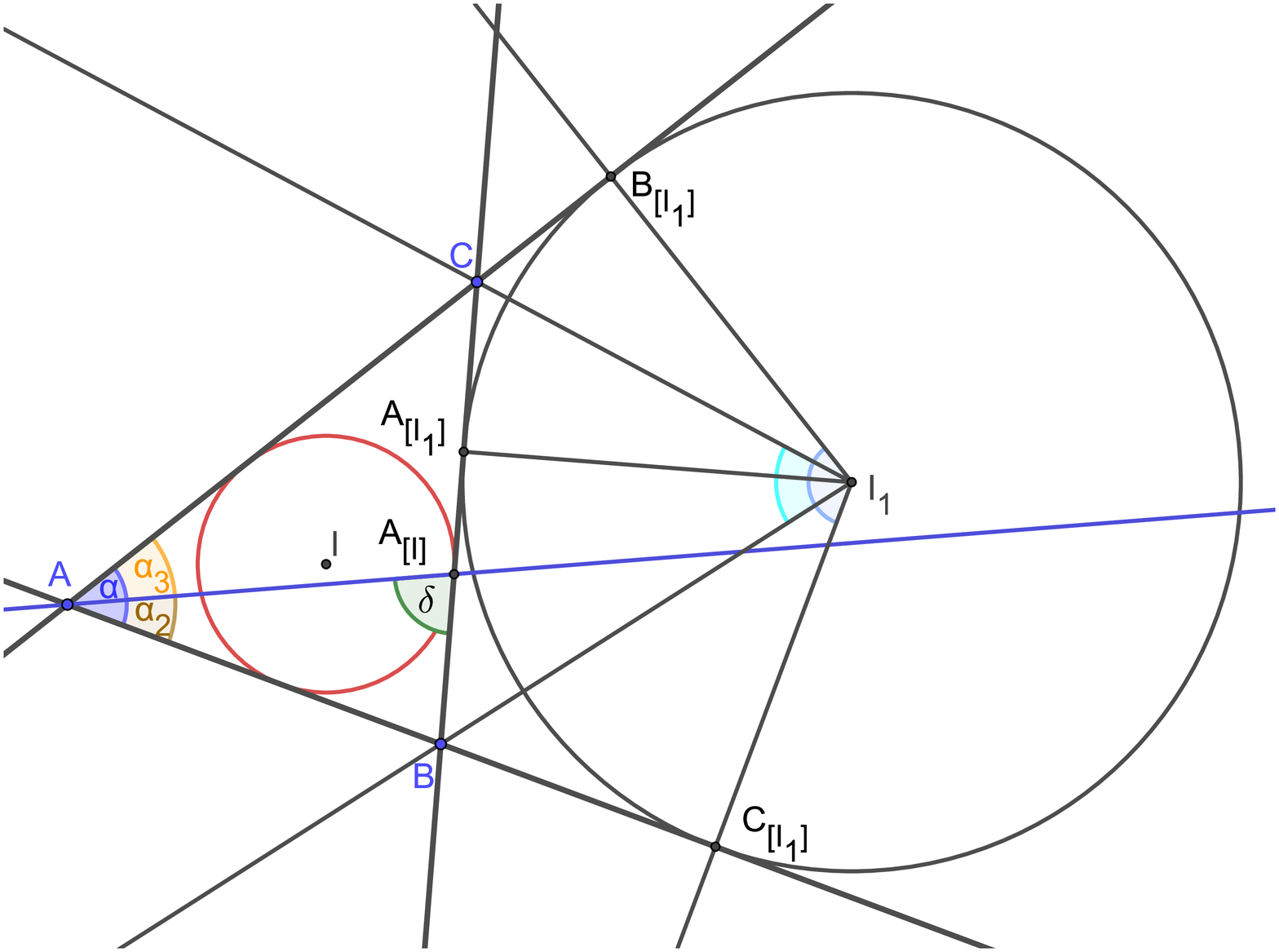}
\caption{
%For example, the pole $A'$ of the line $B \vee C$ with respect to $\mathcal{C}$ can be obtained as follows: Construct the pole of $B \vee C$ with respect to $\tilde{\mathcal{C}}$, then its mirror image in the center $(1:0:0)$ of $\tilde{\mathcal{C}}$ is $A'$.
}\label{fig:eucl}
\end{figure}\vspace*{-0.5mm}

\subsubsection{}It is easy to dualize these two theorems. We present euclidean limit versions of these duali\-zations, see Figure \ref{fig:eucl}:\vspace*{2mm}\\
\centerline{$2\,\textrm{measure}(\angle_+(B I_i C)) = \textrm{measure}(\angle_+(B_{[I_i]} I_i C_{[I_i]})) ,\; i = 0, 1.$}\vspace*{-1mm}\\

Between the angles $\alpha_2, \alpha_3, \delta$ in Figure \ref{fig:eucl}, the following relationship applies:\vspace*{-0.5mm}
\[ \sin^2(\delta) = \frac {2\sin(\alpha_2)\sin(\alpha_3)}{1-\cos(\alpha_2{+}\alpha_3)} = \frac {\sin(\alpha_2)\sin(\alpha_3)}{\sin^2(\frac{1}{2}\alpha)} \vspace*{1.5mm}.\]

A  spherical version of the following theorem together with a proof can be found in \cite{TL} ch. IX.
\subsubsection{A second Tangent-Secant Theorem} Let $\mathcal{S}$ be a nondegenerate circle with center $M$ and radius $r$, and let $P$ be an anisotropic point. We introduce two sets of lines:\vspace*{-0.5mm}\\
\hspace*{1.5mm} $\mathscr{S}_{\!P} := \{ l\,|\,l\,\textrm{is line through P\,},
\textrm{intersecting the circle\,} \mathcal{S}\, \textrm{at two distinct anisotropic points}\,\}$\vspace*{0.5mm}\\
$\hspace*{1.5mm}\,\mathscr{T}_{\!P} := \{ l\,|\,l\,\textrm{is line through P touching } \mathcal{S} \textrm{\,at an anisotropic point\,}\}.$\vspace*{0.5mm}\\
Since there are at most two isotropic points on $\mathcal{S}$, there are infinitely many lines in $\mathscr{S}_{\!P}$, and if $P$ is a point on or outside the circle, then there is at least one line in $\mathscr{T}_{\!P}$.\\
For each line $l$ in $\mathscr{S}_{\!P}$ we define a number $\mathscr{p}(P,l)$ as follows: If $Q$ and $R$ are points of intersection of $l$ and $\mathcal{S}$, then\\
$\hspace*{35mm}\mathscr{p}(P,l) := \tanh\big(\frac{1}{2}\mu([P,Q]_+)\big)\tanh\big(\frac{1}{2}\mu([P,R]_+)\big)$\\
$\hspace*{48.2mm}= \tanh\big(\frac{1}{2}\mu([P,Q]_-)\big)\tanh\big(\frac{1}{2}\mu([P,R]_-)\big)$.\vspace*{0.5mm}\\
For a tangent $l \in \mathscr{T}_{\!P}$ with touchpoint $T$ we put\\
$\hspace*{35mm}\mathscr{p}(P,l) := \big(\tanh(\frac{1}{2}d(P,T))\big)^2.$\vspace*{0mm}\\
The second Tangent-Secant Theorem states that for lines $l_1, l_2 \in \mathscr{S}_{\!P}\cup \mathscr{T}_{\!P}$:\\
\centerline{$\mathscr{p}(P,l_1) = \mathscr{p}(P,l_2 
).$\;\;\;\;\;\;\;\;\;\;\;\;\;\;\;}\vspace*{0.5mm}\\
Thus, there exists a number $\mathscr{p}(P,\mathcal{S})$ such that $\mathscr{p}(P,\mathcal{S}) = \mathscr{p}(P,l)$ for all $l\in \mathscr{S}_{\!P}\cup \mathscr{T}_{\!P}$. This number is called the \textit{power of the point} $P$ \textit{with respect to the circle} $\mathcal{S}$.\\
This power can be calculated by \vspace*{0.5mm}\\
\centerline{$\mathscr{p}(P,\mathcal{S}) =\tanh\big(\frac{1}{2}(d(P,M)+r)\big)\tanh\big(\frac{1}{2}(d(P,M)-r)\big)$}\vspace*{0.5mm}\\
$\displaystyle\hspace*{35mm}=\frac{\cosh(d(P,M))-\cosh(r)}{\cosh(d(P,M))+\cosh(r)}$.\vspace*{0.5mm}
\subsubsection{Radical lines of two circles} Let $\mathcal{S}_1, \mathcal{S}_2$ be two nondegenerate circles  with centers $M_1, M_2$, $d(M_1, M_2) \ne 0$, and radii $r_1, r_2$. 
We want to find out which anisotropic points have the same power with respect to both circles.\\ First, we define a "modified power" of an anisotropic point $P$ with respect to a circle $\mathcal{S}$ having center $M$ and radius $r$ by $\tilde{\mathscr{p}}(P,\mathcal{S}) := \cosh(d(P,M))/\!\cosh(r)$. In comparison to ${\mathscr{p}}$, this modified power is easier to handle, on the other hand, an anisotropic point $P$ has the same power with respect to $\mathcal{S}_1$ and $\mathcal{S}_2$ precisely when $\tilde{\mathscr{p}}(P,\mathcal{S}_1) = \tilde{\mathscr{p}}(P,\mathcal{S}_2)$.\\
We put $\mathscr{R} = \mathscr{R}(\mathcal{S}_1, \mathcal{S}_2):= \{P\,|\,P\, \textrm{is anisotropic and } \tilde{\mathscr{p}}(P,\mathcal{S}_1) = \tilde{\mathscr{p}}(P,\mathcal{S}_2)\}$. One recognizes immediately that the point  $M':= (M_1{\vee}M_2)^\delta$ belongs $\mathscr{R}$,  as well as all anisotropic points of intersection of the two circles. 
Moreover, whenever $P$ is a point in $\mathscr{R}$, different from $M'$, then every other anisotropic point $Q$ on $M'{\vee}P$ is also a point in $\mathscr{R}$, because:\vspace*{0.5mm}\\
\centerline{$\displaystyle\tilde{\mathscr{p}}(Q,\mathcal{S}_1){\,=\,}\tilde{\mathscr{p}}(P,\mathcal{S}_1)\frac{\cosh(\frac{1}{2}\pi i{-}d(Q,M'))}{\cosh(\frac{1}{2}\pi i{-}d(P,M'))}{\,=\,}\tilde{\mathscr{p}}(P,\mathcal{S}_2)\frac{\cosh(\frac{1}{2}\pi i{-}d(Q,M'))}{\cosh(\frac{1}{2}\pi i{-}d(P,M'))}{\,=\,}\tilde{\mathscr{p}}(Q,\mathcal{S}_2)$.}\vspace*{0.5mm}\hspace*{\fill}\\ \newpage 
There are exactly two points in $\mathscr{R}\cap (M_1{\,\vee\,}M_2)$,\vspace*{0.5mm}\\
$\hspace*{11mm}P_1{\,=\,}\big(\!\cosh(d)\cosh(r_2){-}\cosh(r_1)\big)M_1 + \big(\!\cosh(d)\cosh(r_1){-}\cosh(r_2)\big)M_2$\\
\hspace*{5mm}and $P_2{\,=\,}\big(\!\cosh(d)\cosh(r_2){+}\cosh(r_1)\big)M_1 - \big(\!\cosh(d)\cosh(r_1){+}\cosh(r_2)\big)M_2$,\vspace*{0.5mm}\\
\hspace*{5mm}with $d = d(M_1,M_2).$\vspace*{1mm}\\
Thus, $\mathscr{R}$ is the union of the lines $l_1 = P_1\vee M'$ and $l_2 = P_2\vee M'$. These lines are called \textit{radical lines} of the circles $\mathcal{S}_1$ and $\mathcal{S}_2$. The lines $l_1^{\,}, M_1^{\,\delta}, l_2^{\,}, M_2^{\,\delta}$ form a harmonic pencil.
\subsubsection{Radical centers of three circles}\label{subsubsec:Radical centers} We draw a circle around each vertex of triangle $\Delta_0$, around $A$ a circle $\mathcal{S}_1$ with radius $r_1$, etc. Then there are exactly four \textit{radical centers}, points of equal powers with respect to the three circles. One of these points, $R_0$, is a point inside the triangle $\Delta^0$, with coordinates \vspace*{0.5mm}\\
\centerline{$R_0 = R_0(r_1,r_2,r_3) = [\cosh(r_1) \mathfrak{d}_{11}+\cosh(r_2) \mathfrak{d}_{12}+\cosh(r_3) \mathfrak{d}_{13}:\cdots:\cdots]_\Delta$.}\vspace*{0.5mm}
The other three points $R_1, R_2, R_3$ form the anticevian triple of $R_0$ with respect to $\Delta'$. $R_1$ has coordinates:\vspace*{0.5mm}\\
$\;\;R_1 = [\cosh(r_1) \mathfrak{d}_{11}{-}\cosh(r_2) \mathfrak{d}_{12}{-}\cosh(r_3) \mathfrak{d}_{13}{\,:\,}\\
\hspace*{2.5mm}{-}\cosh(r_1)\mathfrak{d}_{21}{+}\cosh(r_2)\mathfrak{d}_{22}{-}\cosh(r_3)\mathfrak{d}_{23}{\,:\,}{-}\cosh(r_1)\mathfrak{d}_{31}{-}\cosh(r_2)\mathfrak{d}_{32}{+}\cosh(r_3)\mathfrak{d}_{33}]_\Delta$.\vspace*{-2.5mm}\\

A \textit{radical circle} $\tilde{\mathcal{S}}_k$  can be drawn around each point $R_k$; this circle meets the circles $\mathcal{S}_1, \mathcal{S}_2, \mathcal{S}_3$ orthogonally. The radius of $\tilde{\mathcal{S}}_k$ is $\tilde{r}_k = \tilde{\mathscr{p}}(R_k,\mathcal{S}_1) = \tilde{\mathscr{p}}(R_k,\mathcal{S}_2) = \tilde{\mathscr{p}}(R_k,\mathcal{S}_3)$.\vspace*{-2.5mm}\\

By taking special values for $r_1, r_2, r_3$, we can find triangle centers of $\Delta_0$.\vspace*{0.5mm}\\ First, independently of the choice of the radii, we have $\lim\limits_{t \to 0} R_0(t\,r_1,t\,r_2,t\,r_3) = O_0.$ We also get the circumcenter $O_0$  as a result for $R_0$ when we take $r_1=r_2=r_3$.\\
If $r_1=\frac 12(-\mathscr{a}+\mathscr{b}+\mathscr{c}), r_2=\frac 12(\mathscr{a}-\mathscr{b}+\mathscr{c}), r_3=\frac 12(\mathscr{a}+\mathscr{b}-\mathscr{c})$, then $R_0 = I_0$.\\
When we choose $r_1=\mathscr{a}, r_2=\mathscr{b}, r_3=\mathscr{c}$, then $R_0$ = $L$ (de{\!} Longchamps point). \\
\textit{Remark}: In the euclidean plane, when taking radii $r_1=t\mathscr{a}, r_2=t\mathscr{b}, r_3=t\mathscr{c}, t\in \mathbb{R}^{\ge 0}$, one gets points
$R_0(t)=[ t^2(\mathscr{a}^2 (-2\mathscr{a}^2{+}\mathscr{b}^2{+}\mathscr{c}^2)+(\mathscr{b}^2 - \mathscr{c}^2)^2)+\mathscr{a}^2 (\mathscr{a}^2{-}\mathscr{b}^2{-}\mathscr{c}^2):\cdots:\cdots]_\Delta$, all lying on the Euler-line.
%I checked if euclidean limits of centers constructed that way are listed in ETC \cite{ETC}, but I could just find one: For $r_1=\mathscr{a}, r_2=\mathscr{b}, r_3=\mathscr{c}$, the euclidean limit of $R_1$ is the de{\!} Longchamps point $X_{20}.$ 
\subsubsection{Centers of similitude of two circles} Given two nondegenerate circles $\mathcal{S}_1, \mathcal{S}_2$ with centers $M_1, M_2$, $d(M_1, M_2) \ne 0$, and radii $r_1, r_2$, then there exist two points $L_1, L_2$ which are the duals of the two radical lines of the duals of $\mathcal{S}_1$ and $\mathcal{S}_2$. These points are called \textit{centers of similitude} of $\mathcal{S}_1, \mathcal{S}_2$, cf. \cite{TL}.\\
These two centers lie on the line $M_1\vee M_2$, and $M_1, L_1, M_2, L_2$ form a harmonic range. If the two circles have common tangents, then each of these tangents passes either through $L_1$ or through $L_2$.
\subsubsection{Dualizing \ref{subsubsec:Radical centers}}\label{subsubsec:similitude} Given circles $\mathcal{S}_1,\mathcal{S}_2,\mathcal{S}_3$ with centers $A,B,C$ and radii $r_1, r_2,r_3$ (respectively), then there exist six centers of similitude of these circles, taken in pairs. Three of these are the vertices of the cevian triangle of the point $\,T = T(r_1,r_2,r_3) = [1/\!\sinh(r_1):1/\!\sinh(r_2):1/\!\sinh(r_3)]_\Delta$, while the other three centers lie on the tripolar of $\,T$. \\For $r_1=r_2=r_3$, we have $\,T=G_0$. When we choose for the radii $r_1=\mathscr{a}, r_2=\mathscr{b}, r_3=\mathscr{c}$, then the point $\,T$ is the isotomic conjugate of the incenter $I_0$, and for $r_1=\frac 12({-}\mathscr{a}{+}\mathscr{b}{+}\mathscr{c}), r_2=\frac 12(\mathscr{a}{-}\mathscr{b}{+}\mathscr{c}), r_3=\frac 12(\mathscr{a}{+}\mathscr{b}{-}\mathscr{c})$ we get $\,T = Ge_0$.\hspace*{-1.5mm}

\begin{figure}[!thpb]
\includegraphics[height=10cm]{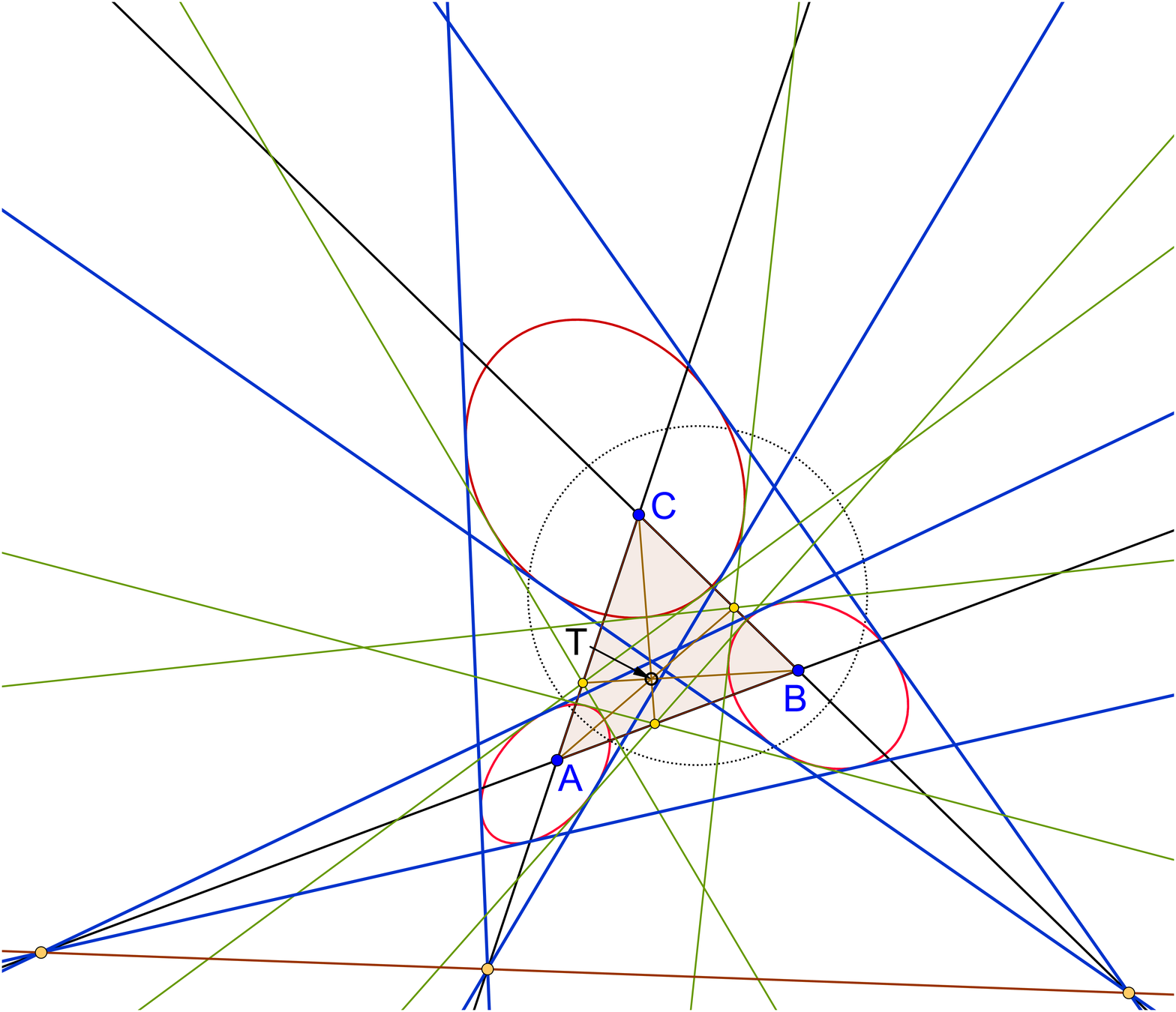}
\caption{Illustration of \ref{subsubsec:similitude}.
\newline
This Figure and the following with the exception of Figure 12 show the situation in the elliptic plane. An indication is the grey dotted circle. Since the absolute conic $\mathcal{C}_{abs}$ has no real points, this circle  $\tilde{\mathcal{C}} = \{(x_0{:}x_1{:}x_2) |\, \rho x_0^{\;2} = x_1^{\;2}+x_2^{\;2}\}$ serves as a substitute for constructions.}\vspace*{-2mm}
\end{figure}
\subsubsection{The excentral triangle and the extangents triangle}\label{subsubsec:extangents triangle} The triangle $(I_1I_2I_3)_0$ is called the excentral triangle of $\Delta_0$. The radical centers of the three excircles are the Spieker center (see \ref{subsubsec:Spieker}) together with its harmonic associates with respect to the dual $I_1'I_2'I_3'$ of the triple $I_1I_2I_3$.\\
The triple $I_1I_2I_3$ is perspective to the orthic triple $\Delta_{[H]}$ at the perspector\vspace*{0.5mm}\\
\centerline{[$\sqrt{\mathfrak{d}_{11}}(\mathfrak{d}_{23}-\mathfrak{d}_{31}-\mathfrak{d}_{12}):\cdots:\cdots]_\Delta=[\sinh(\alpha)(\cosh(\alpha)-\cosh(\beta)-\cosh(\gamma)):\cdots:\cdots]_\Delta,$}\vspace*{0.5mm}
a point with euclidean limit $X_{46}.$

We introduce three points\vspace*{0.5mm}\\
$\hspace*{0.5mm}E_1\!\! := [-\mathfrak{d}_{11} (\sqrt{\mathfrak{d}_{22}} \sqrt{\mathfrak{d}_{33}}{+}\mathfrak{d}_{23}){:}\sqrt{\mathfrak{d}_{22}} (\sqrt{\mathfrak{d}_{11}} \mathfrak{d}_{23}{+}\sqrt{\mathfrak{d}_{33}}\mathfrak{d}_{12} ){:}\sqrt{\mathfrak{d}_{33}} (\sqrt{\mathfrak{d}_{11}} \mathfrak{d}_{23}{+}\sqrt{\mathfrak{d}_{22}} \mathfrak{d}_{31})]_\Delta,$\vspace*{0.5mm}\\
$\hspace*{0.5mm}\;\;\;\;\;\; = [-\sinh(\alpha)(1{+}\cosh(\alpha))\,{:}\sinh(\beta)(\cosh(\alpha){+}\cosh(\beta))\,{:}\sinh(\gamma)(\cosh(\alpha){+}\cosh(\gamma))]_\Delta,$\vspace*{0.5mm}\\
$\hspace*{0.5mm}E_2\! := [\sqrt{\mathfrak{d}_{11}} (\sqrt{\mathfrak{d}_{33}} \mathfrak{d}_{12}{+}\sqrt{\mathfrak{d}_{22}}\mathfrak{d}_{31}){:}{-\mathfrak{d}_{22}} (\sqrt{\mathfrak{d}_{33}} \sqrt{\mathfrak{d}_{11}}{+}\mathfrak{d}_{31} ){:}\sqrt{\mathfrak{d}_{33}} (\sqrt{\mathfrak{d}_{11}} \mathfrak{d}_{23}{+}\sqrt{\mathfrak{d}_{22}} \mathfrak{d}_{31})]_\Delta,$\vspace*{1mm}\\
$\hspace*{0.5mm}E_3\! := [\sqrt{\mathfrak{d}_{11}} (\sqrt{\mathfrak{d}_{33}} \mathfrak{d}_{12}{+}\sqrt{\mathfrak{d}_{22}}\mathfrak{d}_{31}){:}\sqrt{\mathfrak{d}_{22}} (\sqrt{\mathfrak{d}_{11}} \mathfrak{d}_{23}{+}\sqrt{\mathfrak{d}_{33}}\mathfrak{d}_{12} ){:}-\mathfrak{d}_{33} (\sqrt{\mathfrak{d}_{11}} \sqrt{\mathfrak{d}_{22}}{+}\mathfrak{d}_{12})]_\Delta.$\vspace*{1mm}\\
These points are the vertices of the \textit{extangents triangle} of $\Delta_0$, a triangle with following properties:\\
$\bullet$ The sideline $E_j\vee E_k$ of this triangle is a tangent of the excircles $\mathscr{I}_{\!j}$ and $\mathscr{I}_k, 1\leq j < k\leq 3.$\\
$\bullet$ It has the Bevan point $Be$ (s.\ref{subsubsec:Spieker} ) as its incenter.\vspace*{0.5mm}\\
The triple $E_1E_2E_3$ is perspective to $\Delta$ at perspector \vspace*{0.5mm}\\$[\sqrt{\mathfrak{d}_{11}} (\sqrt{\mathfrak{d}_{22}}\mathfrak{d}_{31}+\sqrt{\mathfrak{d}_{33}} \mathfrak{d}_{12}){\,:\,}\sqrt{\mathfrak{d}_{22}} (\sqrt{\mathfrak{d}_{33}}\mathfrak{d}_{12}+\sqrt{\mathfrak{d}_{11}} \mathfrak{d}_{23}   ){\,:\,}\sqrt{\mathfrak{d}_{33}} (\sqrt{\mathfrak{d}_{11}} \mathfrak{d}_{23} + \sqrt{\mathfrak{d}_{22}} \mathfrak{d}_{31})]_\Delta$\vspace*{0.5mm}\\
$= [\sinh(\alpha)(\cosh(\beta)+\cosh(\gamma)):\sinh(\beta)(\cosh(\gamma)+\cosh(\alpha)):\sinh(\gamma)(\cosh(\alpha)+\cosh(\beta))]_\Delta.$\vspace*{0.5mm}\\

This point is the isogonal conjugate of the pseudo- Schiffler point (see \ref{subsubsec:Spieker}) and has euclidean limit $X_{65}$, but in contrast to the euclidean case this perspector differs from the orthocenter of the intouch triangle.\\
The triple $E_1E_2E_3$ is also perspective to the orthic triple $\Delta_{[H]}$, the perspector is the \textit{Clawson point} \vspace*{0.5mm}\\
$[\sqrt{\mathfrak{d}_{11}}\mathfrak{d}_{31}\mathfrak{d}_{12}(\sqrt{\mathfrak{d}_{22}}\sqrt{\mathfrak{d}_{33}}+\mathfrak{d}_{23})(\sqrt{\mathfrak{d}_{11}}\sqrt{\mathfrak{d}_{22}}\sqrt{\mathfrak{d}_{33}}-\sqrt{\mathfrak{d}_{11}}\mathfrak{d}_{23}+\sqrt{\mathfrak{d}_{22}}\mathfrak{d}_{31}+\sqrt{\mathfrak{d}_{33}}\mathfrak{d}_{12}):\cdots:\cdots]_\Delta = [\tanh(\alpha)\cosh^2(\frac {\alpha}{2})(\cosh^2(\frac {\alpha}{2})-\cosh^2(\frac {\beta}{2})-\cosh^2(\frac {\gamma}{2}):\cdots:\cdots]_\Delta$\\ with euclidean limit $X_{19}$.

\subsection{Orthology, pedal-cevian points and cevian-pedal points}
\noindent\vspace*{-0.5mm}\subsubsection{Orthologic triples} A point-triple $PQR$, $P\ne A', Q\ne B',R\ne C'$,  is \textit{orthologic} to the triple $\Delta=ABC$ if the lines $P\vee A', Q\vee B', R\vee C'$ concur at some point $S$, which is then called the \textit{center} of this orthology.\\
If $PQR$ is orthologic to $\Delta$, then $\Delta$ is orthologic to $PQR$; the lines $\textrm{perp}(Q\vee R,A)$, $\textrm{perp}(R\vee P,B)$, $\textrm{perp}(P\vee Q,C) $
meet at some point $\,T$.\\
\textit{Outline of a proof}: If $PQR$ is orthologic to $\Delta$ with center $S=[s_1{:}s_2{:}s_2]_\Delta$, then there exist real numbers $x,y,z$ such that 
\[
\begin{split}
P &= [s_1+x\,d_{11}:s_2+x\,d_{12}:s_3+x\,d_{31}]_\Delta,\\
Q &= [s_1+y\,d_{12}:s_2+y\,d_{22}:s_3+y\,d_{23}]_\Delta,\\
R &= [s_1+z\,d_{31}:s_2+z\,d_{23}:s_3+z\,d_{33}]_\Delta.\\
\end{split}
\]
Define vectors $\boldsymbol{p},\boldsymbol{q},\boldsymbol{r}$ by $\boldsymbol{p}=(s_1{+}x\,d_{11},s_2{+}x\,d_{12},s_3{+}x\,d_{31}), \boldsymbol{q} = (s_1{+}x\,d_{12},s_2{+}x\,d_{22},$ $s_3{+}x\,d_{23})$, $\boldsymbol{r} = (s_1{+}z\,d_{31},s_2{+}z\,d_{23},s_3{+}z\,d_{33})$, then $((\boldsymbol{q}\times\boldsymbol{r})D)\times A$, $((\boldsymbol{r}\times\boldsymbol{p})D)\times B$, ${((\boldsymbol{p}{\times}\boldsymbol{q})D)\times C}$ form a linear dependent system (use CAS to check). The lines $\textrm{perp}(Q{\vee}R,A)$, $\textrm{perp}(R{\vee}P,B), \textrm{perp}(P{\vee}Q,C) $ meet at the point $T = [1/x,1/y,1/z]_\Delta.\;\;\;\;\Box$\vspace*{1.5mm}\\
\textit{Remark:} In euclidean geometry, the coordinates of $S$ with respect to the triple $PQR$ are the same as the coordinates of $\,T$ with respect to $\Delta$ (cf.\cite{DD}). This is not true in elliptic/hyperbolic geometry. On the other hand, still applies the\vspace*{1.5mm}\\
\noindent\textit{Addition}: If $P{\,=\,}A_{[S]},\, Q{\,=\,}B_{[S]},\,R{\,=\,}C_{[S]}$ form the pedal triple of $S$, then $T$ is the isogonal conjugate of $S$.\vspace*{-0.5mm}
\subsubsection{Pedal-cevian points and the Darboux cubic} A point $P$ is a \textit{pedal-cevian point} of $\Delta$ if its pedal triple $\Delta_{[P]}$ is perspective to $\Delta$; the perspector we call \textit{cevian companion} of $P$. A point $P = [p_1{:}p_2{:}p_3]_\Delta$ is a pedal-cevian point precisely when its coordinates satisfy the cubic equation\vspace*{2mm}\\
\centerline{$
\sum\limits_{j=1,2,3} (\mathfrak{d}_{jj} \mathfrak{d}_{j+1,i+2} + \mathfrak{d}_{j,j+1} \mathfrak{d}_{j+2,j}) p_j (\mathfrak{d}_{j+1,j+1}^{\; } p_{j+2}^{\;2}-\mathfrak{d}_{j+2,j+2}^{\; } p_{i+1}^{\;2}) \;=\;0.$}\\

%                             2         2 
%p1 (d11 d23 + d12 d31) (d22 p3  - d33 p2 )
\noindent As in the euclidean case, this cubic is a self-isogonal cubic with pivot point $L$.
On this cubic - we call it \textit{Darboux cubic} as its euclidean limit - lie the points $O_k, I_k (k = 0,1,2,3), O^+\!, H$, $L$, $Be$ and their isogonal conjugates. The cevian companions of $O_k, I_k$, $O^+\!, H, Be, L$ are $G_k, Ge_k$, $G^+\!,$ $H$, $N\!a$ and the isotomic conjugate of $H$, respectively. 
The points $A', B', C'$ are also lying on the Darboux cubic. Their cevian companions are ${[-\,\mathfrak{c}{\,:\,}1/\mathfrak{c}_{12}{\,:\,}1/\mathfrak{c}_{31}]_\Delta}$, $[\,1/\mathfrak{c}_{12}{\,:\,}-\,\mathfrak{c}{\,:\,}1/\mathfrak{c}_{32}]_\Delta,$ $[\,1/\mathfrak{c}_{31}{\,:\,}1/\mathfrak{c}_{23}{\,:\,}-\,\mathfrak{c}]_\Delta$, respectively. 
These three points form a triple which is perspective to $\Delta$ at $G^+$. \vspace*{-1mm}
\subsubsection{Cevian-pedal points and the Lucas cubic}
A point $P$ is a \textit{cevian-pedal point} of $\Delta$ if its cevian triple $\Delta_Q$ is perspective to $\Delta'$; the perspector we call the \textit{pedal companion} of $Q$. The cevian-pedal points $P = [p_1{:}p_2{:}p_3]_\Delta$ form a cubic with the equation \vspace*{2mm}\\
\centerline{$
\sum\limits_{j=1,2,3} \mathfrak{d}_{j+1,j+2} p_j (p_{j+2}^{\;2}-p_{j+1}^{\;2}) \;=\;0.$
} \\

\noindent This cubic, we call it \textit{Lucas cubic}, is a self-isotomic pivotal cubic; the pivot is the isotomic conjugate of $H$.
On this cubic lie the points $G_k, Ge_k, N\!a_k (k = 0,1,2,3),G^+, H, L$ and their isotomic conjugates. The pedal companion of $L$ is \vspace*{1.5 mm}\\
\centerline{$[\mathfrak{d}_{11}(2\mathfrak{d}_{22}\mathfrak{d}_{33}\mathfrak{d}_{31}\mathfrak{d}_{12}+\mathfrak{d}_{23}(-\mathfrak{d}_{11}^{\; }\mathfrak{d}_{23}^{\;2}+\mathfrak{d}_{22}^{\; }\mathfrak{d}_{31}^{\;2}+\mathfrak{d}_{33}^{\; }\mathfrak{d}_{12}^{\;2}+\mathfrak{d}_{11}\mathfrak{d}_{22}\mathfrak{d}_{33})):\cdots:\cdots]_\Delta$}\vspace*{1.5 mm}\\This is another point on the Darboux conic; it has euclidean limit $X_{1498}$. \vspace*{-0.5mm}
\subsubsection{The Darboux cubic and the Lucas conic of $\Delta'$} The pedal-cevian points of $\Delta$ are the cevian-pedal points of $\Delta'$ and vice versa. Therefore, the Darboux cubic and the Lucas cubic of $\Delta$ are the Lucas cubic and the Darboux cubic of $\Delta'$, respectively. \vspace*{-0.5mm}
%- d11 (- d11 d23^3 + 2 d12 d22 d31 d33 + d23 (d11 d22 d33 + d12^2 d33 + d22 d31^2))
\subsubsection{The Thomson cubic} In euclidean geometry, the \textit{Thomson cubic} is the locus of perspectors of circumconics such that the normals at the vertices $A, B, C$ meet at one point. An equation of the elliptic/hyperbolic analog is
\[
\sum\limits_{j=1,2,3} (\mathfrak{d}_{j,j}\mathfrak{d}_{j+1,j+2}-\mathfrak{d}_{j,j+1}\mathfrak{d}_{j+2,j}) p_j (\mathfrak{d}_{j+2,j+2}p_{j+1}^{\;2}-\mathfrak{d}_{j+1,j+1}p_{j+2}^{\;2}) \;=\;0.\;\;\;\;\;\;    (\star)
\] 
This cubic is an isogonal cubic with pivot $G^+$. Besides the vertices of $\Delta$ and the point $G^+$, it passes through the points $H, O^+$ and $I_k, \tilde{K}_k,\tilde{G}_k, k=0,1,2,3$. \\In the above definition of the Thomson cubic, the word "perspectors" may be replaced by "centers" without changing the euclidean curve, see \cite{Gi}. But this is not the case in elliptic and hyperbolic geometry; here we get a different curve of higher degree. A center $[z_1{:}z_2{:}z_3]_\Delta$ of a circumconic belongs to this curve precisely when the coordinates of the corresponding perspector $[p_1{:}p_2{:}p_3]_\Delta$, \vspace*{2mm} \\
 %\[
%\begin{split}
%p_j = &\;
%z_j\big(z_j^{\;2}(\mathfrak{d}_{j+1,j+1}\mathfrak{d}_{j+2,j+2}-\mathfrak{d}_{j+1,j+2}^{\;2}) - z_{j+1}^{\;2}(\mathfrak{d}_{j+2,j+2}\mathfrak{d}_{j,j}-\mathfrak{d}_{j,j+1}^{\;2})\\
%&-z_{j+2}^{\;2}(\mathfrak{d}_{j,j}\mathfrak{d}_{j+1,j+1}-\mathfrak{d}_{j+2,j}^{\;2})+2z_{j+1}z_{j+2}(\mathfrak{d}_{j+1,j+2}\mathfrak{d}_{j,j}-\mathfrak{d}_{j+2,j}\mathfrak{d}_{j,j+1}) \big)\,,
%\end{split}
%\] 
\centerline{$p_j = \;
z_j\big(z_j^{\;2}(\mathfrak{d}_{j+1,j+1}\mathfrak{d}_{j+2,j+2}-\mathfrak{d}_{j+1,j+2}^{\;2}) - z_{j+1}^{\;2}(\mathfrak{d}_{j+2,j+2}\mathfrak{d}_{j,j}-\mathfrak{d}_{j,j+1}^{\;2})$}\\
$\hspace*{22mm} -z_{j+2}^{\;2}(\mathfrak{d}_{j,j}\mathfrak{d}_{j+1,j+1}-\mathfrak{d}_{j+2,j}^{\;2})+2z_{j+1}z_{j+2}(\mathfrak{d}_{j+1,j+2}\mathfrak{d}_{j,j}-\mathfrak{d}_{j+2,j}\mathfrak{d}_{j,j+1}) \big)\,,$\vspace*{2.5mm}\\
satisfy the equation $(\star)$. Points on this curve are: $H$ and $I_k, O_k, k=0,\cdots,3$.

\subsection{Conway's circle, Kiepert perspectors, Hofstadter points, and related objects}
\subsubsection{}\label{subsubsec:lemma}
For all real numbers $x, y, z$ define the points $P_a(x), Q_a(x), P_b(y), Q_b(y)$, $P_c(z)$, $Q_c(z)$ by \vspace*{-1 mm}
\[
\begin{split}
P_a(x) = [0{:}x{:}1]_\Delta,\,  &P_b(y) = [1{:}0{:}y]_\Delta,\, P_c(z) = [z{:}1{:}0]_\Delta,\\
Q_a(x) = [0{:}1{:}x]_\Delta,\,  &Q_b(y) =[y{:}0{:}1]_\Delta,\, Q_c(z) = [1{:}z{:}0]_\Delta.
\end{split}
\]
These six points lie on the conic \vspace*{1 mm}\\
\centerline{ $\{[p_1{:}p_2{:}p_3]_\Delta |\;p_1^2+p_2^2+p_3^2 - (x+\dfrac{1}{x})p_2p_3 - (y+\dfrac{1}{y\,})p_3p_1  - (z+\dfrac{1}{z})p_1p_2 = 0 \}$.}\vspace*{1 mm}
The points $P_a, P_b, P_c$ are collinear precisely when $xyz = -1$.\\

Put
\[
\begin{split}
&X = X(x,y,z) = (Q_c(z)\vee P_a(x)) \wedge (Q_a(x)\vee P_b(y)), \\
&Y = (Q_a(x)\vee P_b(y)) \wedge (Q_b(y)\vee P_c(z)),\,Z = (Q_b(y)\vee P_c(z)) \wedge (Q_c(z)\vee P_a(x)).
\end{split}
\]

\begin{figure} [!thpb]
\includegraphics[height=10cm]{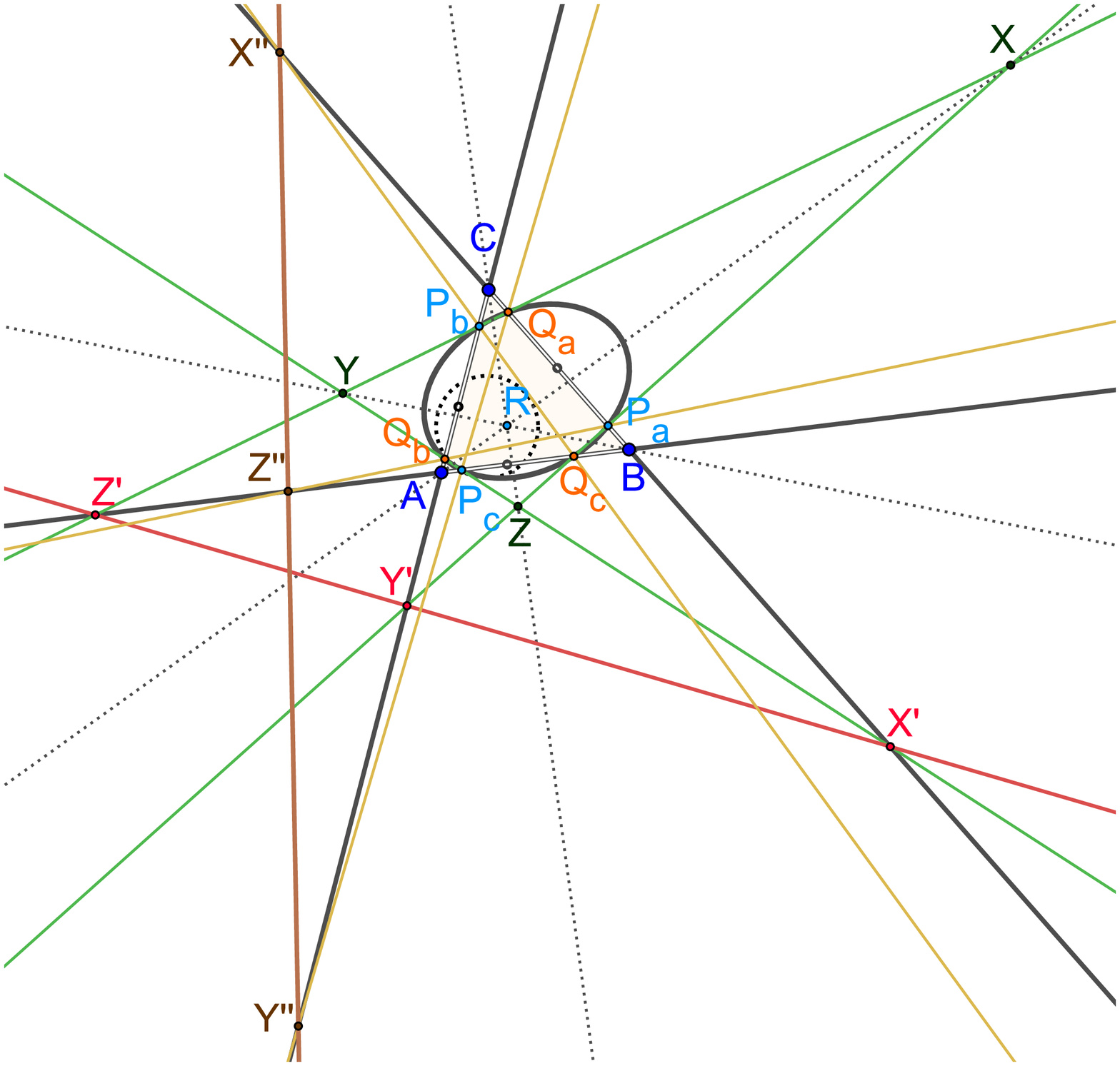}
\caption{Illustration of \ref{subsubsec:lemma}}
\end{figure}
The triple $XYZ$ is perspective to $\Delta$ at $R(x,y,z) = [\dfrac{1}{x+yz},\dfrac{1}{y+zx},\dfrac{1}{z+xy}]_\Delta.$\\
If $xyz\ne 0$, then $R(x,y,z) = R(1/x,1/y,1/z)$.\\

Put \vspace*{-1 mm}
\[
\begin{split}
&X' = (Q_b(y)\vee P_c(z)) \wedge a, \\
&Y' = (Q_c(z)\vee P_a(x)) \wedge b,\,Z' = (Q_a(x)\vee P_b(y)) \wedge c.
\end{split}
\]
The points $X',Y',Z'$ are collinear on the tripolar line of $[x,y,z]_\Delta.$\\

Put \vspace*{-1mm}
\[
\begin{split}
&X'' = (Q_c(z)\vee P_b(y)) \wedge a, \\
&Y'' = (Q_a(x)\vee P_c(z)) \wedge b,\,Z'' = (Q_b(y)\vee P_a(x)) \wedge c.
\end{split}
\]
The points $X'',Y'',Z''$ are collinear on the tripolar line of $[yz,zx,xy]_\Delta.$\vspace*{-2mm}\\

Special cases:\vspace*{1mm}\\ 
$\bullet$ If $x = y = z$, then $R = G$ and $X',Y',Z',X'',Y'',Z''$ lie on the tripolar of $G$.\vspace*{1mm}\\
$\bullet$ Assume 
\[
y = \frac{\sqrt{\mathfrak{d}_{22}}}{\sqrt{\mathfrak{d}_{11}}}x+ \frac{ \mathfrak{d}_{23}\sqrt{{\mathfrak{d}_{11}}} - \mathfrak{d}_{31}\sqrt{\mathfrak{d}_{22}}}{\sqrt{{\mathfrak{d}_{11}}}( \mathfrak{d}_{12}-\sqrt{{\mathfrak{d}_{11}}}\sqrt{{\mathfrak{d}_{22}}}) } \, \text{and}\; z = \frac{\sqrt{\mathfrak{d}_{33}}}{\sqrt{\mathfrak{d}_{11}}}x+ \frac{ \mathfrak{d}_{23}\sqrt{{\mathfrak{d}_{11}}} - \mathfrak{d}_{12}\sqrt{\mathfrak{d}_{33}}}{\sqrt{{\mathfrak{d}_{11}}}( \mathfrak{d}_{31}-\sqrt{{\mathfrak{d}_{11}}}\sqrt{{\mathfrak{d}_{33}}}) }.
\]
\hspace*{2mm}In this case, $d(P_a(x),B) = d(P_b(y),C) = d(P_c(z),A)$ and:\\
\hspace*{2mm}- The points $R(x,y,z), x \in \mathbb{R},$ lie on a circumconic of $\Delta$ through $I$, $Ge$ and $N\!a$. \\
\hspace*{2mm}- For all $x,$ the tripolar of $[yz, zx, xy]_\Delta$ passes through the point $(I\vee Ge)^\delta$.\vspace*{1mm}\\
$\bullet$ Assume 
\[0<k<1 \;\text{and}\; x=\frac{\sinh((1-k)\mathscr{a})}{\sinh(k\mathscr{a})}, y=\frac{\sinh((1-k)\mathscr{b})}{\sinh(k\mathscr{b})},z=\frac{\sinh((1-k)\mathscr{c})}{\sinh(k\mathscr{c})}.
\]
In this case, we get \;$\lim\limits_{k \to 0}R(x,y,z) = \lim\limits_{k \to 1}R(x,y,z) = [\dfrac{\sinh(\mathscr{a})}{\mathscr{a}}:\dfrac{\sinh(\mathscr{b})}{\mathscr{b}}:\dfrac{\sinh(\mathscr{c})}{\mathscr{c}}]_\Delta$.\vspace*{0.5mm}\\

% geprueft 18_06_20.dfw 18_06_20.ggb 18_06_20A.ggb 

\subsubsection{The Conway circle}
Reflect $B$ in $I_1\vee C$ and reflect the mirror image in $B_G$ to get the point $[\mathfrak{d}_{23}\sqrt{{\mathfrak{d}_{11}}}+\mathfrak{d}_{31}\sqrt{{\mathfrak{d}_{22}}}:0:\sqrt{{\mathfrak{d}_{11}}}(\mathfrak{d}_{12}+\sqrt{{\mathfrak{d}_{11}}}\sqrt{{\mathfrak{d}_{22}}})]_\Delta$, denoted by $R_{23}$. Reflect $C$ in $I_1\vee B$ and this point in $C_G$ to get the point $R_{32}=[\sqrt{{\mathfrak{d}_{11}}}(\mathfrak{d}_{31}+ \sqrt{{\mathfrak{d}_{33}}}\sqrt{{\mathfrak{d}_{11}}}):\mathfrak{d}_{23}\sqrt{{\mathfrak{d}_{11}}}+\mathfrak{d}_{12}\sqrt{{\mathfrak{d}_{33}}}:0]_\Delta.$ Both, $R_{23}$ and $R_{32}$, have distance $\mathscr{a}$ from $A$. 
Construct likewise the points $R_{31}, R_{13}$ with distance $\mathscr{b}$ from $B$, and $R_{12}, R_{21}$ with distance $\mathscr{c}$ from $C$.
The six points $R_{12}, R_{21}, R_{23}, R_{32}$, $R_{31}, R_{13}$ lie on a circle with center $I$. The radius $r_{Conway}$ of this circle can be calculated by\vspace*{-2mm}
\[
\begin{split}
 \cosh(r_{Conway}) = \cosh(s)\cosh(r) \;\text{with}\; s &= \frac{1}{2}(\mathscr{a}+\mathscr{b}+\mathscr{c}) =\;  \text{semiperimeter of\;} \Delta_0\\
 \text{and}\, r &= d(I,A_{[I]}) =\; \text{radius of the incircle}\; \mathscr{I}_0.
\end{split}
\]
\textit{Proof}: By reflecting $R_{23}$ in the line $I \vee A$ we get  $R_{32}$, and therefore $d(I,R_{23}) = d(I,R_{32})$. Accordingly, we have $d(I,R_{31}) = d(I,R_{13})$ and $d(I,R_{12}) = d(I,R_{21})$.\\The distance between $A$ and $R_{31}$ is $\mathscr{b}+\mathscr{c}$ and agrees with the distance between $A$ and $R_{21}$. By reflecting $R_{31}$ in the line $I \vee A$ we get  $R_{21}$, thus $d(I,R_{31}) = d(I,R_{21})$. In the same way, we get $d(I,R_{12}) = d(I,R_{32})$ and $d(I,R_{13}) = d(I,R_{23})$.\\
By reflecting $R_{13}$ in the line $I \vee A_{[I]}$ we get  $R_{12}$. Therefore, $A_{[I]}$ is a midpoint of $\{R_{13},R_{12}\}$. The radius $r_{Conway}$ can be calculated by applying the elliptic resp. hyperbolic version of Pythagoras' theorem.\;\;$\Box$
%[d1 w1 + d2 w2,0,w1 (d3 + w1 w2)]
%[w1 (d2 + w1 w3), d1 w1 + d3 w3, 0]

The three points $(R_{23}\vee R_{32}) \wedge a, (R_{31}\vee R_{13}) \wedge b, (R_{12}\vee R_{21}) \wedge c$ are collinear on a line with the equation \vspace*{-1 mm}
\[\sum\limits_{j=1,2,3} (\mathfrak{d}_{j,j+1}\sqrt{\mathfrak{d}_{j+2,j+2}} + \mathfrak{d}_{j+2,j}\sqrt{\mathfrak{d}_{j+1,j+1}})(\mathfrak{d}_{j+1,j+2}-\sqrt{\mathfrak{d}_{j+1,j+1}}\sqrt{\mathfrak{d}_{j+2,j+2}}) x_j = 0\;.\vspace*{0mm}
\] 
%(\text{indices mod}\, 3
The euclidean limit of this line is the tripolar of $X_{86}$.

\begin{figure}[!htbp]
\includegraphics[height=8cm]{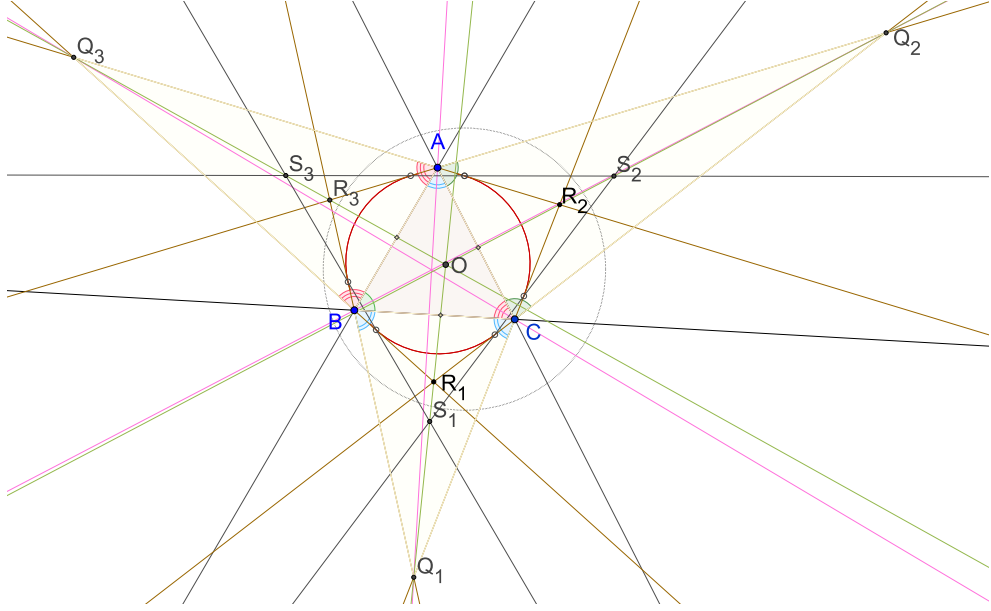}
\caption{
%For example, the pole $A'$ of the line $B \vee C$ with respect to $\mathcal{C}$ can be obtained as follows: Construct the pole of $B \vee C$ with respect to $\tilde{\mathcal{C}}$, then its mirror image in the center $(1:0:0)$ of $\tilde{\mathcal{C}}$ is $A'$.
A description of the red circle is given in \ref{subsubsec:conwaydual}.}\label{fig:dual conway}
\end{figure}
\subsubsection{A dual of Conway's circle theorem}\label{subsubsec:conwaydual} Let $Q_1, Q_2, Q_3$ be the apices of isosceles triangles, erected outwardly on the sides $a, b, c$ of $\Delta_0$ with base angles $\alpha, \beta, \gamma$, respectively. Then there exists a circle with center $O$ which touches the lines $A\vee Q_2, A\vee Q_3$, $B\vee Q_3, B\vee Q_1$, $C\vee Q_1, C\vee Q_2$, see Figure \ref{fig:dual conway}. The euclidean limit of this circle is the circumcircle of $\Delta_0$.\vspace*{-2 mm}\\

Instead of base angles $\alpha, \beta, \gamma$, we can also take base angles $\pi i - \beta -\gamma, \pi i - \gamma -\alpha, \pi i - \alpha -\beta$ to get the same touching circle.\vspace*{-2 mm}\\
 
The triple $Q_1Q_2Q_3$ is perspective to $\Delta$ at a point with euclidean limit $X_{6}$.\vspace*{0.5 mm}
\subsubsection{}\textit{Dualizing} 2.6.1. Let $x, y, z$ be real numbers, and define the points $P_1, P_2, P_3$ by \vspace*{-1 mm}
\[
\begin{split}
&P_1 = [yz \sqrt{\mathfrak{d}_{11}}{:}y\sqrt{\mathfrak{d}_{22}}{:}z\sqrt{\mathfrak{d}_{33}}]_\Delta,\\  &P_2 = [x \sqrt{\mathfrak{d}_{11}}{:}zx\sqrt{\mathfrak{d}_{22}}{:}z\sqrt{\mathfrak{d}_{33}}]_\Delta,\, P_3 = [x \sqrt{\mathfrak{d}_{11}}{:}y\sqrt{\mathfrak{d}_{22}}{:}xy\sqrt{\mathfrak{d}_{33}}]_\Delta.
\end{split}
\]
Then:\\
%P3 ≔ [x √d11, y √d22, x y √d33]
\hspace*{1 mm}$\mu(\angle_{+} BAP_3) = \mu(\angle_{+} CAP_2),\, \mu(\angle_{+} CBP_1) = \mu(\angle_{+} ABP_3),\, \mu(\angle_{+} BCP_1) = \mu(\angle_{+} ACP_2).$\\   
\hspace*{1 mm}The triple $P_1P_2P_3$ is perspective to $\Delta$ at the point $R(x,y,z) = [x\sqrt{\mathfrak{d}_{11}},y\sqrt{\mathfrak{d}_{22}},z\sqrt{\mathfrak{d}_{33}}]_\Delta$.\\

Put  $Q_1 = (P_3\vee B)\wedge(P_2\vee C), Q_2 = (P_1\vee c)\wedge(P_3\vee A), Q_3 = (P_2\vee A)\wedge(P_1\vee B)$. The triple $Q_1Q_2Q_3$ is perspective to $\Delta$
at the isogonal conjugate of $R(x,y,z)$.\\

Special cases:\vspace*{1mm}\\ 
$\bullet\;$ If $x=y=z=1$, then $R = I$.\vspace*{1mm}\\
$\bullet\;$ Assume 
\[
y = \frac{\sqrt{\mathfrak{d}_{11}}\sqrt{\mathfrak{d}_{33}}}{\sqrt{\mathfrak{d}_{22}}\sqrt{\mathfrak{d}_{33}} + x  ( \mathfrak{d}_{31} - \mathfrak{d}_{23})} \, \text{and}\; z = \frac{\sqrt{\mathfrak{d}_{11}}\sqrt{\mathfrak{d}_{22}}}{\sqrt{\mathfrak{d}_{22}}\sqrt{\mathfrak{d}_{33}}+x ( \mathfrak{d}_{12} - \mathfrak{d}_{23})}.\vspace*{-1mm}
\]
In this case, the points $P_1, P_2, P_3$ are the apices of isosceles triangles, erected on the sides
$a, b, c$ of $\Delta_0$ (either all inwardly or all outwardly) with base angles which have all the same measure.\vspace*{-3mm}\\

The points $R(x,y,z), x \in \mathbb{R},$ are called \textit{Kiepert perspectors}; they lie on the conic \vspace*{-1 mm}
\[\big\{[p_1,p_2,p_3]_\Delta \;|\; \sum\limits_{j=1,2,3} (\mathfrak{d}_{j+2,j} - \mathfrak{d}_{j,j+1})p_{j+1}p_{j+2} = 0\big\},\vspace*{-1.5 mm}
\]
which is a circumconic of $\Delta$ through $G$ and $H$. The euclidean limit of this conic is the Kiepert hyperbola. \vspace*{-2.5mm}\\

The isogonal conjugates of the Kiepert perspectors lie on the line through $K$ (= iso\-gonal conjugate of $G$) and $O^+$ (= iso\-gonal conjugate of $H$). This line also passes through the triangle centers $O, \tilde{K},\tilde{K}^\star$. The tripole of $O \vee K$ is a point on the circumcircle $\mathscr{C}_0$ with euclidean limit $X_{110}$.\vspace*{1.5 mm} \\
%Put $S_1 = \tilde{K}^\tau \wedge a$, and let $\mathscr{C}(S_1,A)$ denote the circle with center $S_1$ through vertex $A$. This circle meets the circumcircle $\mathscr{C}_0$ perpendicularly at vertex $A$.  Define cyclically the points $S_2, S_3$ and the circles
%$\mathscr{C}(S_2,B), \mathscr{C}(S_3,C)$. There are two points on $O \vee K$ at which all three circles $\mathscr{C}(S_1,A),\mathscr{C}(S_2,B), \mathscr{C}(S_3,C)$ 
%meet; their euclidean limits are the isodynamic points, cf. \cite{Ev}. The euclidean limit of the line $O \vee K$ is the Brocard axis.\vspace*{1.5 mm}\\ 
$\bullet\;$  Let $k$ be a real number, different from $0$ and $1$, and assume \vspace*{-1mm}
\[x=\frac{\sinh(k\mathscr{\alpha})}{\sinh((1-k)\mathscr{\alpha})},\; y=\frac{\sinh(k\mathscr{\beta})}{\sinh((1-k)\mathscr{\beta})},\;z=\frac{\sinh(k\mathscr{\gamma})}{\sinh((1-k)\mathscr{\gamma})}.\vspace*{-1mm}
\]
In this case, we call the point $R(x,y,z)$ \textit{Hofstadter k-point}, according to the euclidean case. The isogonal conjugate of the Hofstadter k-point is the Hofstadter (1-k)-point. Here are some examples of Hofstadter points: The $\frac{1}{2}$-point is $I$, the (-1)-point $H$ and the 2-point is $O^+$. The limes of the k-point, as k approaches zero, is $[\alpha{:}\beta{:}\gamma]_\Delta$.

\subsection{The Lemoine conic, the Lemoine axis and Tucker circles} \vspace*{-1mm}
%\subsubsection{A conic with Lemoine circle as euclidean limit} 
%Let $P$ be a point not on a sideline of $\Delta$. Put ${P_{11} := \text{par}(a,P) \wedge a}, {P_{12}:= \text{par}(a,P)\wedge b}$ and define the points $P_{13}, P_{22}, P_{23}, P_{21}$, $P_{33}, P_{31}, P_{32}$ , accordingly.\\
%The points $P_{11}, P_{22}, P_{33}$ are collinear on the line $P^\delta$, and the other six points lie on a conic. For $P = \tilde{K}$, this conic has the symmetry line $O\vee K$ (see []). But in general, it is not a circle, as it is in the euclidean case. We call the conic \textit{Lemoine conic}.
\subsubsection{Lemma} 
Let $P = [p_1:p_2:p_3]_\Delta$ be a point not on a sideline of $\Delta$, and let $\ell = \{{[x_1:x_2:x_3]_\Delta} |\, q_1 x_1 +q_2 x_2 + q_3 x_3 = 0 \}$ be a line that does not contain any vertex of $\Delta$. Put $P_{11} := \ell \wedge a, P_{12} := (P_{11} \vee P) \wedge b , P_{13} := (P_{11} \vee P)\wedge c$, and define the points $P_{22}, P_{23}, P_{21}$, $P_{33}$, $P_{31}, P_{32}$ accordingly.\\
The points $P_{12}, P_{13}, P_{21}$, $P_{23}$, $P_{31}, P_{32}$ lie on the conic \vspace*{-1mm}
\[ 
\begin{split}
\big\{[x_1{:}x_2{:}x_3]_\Delta |& \sum\limits_{j=1,2,3} \big(\frac{q_j(p_{j+1}q_{j+1} + p_{j+2}q_{j+2})}{p_j}x_{j}^2\\
 &- \frac{p_j q_j(p_j q_j+p_{j+1} q_{j+1}+p_{j+2} q_{j+2})+2p_{j+1} q_{j+1}p_{j+2} q_{j+2}}{p_{j+1} p_{j+2}}x_{j+1}x_{j+2}\big)  = 0\big\}.\vspace*{-1mm}
\end{split}
\]
We call this conic the \textit{conic associated with} the pair $(P,\ell)$.
Its perspector is \vspace*{-1mm}
\[ [\frac{p_1}{(2 p_1 q_1 (p_2 q_2 + p_3 q_3) + p_2 p_3 q_2 q_3)}:\cdots:\cdots]_\Delta.\vspace*{1mm}
\]
\textit{Examples:}\\
$\bullet\;\;$ The perspector of the conic associated with $(P,P^\tau)$ is the point $P$.\\
$\bullet\;\;$ The conic associated with $(\tilde{K},\tilde{K}^\delta)$ is the Lemoine conic which is described in \cite{Ev}; its euclidean limit is the First Lemoine circle. The line $O\vee K$ $(= {O\vee \tilde{K}})$ is a symmetry axis of this conic.\\
$\bullet\;\;$ The conic associated with $G$ and $H^\delta$ is a conic with perspector \\
\centerline{$[1/(\mathfrak{c}_{23} + 2\mathfrak{c}_{31}+2\mathfrak{c}_{12}):1/(2\mathfrak{c}_{23} + \mathfrak{c}_{31}+2\mathfrak{c}_{12}):1/(2\mathfrak{c}_{23} + 2\mathfrak{c}_{31}+\mathfrak{c}_{12})
]_\Delta$.}\\
$\bullet\;\;$ The conic associated with incenter $I$ and the orthotransversal of $I$ has $I$ as a symmetry point and a perspector with euclidean limit $X_{10390}$.\\
$\bullet\;\;$ In euclidean geometry, the conic associated with the Lemoine point $K$ and the line at infinity is the First Lemoine circle; the Second Lemoine circle is associated with $K$ and the tripolar of $X_{25}$.

\subsubsection{The Lemoine axis and the apollonian circles} 
The tripolar of the Lemoine point $\tilde{K}$ is called \textit{Lemoine axis}. This axis is perpendicular to the line $O\vee K$ $(= {O^+\vee \tilde{K}})$. Let $L_1, L_2, L_3$ be the intersection points of the Lemoine axis with the sidelines $a, b, c$, respectively. The circle $\mathscr{C}(L_1,A)$ with center $L_1$ through vertex $A$ meets the circumcircle $\mathscr{C}_0$ perpendicularly. 
Mutatis mutandis, this is also true for the circles $\mathscr{C}(L_2,B), \mathscr{C}(L_3,C)$. We will call these circles \textit{apollonian circles}, as they are called in the euclidean case. There are two points, $J_+$ and $J_+$,\vspace*{0.5mm}\\
\centerline{$J_\pm = [(\mathfrak{c}_{23}-\mathfrak{c})(\,\mathfrak{c}+\mathfrak{c}_{23}-\mathfrak{c}_{31}-\mathfrak{c}_{12}\;\pm\,\sqrt{\frac13|\det(\mathfrak{C})|\,}\,):\cdots:\cdots]_\Delta$,}\vspace*{0.5mm}\\
on the line $O\vee K$ at which all three apollonian circles meet; their euclidean limits are the isodynamic points. The points $\tilde{K}^{\star}$ and $\tilde{K}^{\tau}\wedge(K\vee O)$ are the midpoints of $J_-$ and $J_+$, and $J_-, O, J_+, \tilde{K}$ form a harmonic range.\\
\textit{Remark}: An elliptic/hyperbolic version of Apollonius' Theorem (see \cite{TL} ch. VIII for the elliptic/spherical version together with its proof) states that\vspace*{1mm}\\ 
\centerline{$\displaystyle\mathscr{C}(L_1,A) = \big\{P\;|\;\frac{\sinh(\frac{1}{2} \mu([P,B]_+))}{\sinh(\frac{1}{2} \mu([P,C]_+)}= \frac{\sinh(\frac{1}{2} \mathscr{c})}{\sinh(\frac{1}{2} \mathscr{b})}\big\}$.}
\subsubsection{Tucker hexagons and Tucker circles} 
Let $P_1Q_2P_3Q_1P_2$ $Q_3P_1$ be a closed poly\-gonal chain with vertices on the sidelines of $\Delta_0$, $P_1, Q_1$ on $a$, $P_2, Q_2$ on $b$, $P_3, Q_3$ on $c$, and assume that none of these points is a vertex of $\Delta$.
In euclidean geometry, the polygon $P_1Q_2P_3Q_1P_2Q_3$ is called a \textit{Tucker hexagon} of $\Delta_0$, if its line segments are alternately parallel\,/\,antiparallel to the sidelines of $\Delta$, $P_1Q_2$ parallel (or antiparallel) to $c$, $Q_2P_3$ antiparallel (parallel) to $a$, etc. The vertices of a Tucker hexagon are always concyclic and the corresponding circle is called a \textit{Tucker circle}.\vspace*{-2mm}\\
%Two examples: The First and the Second Lemoine circle are Tucker circles.\\
%In elliptic and in hyperbolic geometry, (anti-) parallelism of lines is not defined. As a substitute, we take the following definition introduced by Akopyan []:\\

We will show that Tucker hexagons and Tucker circles also exist in the elliptic and in the hyperbolic plane.\\
The concept of \textit{antiparallelism of lines with respect to an angle} can be transferred from euclidean geometry to elliptic and to hyperbolic geometry.  Let us explain this for the angle $\angle_+ BAC$ of the triangle $\Delta_0$. (See Akopyan \cite{Ak} for a more detailed description).\\
Given two lines $g, h$ such that $\#\{b,c,g,h\} = 4$, then $g$ is \textit{antiparallel} to $h$ with respect to $\angle_+ BAC$  $\,$ iff one of the following two conditions holds:\\
condition 1: $g\wedge h = A$ and $h$ is the mirror-image of $g$ in the angle-bisector $A \vee I_0$ of $\angle_+ BAC$.\\
condition 2: Define $P_1{:=\,}b\wedge g, P_2{:=\,}g\wedge c,P_3{:=\,}c\wedge h, P_4{:=\,}h\wedge b$. $\#\{P_1,P_2,P_3,P_4\} \ge 3$, and the segments $[P_1, P_2]_+,[P_2, P_3]_+, [P_3, P_4]_+, [P_4, P_1]_+$ are cords of a circle. (We recall that a cord of a circle is a closed segment whose boundary points lie on the circle and whose inner points lie inside the circle.)\\
It can be easily verified that, if $A$ is neither a point on $g$ nor on $h$, then $g$ is antiparallel to $h$ with respect to $\angle_+ BAC$ precisely when $\mu(\angle_+(P_1P_2A))-\mu(\angle_+(P_2P_1A)) = \mu(\angle_+(P_3P_4A))-\mu(\angle_+(P_4P_3A))$. \\
Now we define parallelism between lines with respect to the angle $\angle_+ BAC$:
Two lines $g, h$ are parallel with respect to  $\angle_+ BAC$ precisely when the mirror image of $h$ in the angle bisector $A \vee I_0$ is antiparallel to $g$. 
\begin{figure}[!htbp]
\includegraphics[height=9cm]{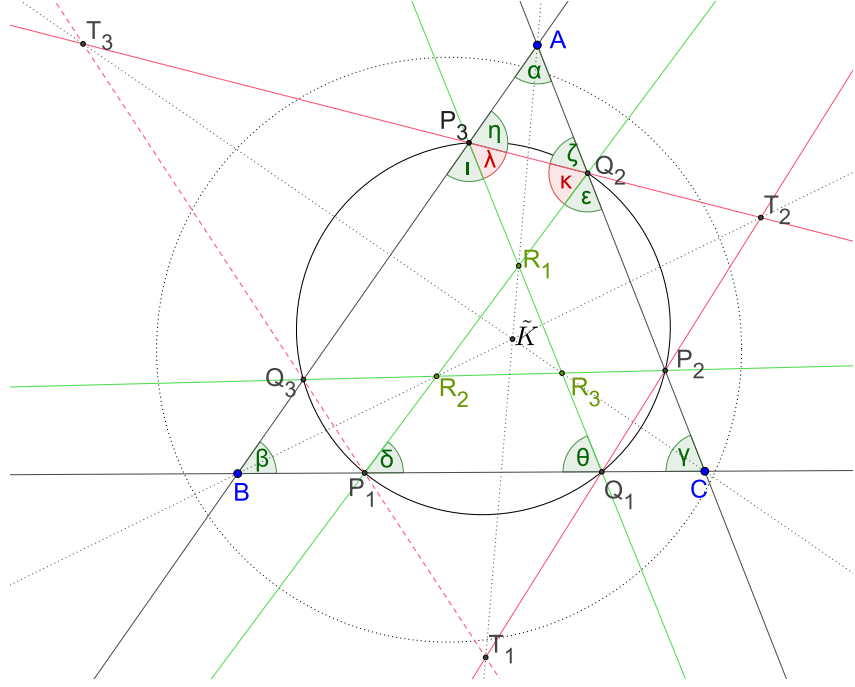}
\caption{A Tucker hexagon and its circle.}\label{fig:T_hexagon}
\end{figure}

Choose a point $P_1$ on $a$, different from $B$ and $C$, and construct successively points $Q_2$ on $b$, $P_3$ on $c$, $Q_1$ on $a$, $P_2$ on $b$ and $Q_3$ on $c$ such that $P_1\vee Q_2$ is parallel to $c$, $Q_2\vee P_3$ antiparallel to $a$, $P_3\vee Q_1$ parallel to $b$, $Q_1\vee P_2$ antiparallel to $c$, $P_2\vee Q_3$ parallel to $a$. Then, $Q_3\vee P_1$ is antiparallel to $b$ and the points $P_1,\cdots, Q_3$ are concyclic.\\
\textit{Proof}: Since $P_1\vee Q_2$ is parallel to $c$, $Q_2\vee P_3$ antiparallel to $a$, and $P_3\vee Q_1$ parallel to $b$, we have the equations $\alpha - \beta = \varepsilon - \delta,\, \beta - \gamma = \zeta - \eta,\, \gamma - \alpha = \theta - \iota$, see Figure \ref{fig:T_hexagon}. We can prove that the four points $P_1, Q_2, P_3, Q_1$ are concyclic on a circle $\mathcal{C}$ by showing$\;\;\;\kappa - \lambda = \theta-\delta\;\;\;(\star)$.\vspace*{1mm}
\centerline{$\text{Proof of} \,\;(\star):\;\; \kappa -\lambda = (\pi i - \varepsilon -\zeta) - (\pi i - \eta - \iota)\,$}\\
$ \hspace*{59mm}      = (\eta -\zeta) + \iota - \varepsilon$\\
$	\hspace*{59mm}			= (\gamma - \beta) + (\theta + \alpha - \gamma) - (\delta + \alpha -\beta)$\\
$ \hspace*{59mm}      = \theta - \delta$ \vspace*{1.5mm}

$Q_1\vee P_2$ is antiparallel to $P_1\vee Q_2$ and to c precisely when $P_2$ is a point on $\mathcal{C}$, concyclic with $P_1, Q_2, Q_1$. $P_2\vee Q_3$ is antiparallel to $Q_2\vee P_3$ and therefore parallel to $a$ precisely when $Q_3$ is a point on $\mathcal{C}$, concyclic with $P_1, Q_1, P_2$. The point $Q_3$ is a point on $\mathcal{C}$, together with $P_1, Q_1, P_3$. Therefore, $Q_3\vee P_1$ is antiparallel to $P_3\vee Q_1$ and to $b$. The polygon $P_1Q_2P_3Q_1P_2Q_3$ is a Tucker hexagon of $\Delta_0$; its vertices are concyclic.$\;\;\Box$\vspace*{-1mm}\\

Starting from $Q_3 = [q_{31}{:}1{:}0]_\Delta$, we calculate the coordinates of the other vertices of the Tucker hexagon. We get, for example,\vspace*{1 mm}\\
\noindent$\displaystyle\hspace*{15 mm} P_1 = [0{:}1{:}p_{13}]_\Delta,\;\;\; p_{13} = \frac{q_{31}(s - \mathfrak{c}_{12})}{s -\mathfrak{c}_{23}} ,$\vspace*{1 mm}\\
\noindent$\hspace*{46.5 mm} s = \sqrt{\mathfrak{c}}\,\sqrt{(q_{31},1,0){\,\scriptscriptstyle{[\mathfrak{C}]}\,}(q_{31},1,0)}\; - \mathfrak{c} q_{31}$\vspace*{1 mm}\\
\noindent$\hspace*{49 mm}  = \sqrt{\mathfrak{c}}\,\sqrt{\mathfrak{c}(q_{31}^{\;2}+1)+2\mathfrak{c}_{12}q_{31}}\;-\mathfrak{c}\,q_{31},$\\
\noindent$\displaystyle\hspace*{15 mm} P_2 = [p_{21}{:}0{:}1]_\Delta,\;\;\; p_{21} = \frac{(s-{\mathfrak{c}})(s\,(\mathfrak{c}_{31}-\mathfrak{c})+\mathfrak{c}\,(2\mathfrak{c}_{12}-\mathfrak{c}_{31})-1)}{2\mathfrak{c}\,(\mathfrak{c}-\mathfrak{c}_{12})(s-\mathfrak{c}_{12})}.$\vspace*{1 mm}\\
The center $T$ of the Tucker circle can be calculated by $T = ((P_1-Q_3) \vee (P_2-Q_3))^\delta$. By using CAS, it can be verified that $T$ lies on the line $ O \vee K \;(=O \vee \tilde{K})$. The coordinates of $T$ can then be expressed by\\
$T = [(\mathfrak{c}{-}\mathfrak{c}_{23})(\mathfrak{c}{+}\mathfrak{c}_{23}{-}\mathfrak{c}_{31}{-}\mathfrak{c}_{12}{+\,}t){:}(\mathfrak{c}{-}\mathfrak{c}_{31})(\mathfrak{c}{+}\mathfrak{c}_{31}{-}\mathfrak{c}_{12}{-}\mathfrak{c}_{23}{+\,}t){:}(\mathfrak{c}{-}\mathfrak{c}_{12})(\mathfrak{c}{+}\mathfrak{c}_{12}{-}\mathfrak{c}_{23}{-}\mathfrak{c}_{31}{+\,}t)]_\Delta$\\
$ \; \textrm{with}$\\
\noindent$\hspace*{20 mm}\displaystyle t = \frac{(2\mathfrak{c}_{12}\mathfrak{c}_{23}\mathfrak{c}_{31}-\mathfrak{c}\,(\mathfrak{c}_{12}^{\;2}{+}\mathfrak{c}_{23}^{\;2}{+}\mathfrak{c}_{31}^{\;2}{-1}))((p_{13}{+}1)\,r_{31}{\;-\;}(q_{31}{+}1)r_{13})}{r_{31}\,s_{13} + r_{13}\,s_{31}},$\vspace*{1mm}\\
$\hspace*{17mm} r_{31} = \sqrt{\mathfrak{c}(q_{31}^{\;2}{+}1)+2\mathfrak{c}_{12}q_{31}} ,\;\;\;\; r_{13} = \sqrt{\mathfrak{c}(p_{13}^{\;2}+1)+2\mathfrak{c}_{23}p_{13}},$\vspace*{0.5mm}\\

$\hspace*{12.6mm} s_{31} = q_{31}(\mathfrak{c}(\mathfrak{c}_{12}{-}\mathfrak{c}_{23}{+}\mathfrak{c}_{31})-2\mathfrak{c}_{12}\mathfrak{c}_{31}{+}1)+\mathfrak{c}(\mathfrak{c}_{12}{+}\mathfrak{c}_{23}{-}\mathfrak{c}_{31})-2\mathfrak{c}_{12}\mathfrak{c}_{23}{+}1\;,$\vspace*{0.5mm}\\
$\hspace*{16.7mm} s_{13} = p_{13}(\mathfrak{c}(\mathfrak{c}_{12}{-}\mathfrak{c}_{23}{-}\mathfrak{c}_{31})+2\mathfrak{c}_{23}\mathfrak{c}_{31}{-}1)-\mathfrak{c}(\mathfrak{c}_{12}{+}\mathfrak{c}_{23}{-}\mathfrak{c}_{31})+2\mathfrak{c}_{12}\mathfrak{c}_{23}{-}1\;.$\vspace*{2mm}\\
Constructions with GeoGebra show that the perspector of the Tucker circle is, in general, not a point on the line $O\vee \tilde{K}$.\vspace*{-1mm}\\

\begin{figure}[!b]
\includegraphics[height=8.5cm]{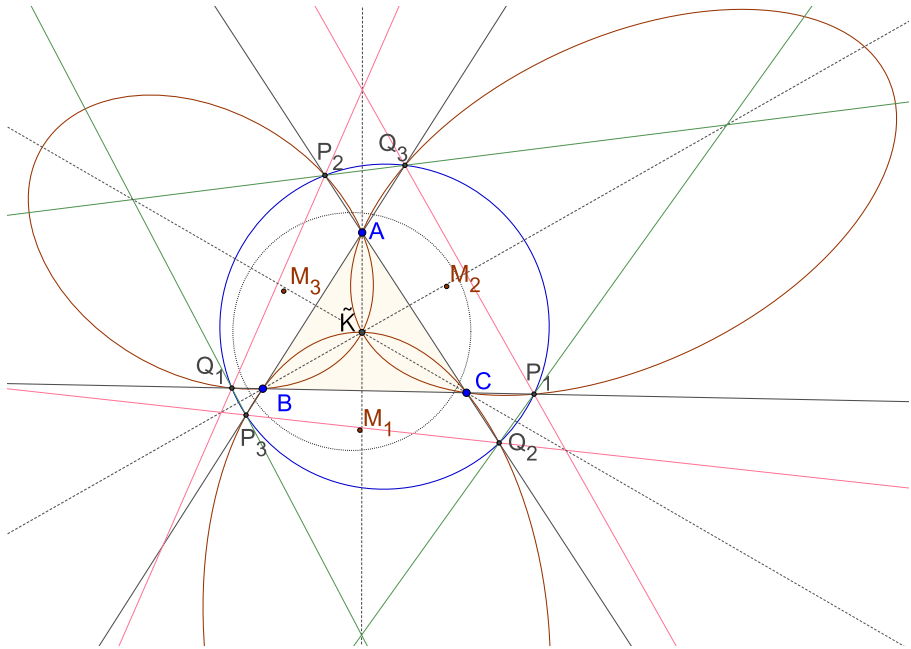}\label{fig:LC}
\caption{The $3^{rd}$ Lemoine circle (blue)\newline}\label{fig:Lemoine_circle}
\end{figure}
\begin{figure}
\includegraphics[height=8.5cm]{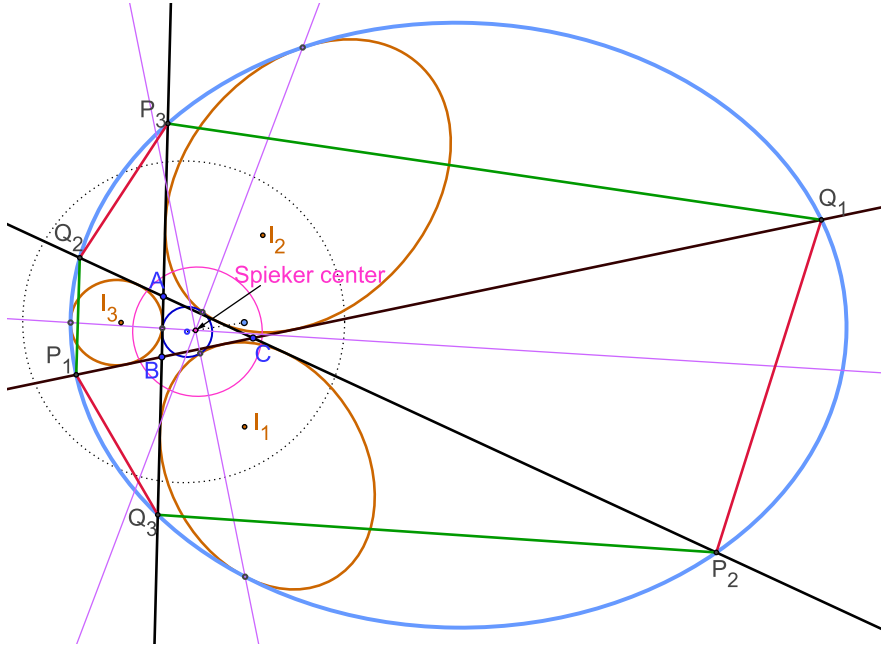}
\caption{The Apollonius circle (light blue). The smaller blue circle is the Hart-Feuerbach circle which touches the (brown) excircles externally. The Spieker center is the radical center of the excircles.}\label{fig:Apollonius_circle}
\end{figure}

\begin{figure}[!t]
\includegraphics[height=9cm]{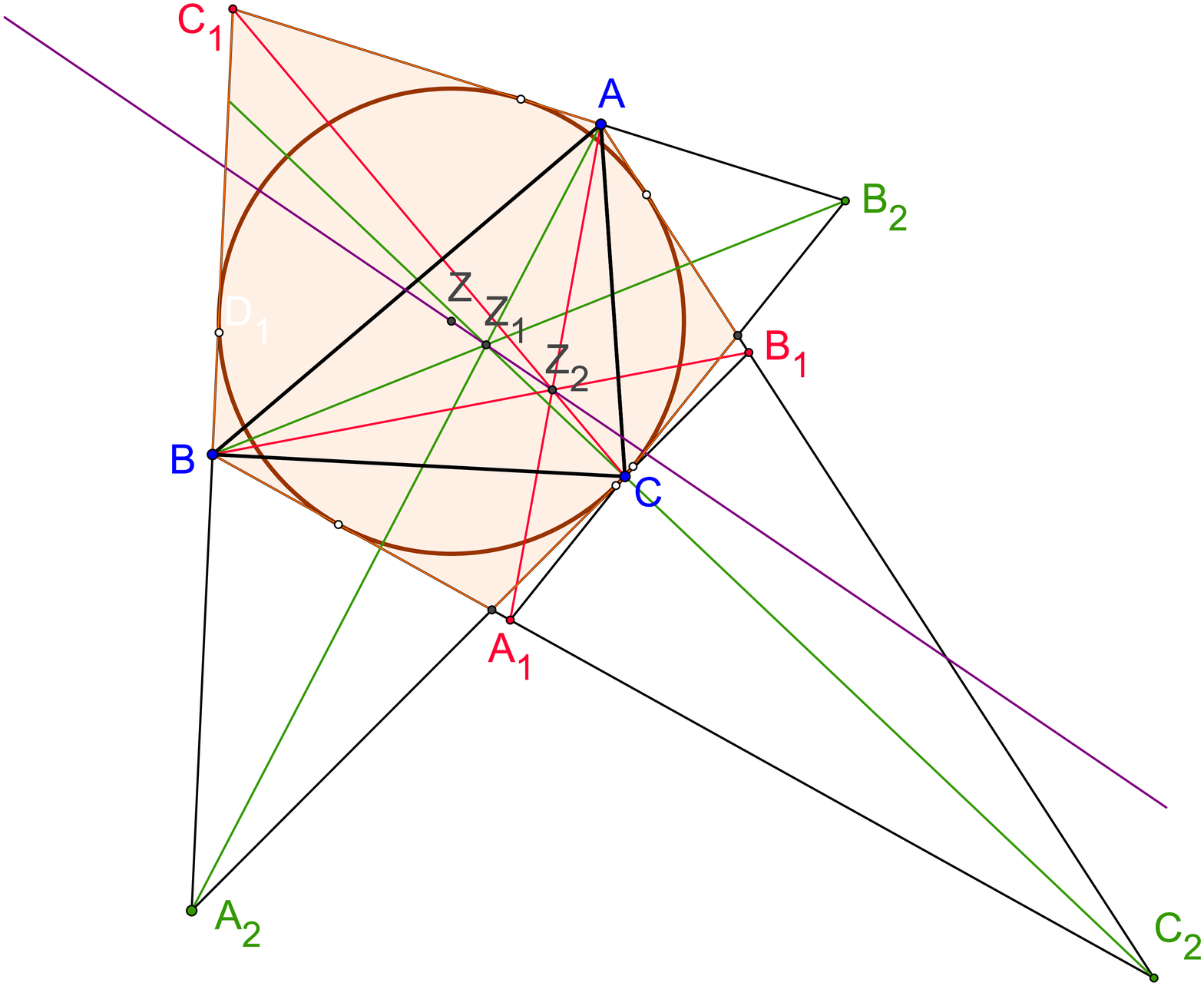}
\caption{A triangle in the euclidean plane with duals of a Tucker hexagon and a Tucker circle. }\label{fig:DualTucker}
\end{figure}

Define the points $R_1,R_2,R_3,T_1,T_2,T_3$ by \vspace*{1mm}\\
\centerline{$
R_1 = (P_1\vee Q_2) \wedge (Q_1\vee P_3) , R_2 = (P_2\vee Q_3) \wedge (Q_2\vee P_1) ,  R_3 = (P_3\vee Q_1) \wedge (Q_3\vee P_2) ,$}\vspace*{0.5mm}\\
\centerline{$
\;\;T_1 = (P_1\vee Q_3) \wedge (Q_1\vee P_2) , \,T_2 = (P_2\vee Q_1) \wedge (Q_2\vee P_3) , \; T_3 = (P_3\vee Q_2) \wedge (Q_3\vee P_1).$
}\vspace*{2mm}\\
The triples $R_1R_2R_3$ and $T_1T_2T_3$ are both perspective to $\Delta$, the perspector being the Lemoine point $\tilde{K}$.
This can be shown by calculations with the help of CAS.\vspace*{-1mm}\\

Constructions indicate: The points $a\wedge(P_2 \vee Q_3), b\wedge(P_3 \vee Q_1), c\wedge(P_1 \vee Q_2)$ are collinear on a line which we name $g$, and the points
$a\wedge(Q_2 \vee P_3), b\wedge(Q_3 \vee P_1), c\wedge(Q_1 \vee P_2)$ are collinear on a line $r$. The lines $g$ and $r$ meet at a point on the dual of $T$.\\

%p23 = (√\mathfrak{c}√(\mathfrak{c}(p12 + 1)^2 - 2 p12 w1^2) (w1^2 - w2^2) - \mathfrak{c}(p12 (w1^2 - w2^2) - w1^2 - w2^2) - w1^2 w2^2)/(p12 w1^2 (2 \mathfrak{c}- w1^2))
%p12 (√\mathfrak{c}(w1^2 - w2^2) √(\mathfrak{c}(p12^2 + 2 p12 + 1) - 2 p12 w1^2) + \mathfrak{c}(p12 (w1^2 - w2^2) - w1^2 - w2^2) + w1^2 w2^2)/(2 \mathfrak{c}(p12 (w1^2 - w2^2) - w2^2) + w2^4)
%2(w1^2 - w2^2) (√( (p12^2 + 2 p12 + 1) - 2 \mathfrak{c}p12 w1^2) +\mathfrak{c}p12) - ( -\mathfrak{c}w1^2 - \mathfrak{c}w2^2 + w1^2 w2^2)
Let $l_1, l_2, l_3$ be three lines passing through $\tilde{K}$ and  parallel to $a,b,c$ with respect to $\angle_+ BAC$, $\angle_+ CBA$, $\angle_+ ACB$, respectively. 
The six points $l_1\wedge b, l_1\wedge c, l_2\wedge c, l_2\wedge a, l_3\wedge a, l_3\wedge b$ lie on a circle, the \textit{$1^{st}$ Lemoine circle}. The point $Q_3 =
l_1 \wedge c$ has coordinates\vspace*{1mm}\\
$\displaystyle [q_{31}{:}1{:}0]_\Delta = [\frac{(\mathfrak{c}_{12}-\mathfrak{c}_{23})(\mathfrak{c}_{\,}(\mathfrak{c}_{12}+\mathfrak{c}_{23}-\mathfrak{c}_{31})+\sigma\sqrt{c}-1)+2(1-\mathfrak{c}_{12}\mathfrak{c}_{23})(\mathfrak{c}-\mathfrak{c}_{23})}{4(1-\mathfrak{c}_{12}\mathfrak{c}_{23})(\mathfrak{c}-\mathfrak{c}_{31})}{:}1{:}0]_\Delta,$\vspace*{0.5mm}\\
\noindent$\hspace*{2mm}\sigma = \sqrt{\mathfrak{c}\,\big(\mathfrak{c}_{12}^{\;2}{+}\mathfrak{c}_{23}^{\;2}{+}\mathfrak{c}_{31}^{\;2}-2(\mathfrak{c}_{12}\mathfrak{c}_{23}{+}\mathfrak{c}_{23}\mathfrak{c}_{31}{+}\mathfrak{c}_{31}\mathfrak{c}_{12})+5\big)+4\mathfrak{c}_{12}\mathfrak{c}_{23}\mathfrak{c}_{31}-2(\mathfrak{c}_{12}{+}\mathfrak{c}_{23}{+}\mathfrak{c}_{31})}$.\vspace*{0.5mm}\\
Constructions indicate: The first Lemoine circle is the conic associated with $(\tilde{K}, g)$, where $g$ is the line described above.
\vspace*{-1.5mm}\\

Let $l_1, l_2, l_3$ be three lines passing through $\tilde{K}$ and antiparallel to $a,b,c$, respectively.
The six points $l_1\wedge b, l_1\wedge c, l_2\wedge c, l_2\wedge a, l_3\wedge a, l_3\wedge b$ lie on a circle, the $2^{nd}$\textit{Lemoine circle}.The point $Q_3 =
l_2 \wedge c$ has coordinates\vspace*{0.5mm}\\
$\displaystyle [q_{31}{:}1{:}0]_\Delta = [\frac{(\mathfrak{c}_{12}-\mathfrak{c}_{23})(\mathfrak{c}_{\,}(\mathfrak{c}_{12}+\mathfrak{c}_{23}-\mathfrak{c}_{31})+\sigma\sqrt{c}-1)+2(1-\mathfrak{c}_{12}\mathfrak{c}_{23})(\mathfrak{c}-\mathfrak{c}_{31})}{\mathfrak{c}_{31}^{\;2}(\mathfrak{c}+\mathfrak{c}_{12})-2(\mathfrak{c}\,\mathfrak{c}_{23}\mathfrak{c}_{31})-\mathfrak{c}_{23}+\mathfrak{c}_{31})+\mathfrak{c}-\mathfrak{c}_{12}}{:}1{:}0]_\Delta,$\\
\noindent$\hspace*{2mm}\sigma$ as above.\vspace*{0.5mm}\\
Constructions indicate: The second Lemoine circle is the conic associated with $(\tilde{K}, r)$.
\vspace*{-1.5mm}\\

Let $Q_2$ and $P_3$ be the second intersections of the circumcircle of the triangle $(B\tilde{K}C)_0$ with the sidelines $b$ and $c$, respectively. Define the 
points $Q_3,P_1$ and $Q_1,P_2$ likewise. The hexagon $P_1Q_2P_3Q_1P_2Q_3$ is a Tucker hexagon. The associated Tucker circle is the \textit{$3^{rd}$ Lemoine circle}. See Figure \ref{fig:Lemoine_circle}. We do not present the coordinates of $Q_3$, they are very complicated. \vspace*{-2mm}\\

In euclidean geometry, a circle which internally touches the three excircles of a triangle is called the \textit{Apollonius circle} of this triangle. Grinberg and Yiu \cite{GY} showed that this Apollonius circle is a Tucker circle. As constructions indicate, this seems to be true also in elliptic and in hyperbolic geometry, see Figure \ref{fig:Apollonius_circle}.\vspace*{-5.5mm}\\

\subsubsection{Dualizing Tucker hexagons and Tucker circles} 
We can easily formulate a definition of parallelism and antiparallelism of points $P, Q$ with respect to a segment which is dual to the definition of parallelism and antiparallelism of lines with respect to an angle. (We use the names "parallelism/antiparallelism", even though they do not fit well.) \\
A dual version of a Tucker hexagon and a Tucker circle is shown in Figure \ref{fig:DualTucker}.\\

\subsection{Pseudoparallels of the sidelines and their duals} 
\subsubsection{Lemma / part 1}
Let $Q = [q_1{:}q_2{:}q_3]_\Delta$ be a point not on a sideline of $\Delta$, and let $Q_1, Q_2,Q_3$ be the intersections of the tripolar $Q^\tau$ with the sidelines $a, b, c$, respectively. Let $R_1, R_2, R_3$ be any points in $\mathcal{P}$, with coordinates $R_k = [r_{k,1}{:}r_{k,2}, {:}r_{k,3}]_\Delta, k=1,2,3$. Define three lines $\ell_{\!1}, \ell_{\!2},\ell_{\!3}$ by $\ell_{\!k} = R_k\vee Q_k$.\\
Then: The points $\ell_{\!2}\wedge \ell_{\!3}, \ell_{\!3}\wedge \ell_{\!1},\ell_{\!1}\wedge \ell_{\!2}$ form a triple perspective to $\Delta$.
The perspector is the point\vspace*{-1mm}\\
$\hspace*{20mm}P =\displaystyle[\frac{r_{11}}{\tilde{r}_{1}}:\frac{r_{22}}{\tilde{r}_{2}}:\frac{r_{33}}{\tilde{r}_{3}}]_\Delta$,  $\;\;\displaystyle\tilde{r}_{k} = \frac{r_{k1}}{q_1}+\frac{r_{k2}}{q_2}+\frac{r_{k3}}{q_3}, \;k=1,2,3$.\vspace*{1.5mm}\\
We look at special cases:\\
$\bullet$ If $R_1 R_2 R_3$ is the anticevian triple of $Q$, then $P=Q$.\\
$\bullet$ Let $R = [r_1{:}r_2{:}r_3]_\Delta$ be a point not on a sideline of $\Delta$ and $R_1 R_2 R_3$ its anticevian triple.\\
\noindent \hspace*{2mm}If $Q$ is a point on the circumconic of $\Delta$ with perspector $R$, then $P=Q$.\\
$\bullet$ If $\,R_1 R_2 R_3\,$ is the anticevian triple of $\,G^+$ and $\,Q = H\,$,$\;$then $P$ is a triangle center with\\
\noindent \hspace*{2mm}euclidean limit $X_{2996}$.
\subsubsection{A special case: pseudoparallels of the sidelines.} If in the previous lemma $Q = G$, we call the lines $\ell_{\!1}, \ell_{\!2},\ell_{\!3}$ pseudoparallels of the sidelines $a, b, c$. We look at different cases for the triple $R_1 R_2 R_3$.\vspace*{1mm} \\
$\bullet\,$ We put $P_1\! := \sinh(\mathscr{a})\,A' +\, \cosh(\mathscr{a})\,B, P_2\! := \sinh(\mathscr{b})\, B' +\, \cosh(\mathscr{b})\,C, P_3\! := \sinh(\mathscr{c})\, C' + \cosh(\mathscr{c})\,A$. In this case, $ d(\ell_{\!1},B) = d(\ell_{\!1},C) = \mathscr{a},\, d(\ell_{\!2},C) = d(\ell_{\!2},A) = \mathscr{b},\, d(\ell_{\!3},A) = d(\ell_{\!3},B) = \mathscr{c}$ and\vspace*{-0.5mm} \\
\centerline{\hspace*{10mm}$P = \displaystyle[\frac {\sinh^2(\mathscr{a})}{\cosh(\mathscr{a})\det(\mathfrak{D})+\mathfrak{d}_{11}+\mathfrak{d}_{12}+\mathfrak{d}_{31}}):\cdots:\cdots]_\Delta.$}\vspace*{1mm}\\
The euclidean limit of this point is $X_6$.\\
$\bullet\,$ We use the abbreviations $t_1{\,:=\,}1-\mathfrak{c} \mathfrak{c}_{23}, t_2{\,:=\,}1-\mathfrak{c} \mathfrak{c}_{31}, t_3{\,:=\,}1-\mathfrak{c} \mathfrak{c}_{12},$  $t_{123\!}:= t_1{-}t_2{-}t_3,$\\ \noindent$\,\hspace*{3mm} t_{231\!}:= t_2{-}t_3{-}t_1, t_{312\!}:= t_3{-}t_1{-}t_2$ and  define the points $P_1,P_2,P_3$ by\\
\centerline{$P_1 = [t_1t_{123}:t_2t_{231}-t_2 t_3-k\sqrt{t_2}\sqrt{t_3}\,t_{123}:t_3t_{312}-t_2 t_3-k\sqrt{t_2}\sqrt{t_3}\,t_{123}]_\Delta$,}
\centerline{$P_2 = [t_1t_{123}-t_3 t_1-k\sqrt{t_3}\sqrt{t_1}\,t_{231}:t_2t_{231}:t_3t_{312}-t_3 t_1-k\sqrt{t_3}\sqrt{t_1}\,t_{231}]_\Delta$,}
\centerline{$P_3 = [t_1t_{123}-t_1 t_2-k\sqrt{t_1}\sqrt{t_2}\,t_{312}:t_2t_{231}-t_1 t_2-k\sqrt{t_1}\sqrt{t_2}\,t_{312}:t_3t_{312}]_\Delta$.}\vspace*{1mm} \\ 
$P_1, P_2, P_3$ are the vertices of the first circumcircle-midarc-triangle of $\Delta_0$ for $k{\,=\,}1$ and the vertices  of the second circumcircle-midarc-triangle for $k{\,=\,}{-1}$. The perspector $P$ has euclidean limit $X_{56}$ {resp.\,}$X_{55}$ if $k{\,=\,}1$ {resp.\;}$k{\,=\,}{-1}$, but $P$ need not lie on the line $I\vee O$. \vspace*{-2.5mm}\\

\subsubsection{Lemma / part 2} Let $R = [r_1{:}r_2{:}r_3]_\Delta$ be a point not on a sideline of $\Delta$ and let 
$R_1 = [r_1{+\,}t_1{:\,}r_2{\,:\,}r_3]_\Delta, R_2 = [r_1{:\,}r_2{+\,}t_2{\,:\,}r_3]_\Delta, R_3= [r_1{\,:\,}r_2{\,:\,}r_3{+\,}t_3]_\Delta$ be points on the lines $R\vee A, R\vee B, R\vee C$, respectively. Then the three points $(R_2\vee R_3)\wedge a$, $(R_3\vee R_1)\wedge b, (R_1\vee R_2)\wedge c$ are collinear on the tripolar of the point $P = [t_1:t_2:t_3]_\Delta$.\\
Examples:\\
\noindent$\bullet\;\;$If Q is the perspector of a circumconic of $\Delta$ and $R_1, R_2, R_3$ are the second intersections of the lines $R\vee A, R\vee B, R\vee C$ with this circumconic, then $P$ is the Ceva point $\displaystyle[{1}/{(r_2 p_3{+}r_3 p_2)}:{1}/{(r_3 p_1{+}r_1 p_3)}:{1}/{(r_1 p_2{+}r_2 p_1})]_\Delta$ of $Q$ and $R$. \\
\noindent$\bullet\;\;$ A line $l$ is a pseudoparallel of the sideline $a'$ of the triangle $\Delta^0$ precisely when its dual $l^\delta$ is a point on the bisector of the inner angle $\angle_+CAB$ of $\Delta_0$. If we choose for $R_1, R_2, R_3$ the vertices of the first resp. second tangent-midarc-triangle, thus points on the angle bisectors,
then the coordinates of the point $P$ are\\\centerline{ $\displaystyle P=[\frac{1}{\sinh(d(I,A)-k r)}:\frac{1}{\sinh(d(I,B)-k r)}:\frac{1}{\sinh(d(I,C)-k r)}]_\Delta$},\vspace*{1mm}\\
with $r$ = radius of the incenter and $k=1$ in case of the first and $k=-1$ in case of the second tangent-midarc-triangle. The euclidean limits of these two centers are $X_{8091}$ and $X_{8092}$, but $R_1 R_2 R_3$ and $\Delta$ are, in general, not orthologic triples.\\

\subsection{Isoptics and isogonic points} 
\subsubsection{Isoptics and thaloids} Given three noncollinear points $P, Q, R$, the \textit{isoptic} (\textit{curve}) of the segment $[P;Q]_+$ through $R$ is the point set \\
\centerline{$\text{isoptic}(P,Q;R) = \{ X \,| \;d(X\vee P, X\vee Q) = d(R\vee P,R\vee Q)\}$.}\\ In the euclidean plane, such an isoptic is, in general, the union of two circles. But these two circles merge into one single circle when $\angle PRQ$ is a right angle (Thales' theorem). The situation is similar in the elliptic and in the extended hyperbolic plane. An isoptic of a segment is an algebraic curve of degree 4: Let $\boldsymbol{p} = (p_1,p_2,p_3), \boldsymbol{q}, \boldsymbol{r}, \boldsymbol{x}$ be vectors with $P = [p_1{:}p_2{:}p_3]_\Delta$, etc,  then the equation of the isoptic is \vspace*{-2mm}
\[ 
\begin{split}
&\big( (\boldsymbol{x} \times \boldsymbol{p}){{\,\scriptscriptstyle{[\mathfrak{D}]}\,}}(\boldsymbol{x} \times \boldsymbol{q})\big)^2 \big((\boldsymbol{r} \times \boldsymbol{p}){{\,\scriptscriptstyle{[\mathfrak{D}]}\,}}(\boldsymbol{r} \times \boldsymbol{p})\big)\big((\boldsymbol{r} \times \boldsymbol{q}){{\,\scriptscriptstyle{[\mathfrak{D}]}\,}} (\boldsymbol{r} \times \boldsymbol{q})\big)\\
=\;& \big( (\boldsymbol{r} \times \boldsymbol{p}){{\,\scriptscriptstyle{[\mathfrak{D}]}\,}}(\boldsymbol{r} \times \boldsymbol{q})\big)^2 \big((\boldsymbol{x} \times \boldsymbol{p}){{\,\scriptscriptstyle{[\mathfrak{D}]}\,}}(\boldsymbol{x} \times \boldsymbol{p})\big)\big((\boldsymbol{x} \times \boldsymbol{q}{{\,\scriptscriptstyle{[\mathfrak{D}]}\,}}(\boldsymbol{x} \times \boldsymbol{q})\big).\vspace*{-1mm}
\end{split}
\]
If the angle $\angle PRQ$ is right, this equation reduces to 
$(\boldsymbol{x} \times \boldsymbol{p}){{\,\scriptscriptstyle{[\mathfrak{D}]}\,}}(\boldsymbol{x} \times \boldsymbol{q}) = 0$, which is the equation of a conic but, in general, not the equation of a circle. This conic is called \textit{orthoptic} or \textit{Thales conic} or \textit{thaloid}.
\begin{figure}[!htbp]
\includegraphics[height=7.5cm]{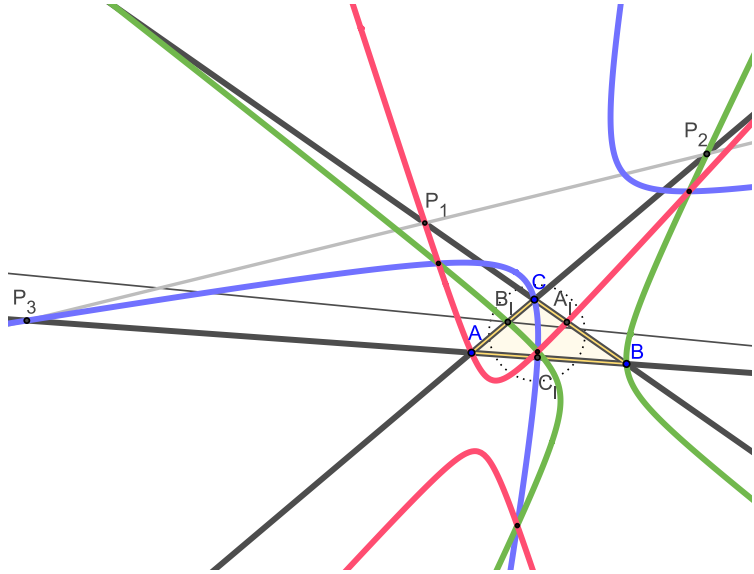}
\caption{The three apollonian thaloids of $\Delta_0$.}
\end{figure}
\begin{figure}[!hbtp]
\includegraphics[height=9.7cm]{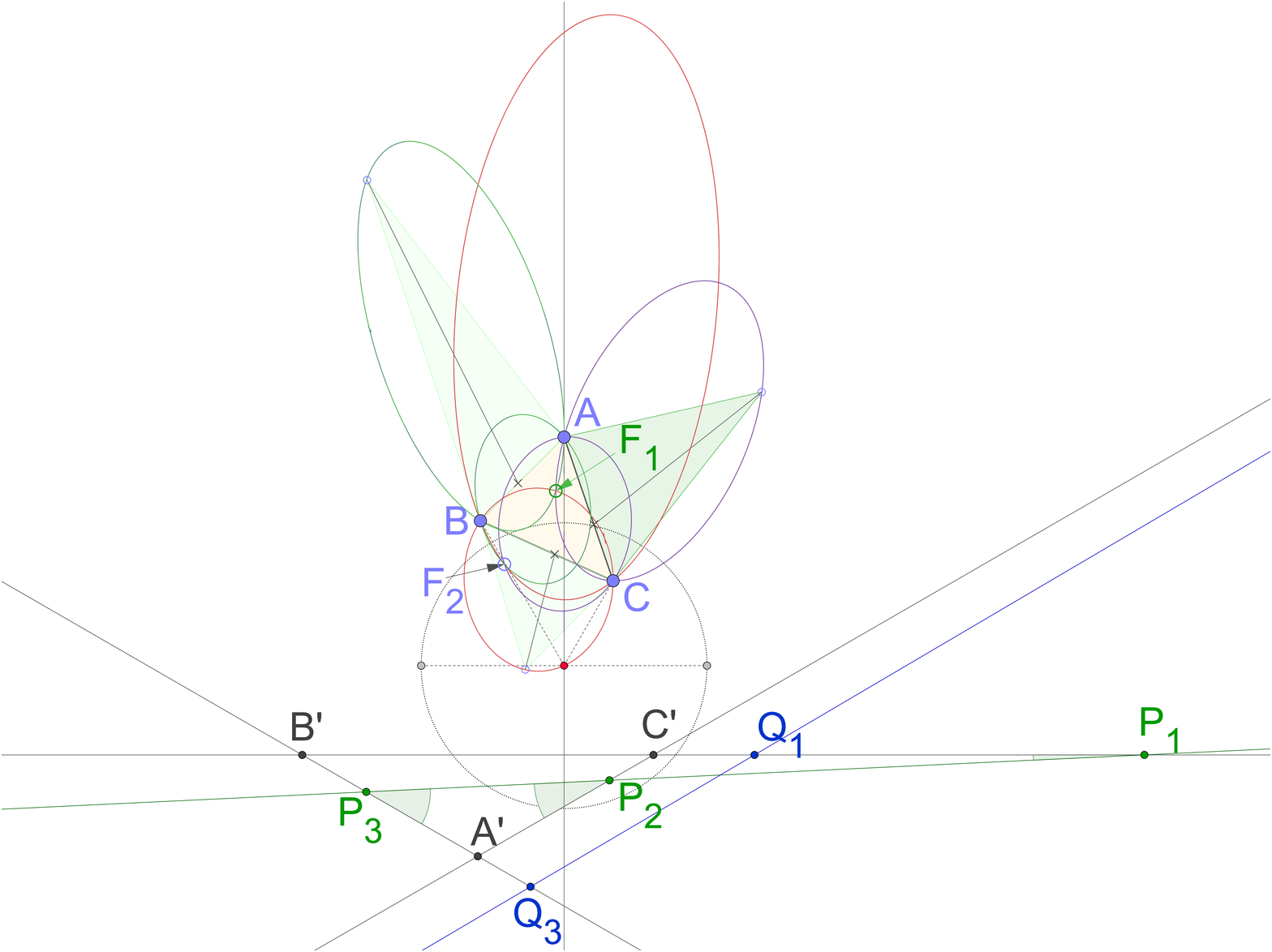}
\caption{Upper part: The isogonic points $F_1$ and $F_2$ are constructed with the help of isoptics.\;\,
 Lower part: The dual lines of $F_1$ (green) and $F_2$ (blue).}\label{fig:Fermat}
\end{figure}
\subsubsection{The apollonian thaloids of a triangle} Let $P_1,P_2,P_3$ be the points where the tripolar of the incenter $I$ meets the sidelines $a, b, c$, respectively. The three curves $\text{isoptic}(A_I,P_1;A)$, $\text{isoptic}(B_I,P_2;B)$, $\text{isoptic}(C_I,P_3;C)$ are thaloids. If two of these thaloids meet at a point, then this point is also a point of the third, see \cite{W1}. The euclidean limits of these thaloids are the apollonian circles of $\Delta_0$, so these thaloids are called \textit{apollonian thaloids} of $\Delta_0$ in \cite{W1}.\vspace*{-0.5 mm}
\subsubsection{Isogonic points of a triangle} A point $P$ is called \textit{isogonic point} of $\Delta_0$ iff $d({P\vee A}$, ${P\vee B}) = d({P\vee B},{P\vee C}) = d(P\vee C,P\vee A) = \frac{1}{3}\pi i.$ Not every triangle has an isogonic point. But in the elliptic plane every triangle has two isogonic points, and the same is true for a proper hyperbolic triangle, i.e. a triangle lying inside the absolute conic. In the euclidean case, there is a well known procedure how to find the isogonic points of a triangle by a geometric construction: The $1^{st} (2^{nd})$ isogonic point is the intersection of the circumcircles of the equilateral triangles erected outwardly (inwardly) on the sides of $\Delta_0$. A similar construction can be used in the elliptic and in the hyperbolic plane, but now using isoptics (instead of circles) and isosceles triangles with angles of measure $\frac{1}{3}\pi i$ at the apex, see Figure \ref{fig:Fermat}.\vspace*{-2.5 mm}\\

Suppose that all the sides of $\Delta_0$ are smaller than $\frac12 \pi i$ and all the angles smaller than $\frac23 \pi i$ (with respect to $\prec$) and that isogonic points exist for $\Delta_0$. Then one isogonic point lies inside the triangle and is a \textit{Fermat point} of $\Delta_0$, a point which minimizes the function $X \rightarrow d(X,A) + d(X,B) + d(X,C)$. See \cite{GH} for a proof in the elliptic case.\vspace*{-2.5 mm}\\

The lower part of Figure 14 shows a dual version of the upper. The green and the blue line are the duals of the isogonic points $F_1$ and $F_2$. Each of these lines is cut by the sidelines $a', b', c'$ of $\Delta'$ into equidistant parts. Since $F_1$ is a Fermat point of $\Delta_0$, the green line, the dual of $F_1$, minimizes the sum of the angle distances to the sidelines of $\Delta{'}$.\\ \vspace*{-1.5mm}

%the function $\ell \rightarrow d(\ell,a') + d(\ell,b') + d(\ell,c'),$
\noindent\textit{Problems}:\\
1. Determine the coordinates of the isogonic points.
2. Experiments with GeoGebra suggest that, as in the euclidean case, a point $P$ is isogonic precisely when its orthocorrespondent is identical with $P$.
A proof is missing.

\subsection{Reflection triangles}\hspace*{\fill}\vspace*{1mm} \\ 
Euclidean analogues of the theorems presented in this section can be found in \cite{HY}.
\subsubsection{The reflection-triple of a point  with respect to $\Delta$}
Given a point $P$ and a line $l$, the reflection of $P$ in $l$ will be denoted by $P^{(l)}\!$.\,\\
The reflections of $P$ in the sidelines of $\Delta$, the points $P^{(a)}, P^{(b)}, P^{(c)}$, are collinear precisely when $P$ is a point on the cubic\vspace*{0.5mm}\\
$\hspace*{10mm} x_1x_2x_3\big(\mathfrak{d}_{11}\mathfrak{d}_{22}\mathfrak{d}_{33}-4\,\mathfrak{d}_{23}\mathfrak{d}_{31}\mathfrak{d}_{12} + \sum\limits_{j=1,2,3} \mathfrak{d}_{j,j}\mathfrak{d}_{j+1,j+2}^{\;\,2}\big)\;\; \\
\hspace*{20mm}+
\;\; \sum\limits_{j=1,2,3} x_j^2\;\big(\;x_{j+1}\mathfrak{d}_{j+2,j+2}(\mathfrak{d}_{j,j+1}\mathfrak{d}_{j+1,j+2}-\mathfrak{d}_{j+1,j+1}\mathfrak{d}_{j+2,j})\\
\hspace*{40mm}             -x_{j+2}\mathfrak{d}_{j+1,j+1}(\mathfrak{d}_{j+1,j+2}\mathfrak{d}_{j+2,j}-\mathfrak{d}_{j+2,j+2}\mathfrak{d}_{j,j+1})\big)\;=\;0$\vspace*{1.5mm}\\
The euclidean limit of this cubic is the union of the circumcircle and the line at infinity.\vspace*{-1mm}\\

The triple $P^{(a)}P^{(b)}P^{(c)}$ is called \textit{the reflection-triple of} $P$ \textit{with respect to} $\Delta$.\\
The circumcenter of the reflection-triangle $(P^{(a)}P^{(b)}P^{(c)})_0$ of $P$ is the isogonal conjugate of $P$.\vspace*{-1mm}\\

The anticevian triple $\Delta^{\!P}$ and the reflection-triple of $P$ are perspective at the cevian quotient of $H$ and $P$,\vspace*{0.5mm}\\
\centerline{$[p_1(-\mathfrak{d}_{23} p_1 + \mathfrak{d}_{31} p_2 + \mathfrak{d}_{12} p_3):p_2(\mathfrak{d}_{23} p_1 - \mathfrak{d}_{31} p_2 + \mathfrak{d}_{12} p_3):p_3(\mathfrak{d}_{23} p_1 + \mathfrak{d}_{31} p_2 - \mathfrak{d}_{12} p_3)]_\Delta$.}\vspace*{0mm}\\

If $P$ is not a vertex of $\Delta$, the reflection triple of $P$ with respect to the cevian triple $\Delta_P$ is perspective to $\Delta$ at the point\vspace*{0.5mm}\\
\centerline{Q = $[p_1\big(p_1^{\;2}(\mathfrak{d}_{33}p_2^{\;2} -2\mathfrak{d}_{23}p_2p_3 +\mathfrak{d}_{22}p_3^{\;2})- \mathfrak{d}_{11}p_2^{\;2}p_3^{\;2}\big):\cdots:\cdots]_\Delta.$}\vspace*{1.5mm}\\
For $P=H$ or $P=L$, the perspector Q is a point on the orthoaxis.
\subsection{The Neuberg cubic and two related cubics} 
The reflection triple of a point $Q$ is perspective to $\Delta$ precisely when $Q$ lies on the \textit{Neuberg cubic} $\textrm{p}\mathcal{K}(K,H^{\star}),$\vspace*{0.5mm}\\
\centerline{$\sum\limits_{j=1,2,3} (\mathfrak{d}_{j,j}\mathfrak{d}_{j+1,j+2})+2\mathfrak{d}_{j+2,j}\mathfrak{d}_{j,j+1})\,x_j\,(\mathfrak{d}_{j+2,j+2}x_{\!j+1}^{\;\,2}-\mathfrak{d}_{j+1,j+1}x_{\!j+2}^{\;\,2}) = 0$.}\vspace*{0.5mm}\\
This is a self-isogonal cubic with pivot $H^{\star\!}$. On this cubic lie the incenter, the excenters, the triangle centers $H, O^+, H^\star, J_+, J_-$, the points\\
$\hspace*{4mm}W_1 = [\sinh(\alpha)(1+2(\cosh(\alpha)-\cosh(\beta)-\cosh(\gamma))):\cdots:\cdots]_\Delta.$\\
$\hspace*{4mm}W_2 = [\sinh(\alpha)(4\cosh(\alpha)(\cosh^2(\alpha)-\cosh^2(\beta)-\cosh^2(\gamma))-\cosh(\alpha)+4\cosh(\beta)\cosh(\gamma))\\
\hspace*{12mm}:\cdots:\cdots]_\Delta,$\\
$\hspace*{4mm}W_3 = [\sinh(\alpha)(4(\cosh^2(\alpha)-\cosh^2(\beta)-\cosh^2(\gamma))+1)/(\cosh(\alpha)+2\cosh(\beta)\cosh(\gamma))\\
\hspace*{12mm}:\cdots:\cdots]_\Delta$\\
and their isogonal conjugates.
$W_1$ is the \textit{Evans perspector}, the perspector of the two triples $I_1 I_2 I_3$ and $A^{(a)}B^{(b)}C^{(c)}$; its euclidean limit is $X_{484}$.
$W_2$ has euclidean limit $X_{399}$, and $W_3$ is a point on the line through $O^+$ and the isogonal conjugate of $N^+$ and has euclidean limit $X_{1157}$.\\

We introduce another cubic which consists of all points $R$ that satisfy one of the following three equivalent conditions:\\
(1) The reflection triple of $R$ and the cevian triple of $R$ are perspective.\\
(2) The reflection triple of $R$ and the triple $A^{(a)}B^{(b)}C^{(c)}$ are perspective.\\
(3) The cevian triple of $R$ and the triple $A^{(a)}B^{(b)}C^{(c)}$ are perspective.\vspace*{0.5mm}\\
This cubic is $\textrm{p}\mathcal{K}(W_4,W_5)$, a $W_4$-isocubic with pivot $W_5$,\\
$\hspace*{4mm}W_4 = [\sinh^2(\alpha)/(1-4\cosh^2(\alpha)):\cdots:\cdots]_\Delta,$\\
$\hspace*{4mm}W_5 = [\sinh(\alpha)\cosh(\alpha)/(1-4\cosh^2(\alpha)):\cdots:\cdots]_\Delta.$\\
The euclidean limits of $W_4, W_5$ are $X_{1989}, X_{265}$, respectively.\vspace*{0.5mm}

There is a bijective mapping $\textrm{p}\mathcal{K}(K,H^{\star}) \rightarrow \textrm{p}\mathcal{K}(W_4,W_5)$,  given by  $Q \mapsto R$, where $R$ is the perspector of the two triples $Q^{(a)}Q^{(b)}Q^{(c)}$, $\Delta$ and $Q $ is the perspector of the two triples $R^{(a)}R^{(b)}R^{(c)}, \Delta_R$.\vspace*{0mm}
\begin{center}
\begin{tabular} [h] {|c||c|c|c|c|c|c|c|c|c|}
\hline
\rule{0pt}{10pt}$Q$&$H$&$O^{+}$&$H^\star$&$I$&$W_{1_{\,}}$&$W_2$&$W_3$\\
\hline
\rule{0pt}{10pt}$R$&$H$&$N^+$&$W_5$&$W_6$&$W_{7_{\,}}$&$H^\star$&$W_8$\\
\hline
\end{tabular}
\end{center}\vspace*{1mm}
$\hspace*{4mm}W_6 = [\sinh(\alpha)/(1+2\cosh(\alpha)):\cdots:\cdots]_\Delta$ with euclidean limit $X_{79},$\\
$\hspace*{4mm}W_7 = [\sinh(\alpha)/(1-2\cosh(\alpha)):\cdots:\cdots]_\Delta$ with euclidean limit $X_{80},$\\
$\hspace*{4mm}W_8 = [\sinh(\alpha)/((1-4\cosh^2(\alpha))(\cosh(\alpha)+2\cosh(\beta)\cosh(\gamma)):\cdots:\cdots]_\Delta$ \\
\hspace*{13mm}with euclidean limit $X_{1141},$\\\vspace*{-1mm}

The points $P$ whose anticevian triple $\Delta^{\!P}$ is perspective to the triple $A^{(a)}B^{(b)}C^{(c)}$ form the cubic $\textrm{p}\mathcal{K}(K,N^+)$, the isogonal cubic with pivot $N^+$. There is a bijective mapping $\textrm{p}\mathcal{K}(K,N^+) \rightarrow \textrm{p}\mathcal{K}(K,H^{\star})$, 
$P \mapsto Q =$ perspector of $\Delta^{\!P}$ and $A^{(a)}B^{(b)}C^{(c)}$.\vspace*{-1mm}\\
%\centerline{$\begin{pmatrix}q_1\\q_2\\ q_3\end{pmatrix} = \begin{pmatrix}0&\;\mathfrak{d}_{33}\mathfrak{e}_{12}&-\mathfrak{d}_{22}\mathfrak{e}_{13}\\-\mathfrak{d}_{33}\mathfrak{e}_{21} &0&\;\mathfrak{d}_{11}\mathfrak{e}_{23}\\\;\mathfrak{d}_{22}\mathfrak{e}_{31}&-\mathfrak{d}_{11}\mathfrak{e}_{32}&0 \end{pmatrix}\begin{pmatrix}p_1\\p_2\\ p_3\end{pmatrix},$}\vspace*{1mm}\\
%\centerline{with $\mathfrak{e}_{i,j} = \mathfrak{d}_{i,i}\mathfrak{d}_{j,j}-2\mathfrak{d}_{i,i}\mathfrak{d}_{i+1,i+2}+2\mathfrak{d}_{j,j}\mathfrak{d}_{j+1,j+2}.$}\vspace*{1mm}\\  
%$Q$ is the perspector of $\Delta^{\!P}$ and $A^{(a)}B^{(b)}C^{(c)}$.\vspace*{1mm}
%$\begin{pmatrix}0&\mathfrak{d}_{33}(\mathfrak{d}_{11}\mathfrak{d}_{22}-2\mathfrak{d}_{11}\mathfrak{d}_{23}+2(\mathfrak{d}_{22}\mathfrak{d}_{31})&-\mathfrak{d}_{22}(\mathfrak{d}_{11}\mathfrak{d}_{33}-\mathfrak{d}_{11}\mathfrak{d}_{23}+2(\mathfrak{d}_{33}\mathfrak{d}_{12})\\-\mathfrak{d}_{33}(\mathfrak{d}_{22}\mathfrak{d}_{11}-2\mathfrak{d}_{11}\mathfrak{d}_{23}+2(\mathfrak{d}_{11}\mathfrak{d}_{23}) &0&\mathfrak{d}_{11}(\mathfrak{d}_{22}\mathfrak{d}_{33}-2\mathfrak{d}_{22}\mathfrak{d}_{31}+2\mathfrak{d}_{33}\mathfrak{d}_{12})\\ \mathfrak{d}_{22}(\mathfrak{d}_{33}\mathfrak{d}_{11}-2\mathfrak{d}_{33}\mathfrak{d}_{12}+2\mathfrak{d}_{11}\mathfrak{d}_{23})&-\mathfrak{d}_{11}(\mathfrak{d}_{33}\mathfrak{d}_{22}-2\mathfrak{d}_{33}\mathfrak{d}_{12}+2\mathfrak{d}_{22}\mathfrak{d}_{31})&0 \end{pmatrix}.
\begin{center}
\begin{tabular} [h] {|c||c|c|c|c|c|c|c|c|}
\hline
\rule{0pt}{10pt}$P$&$H$&$O^{+}$&$W_{9_{\,}}$&$I$&$W_{10}$&$N^+$&$W_{11}$\\
\hline
\rule{0pt}{10pt}$Q$&$H$&$W_2$&$O^+$&$W_{1_{\,}}$&$I$&$H^\star$&$W_3$\\
\hline
\end{tabular}
\end{center}\vspace*{1mm}
$\hspace*{4mm}W_9\;\, = [\sinh(\alpha)(4\cosh(\alpha)(\cosh^2(\alpha)-\cosh^2(\beta)-\cosh^2(\gamma))-\cosh(\alpha)-4\cosh(\beta)\cosh(\gamma))\\
\hspace*{12mm}\;\,:\cdots:\cdots]_\Delta$ with euclidean limit $X_{195},$\\
$\hspace*{4mm}W_{10} = [\sinh(\alpha)(1+2(-\cosh(\alpha)+\cosh(\beta)+\cosh(\gamma))):\cdots:\cdots]_\Delta$ \\
$\hspace*{15mm}$with euclidean limit $X_{3336}$,\vspace*{0.5mm}\\
$\hspace*{4mm}W_{11} = \;$ isogonal conjugate of $N^+$.\\

\subsection{Substitutes for the Euler line}
\subsubsection{The Euler line in euclidean geometry} In euclidean geometry, the centroid $G$ and the orthocenter $H$ of a triangle have a common cevian circle, the nine-point-circle (or Euler circle). $G, H$ and the center $N$ of this circle are collinear. They lie on the Euler line, together with the circumcenter of the triangle and several other interesting triangle centers.\\
In elliptic and in hyperbolic geometry, the points $G$ and $H$ do not have a common cevian circle, so there is no direct analog of an Euler line in these geometries,
but there are several central lines that can serve as a substitute. We list four of these.
\subsubsection{The orthoaxis}One of these lines, the orthoaxis, we treated already in subsection \ref{subsubsec:orthoaxis}. As was shown by N. Wildberger \cite{W1}, there are triangle centers with euclidean limits $X_i, i = 2,3,4,20,30$ lying on this line. So Vigara proposed to call it Euler-Wildberger line.
\subsubsection{The orthoaxis of the medial triangle}The orthoaxis of the medial triangle $\Delta_{G,0}$ is the line $G\vee O$. Besides $G$ and $O$, it passes through the isogonal conjugate of O, through the Lemoine-isoconjugate of O, through L and several other points of interest.
A more detailed description of this line is given in \cite{Ev}.\vspace*{-1.5mm}\\

The line $G\vee O$ can also be a substitute for an Euler line of the anticevian triangle $\Delta^G_{\;\,0}$ of $G$. $G$ can serve as a pseudo-centroid, $\mathscr{C}_0$ as a pseudo-Euler-circle and the lines $A^G\vee A', B^G\vee B', C^G\vee C'$ as pseudo-altitudes of $\Delta^G_{\;\,0}$. These pseudo-altitudes meet at a point on $G\vee O$ which is the circumcevian conjugate of $G$ with respect to $\Delta^{G}_{\;\,0}$.\vspace*{0.5mm}\\
\textit{Problem}: Taking $\Delta^G$ as a reference triple, what are the coordinates of the point $G$?\vspace*{0mm}
\subsubsection{The line through $G$ and $H$}The line $G\vee H$ has the equation
\[
\begin{split}
\sum \limits_{j=1,2,3} (\mathfrak{c}_{j+2,j} \mathfrak{c}_{j,j+1} - \mathfrak{c}\,\mathfrak{c}_{j+1,j+2})(1 + \mathfrak{c}\,\mathfrak{c}_{j+1,j+2}) (\mathfrak{c}_{j,j+1} - \mathfrak{c}_{j+2,j}) x_j   &=  0,\\
\!\!\!\!\!\!\!\!\!\!\!\!\!\!\!\!\!\!\!\!\!\!\!\!\!\!\!\!\!\!\!\!\!\!\!\;\;\;\;\;\;\;\;\textrm{which is equivalent to}  \;\;\;\;\;\;\;\;\;\;\;\sum \limits_{j=1,2,3} \mathfrak{d}_{j+1,j+2}(\mathfrak{d}_{j+2,j} - \mathfrak{d}_{j,j+1}) x_j  &=  0.
\end{split}
\]
It passes through the center of a circle which touches internally the incircle and externally the excircles of $\Delta_0$. This was shown for a triangle on a sphere in 1864 by G. Salmon \cite{Sa}, and Salmon also calculated the trilinear coordinates of this center, which we denote by $\hat{N}$.  It has barycentric coordinates:\vspace*{1mm}\\
$\hspace*{9.3mm} \hat{N} = [\mathfrak{d}_{31} \mathfrak{d}_{12}-{\det}(\mathfrak{D}), \mathfrak{d}_{12} \mathfrak{d}_{23}-{\det}(\mathfrak{D}) ,\mathfrak{d}_{23} \mathfrak{d}_{31} - {\det}(\mathfrak{D})]_\Delta$\vspace*{0.5mm}\\
$\hspace*{18mm} = [(1 + \mathfrak{c}\,\mathfrak{c}_{23}^{\;}) (\, \mathfrak{c}\,(\mathfrak{c}_{31}^{\;2} - \mathfrak{c}_{31}^{\;} \mathfrak{c}_{12}^{\;} + \mathfrak{c}_{12}^{\;2} + 1) + \mathfrak{c}_{23}^{\;} (\mathfrak{c}_{31}^{\;} \mathfrak{c}_{12}^{\;} - 1)):\cdots :\cdots]_\Delta$\vspace*{0.8mm}\\
$\hspace*{18mm} = [\sinh(\alpha)\cosh(\beta-\gamma):\sinh(\beta)\cosh(\gamma-\alpha):\sinh(\gamma)\cosh(\alpha-\beta)]_\Delta.$\\

An equation of the circle in trilinear coordinates was given for the elliptic (spherical) case in 1861 by A.S. Hart \cite{Ha}; Salmon therefore used the name \textit{Hart's circle}.\vspace*{1mm}\\
The Hart circle is a conic $\mathcal{C}(\mathfrak{M})$ with $\displaystyle \mathfrak{m}_{11} = \frac{\mathfrak{c} (\mathfrak{c}_{12} - 2 \mathfrak{c}_{23} + \mathfrak{c}_{31}) + \mathfrak{c}_{12} \mathfrak{c}_{31} - 1}{\mathfrak{c}\,\mathfrak{c}_{23} + 1}, \mathfrak{m}_{23} = 1- \mathfrak{c} \mathfrak{c}_{23}, \cdots.$\vspace*{1mm}\\
Salmon has shown the following relation between the radius $\rho_H$ of the Hart circle and the radius $R$ of the circumcircle: $\tanh(\rho_H) = \frac{1}{2} \tanh(R)$.\vspace*{-0.5mm}\\

The touchpoint of the incircle $\mathscr{I}_k$ with the Hart circle we denote by $F_{0,k}$;  these \textit{Feuerbach points} of triangle $\Delta_0$ have coordinates\vspace*{1mm}\\
$F_{0,0} = [\mathfrak{d}_{31} \mathfrak{d}_{12}{\,-\,}{\det}(\mathfrak{D}){\,-}f\sqrt{\mathfrak{d}_{11}}:\mathfrak{d}_{12} \mathfrak{d}_{23}{\,-\,}{\det}(\mathfrak{D}){\,-}f\sqrt{\mathfrak{d}_{22}}:\mathfrak{d}_{23} \mathfrak{d}_{31}{\,-\,}{\det}(\mathfrak{D}){\,-}f\sqrt{\mathfrak{d}_{33}}]_\Delta,\vspace*{0mm}$\\ 
$\displaystyle \hspace*{6.8mm}= [\sinh(\alpha)\sinh^2(\frac{\beta-\gamma}{2}):\sinh(\beta)\sinh^2(\frac{\gamma-\alpha}{2}):\sinh(\gamma)\sinh^2(\frac{\alpha-\beta}{2})]_\delta,\vspace*{0.5mm}$\\
$\hspace*{9.5mm}f = \sqrt{\mathfrak{d}_{12}}\sqrt{\mathfrak{d}_{23}}\sqrt{\mathfrak{d}_{31}}\; {\det}^2(\mathfrak{D}),\vspace*{0.5mm}$\\
{$F_{0,1} = [\mathfrak{d}_{31} \mathfrak{d}_{12}{\,-\,}{\det}(\mathfrak{D}){\,-\,}f\sqrt{\mathfrak{d}_{11}}:\mathfrak{d}_{12} \mathfrak{d}_{23}{\,-\,}{\det}(\mathfrak{D}){\,+\,}f\sqrt{\mathfrak{d}_{22}}:\mathfrak{d}_{23} \mathfrak{d}_{31}{\,-\,}{\det}(\mathfrak{D}){\,+\,}f\sqrt{\mathfrak{d}_{33}}]_\Delta,$}\vspace*{1mm}\\
 $\hspace*{5.5mm}$etc.\vspace*{1mm}\\

\noindent The triple $F_{0,1}F_{0,2}F_{0,3}$ is perspective to $\Delta$; the perspector is the point\vspace*{1mm}\\
\centerline{\hspace*{5mm}$[\mathfrak{d}_{31} \mathfrak{d}_{12}{\,-\,}{\det}(\mathfrak{D}){\,+}f\sqrt{\mathfrak{d}_{11}}:\mathfrak{d}_{12} \mathfrak{d}_{23}{\,-\,}{\det}(\mathfrak{D}){\,+}f\sqrt{\mathfrak{d}_{22}}:\mathfrak{d}_{23} \mathfrak{d}_{31}{\,-\,}{\det}(\mathfrak{D}){\,+}f\sqrt{\mathfrak{d}_{33}}]_\Delta,\,$}\vspace*{0.5mm}\\
$\hspace*{0.7mm}\,=\,[\sinh(\alpha)\cosh^2(\frac{\beta-\gamma}{2}):\sinh(\beta)\cosh^2(\frac{\gamma-\alpha}{2}):\sinh(\gamma)\cosh^2(\frac{\alpha-\beta}{2})]_\Delta.$\\
It is the harmonic conjugate of $F_{0,0}$ with respect to $\{I_0 , \hat{N}\}$ and has euclidean limit $X_{\!12}$.\vspace*{-2mm}\\

The dual of the Hart circle is a circle with center $\hat{N}$ which internally touches the circumcircle $\mathscr{C}_0$  and externally the circumcircles $\mathscr{C}_k, k = 1,2,3$ , with touchpoints lying on the lines $\hat{N} \vee O_k, k = 0,\cdots,3$, respectively.

\begin{figure}[!htbp]
\includegraphics[height=9cm]{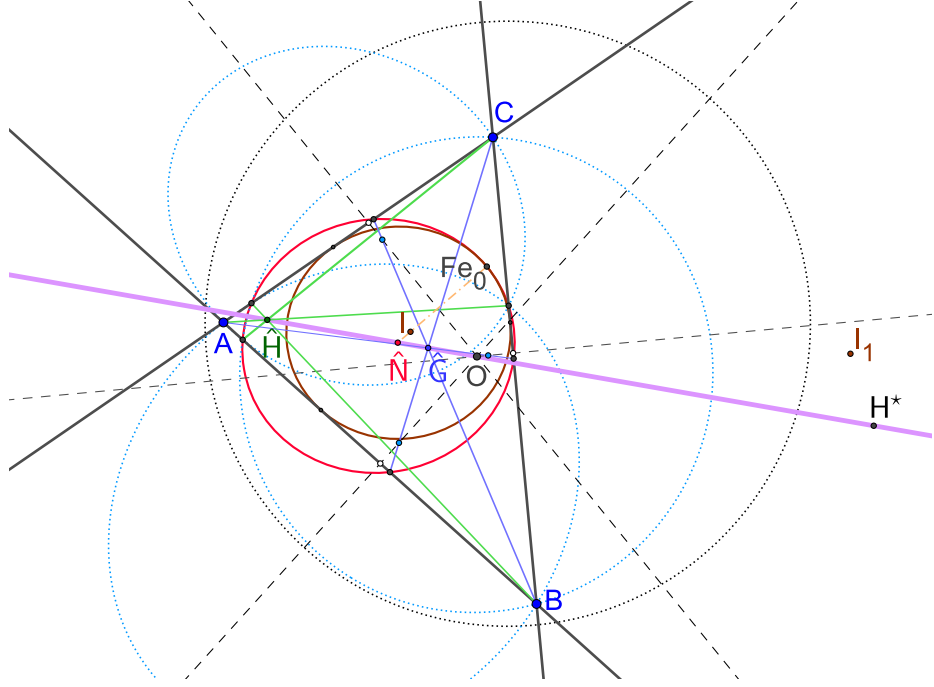}
\caption{The Akopyan line, the Hart circle (red) and incircle (brown)}
\end{figure}
\subsubsection{The Akopyan line}
The line $O\vee H^{\!\star}$ has the equation\vspace*{0.5mm}\\
%\[ (\mathfrak{c}_{31} - \mathfrak{c}_{12}) (2 \mathfrak{c}_{23} - \mathfrak{c}_{31} - \mathfrak{c}_{12} + c(\mathfrak{c}_{12} \mathfrak{c}_{31}-1)) x_1 + (\mathfrak{c}_{12} - \mathfrak{c}_{23}) (2 \mathfrak{c}_{31} - \mathfrak{c}_{12} - \mathfrak{c}_{23} + c(\mathfrak{c}_{23} \mathfrak{c}_{12}-1)) x_2
%+ (\mathfrak{c}_{23} - \mathfrak{c}_{31}) (2 \mathfrak{c}_{12} - \mathfrak{c}_{23} - \mathfrak{c}_{31} + c(\mathfrak{c}_{12} \mathfrak{c}_{23}-1)) x_3 = 0\]
%\[ \sum \limits_{j=1,2,3} \big(\mathfrak{c}_{j+2,j} - \mathfrak{c}_{j,j+1}\big) \big(2 \mathfrak{c}_{j+1,j+2} - \mathfrak{c}_{j+2,j} - \mathfrak{c}_{j,j+1} + \mathfrak{c}\,(\mathfrak{c}_{j,j+1} \mathfrak{c}_{j+2,j}-1)\big) x_j  = 0
%\]
\centerline{$\sum \limits_{j=1,2,3} \big(\mathfrak{c}_{j+2,j} - \mathfrak{c}_{j,j+1}\big) \big((1-\mathfrak{c}\,\mathfrak{c}_{j+2,j})(1-\mathfrak{c}\,\mathfrak{c}_{j,j+1})-2(1-c\,\mathfrak{c}_{j+1,j+2})\big) x_j  = 0.$}\\

There are several important triangle centers on this line. First of all, a point $\hat{G}$ whose cevian lines divide the triangle area in equal parts. The existence of such a point was already shown (for a spherical triangle) 1827 by J. Steiner \cite{St}. 
The coordinates of this point are\vspace*{1.5mm}
%	\[ 
%\begin{split}		
%	&[\frac {\sqrt{1{+}c\,\texttt{c}_{23}}}{\sqrt{2}\,\sqrt{1{+}c\,\texttt{c}_{23}}+\sqrt{1{+}c\,\texttt{c}_{31}}\,\sqrt{1{+}c\,\texttt{c}_{12}}}:\frac {\sqrt{1{+}c\,\texttt{c}_{31}}}{\sqrt{2}\,\sqrt{1{+}c\,\texttt{c}_{31}}+\sqrt{1{+}\texttt{c}_{12}}\,\sqrt{1{+}c\,\texttt{c}_{23}}\\
%	&\hspace*{45mm}:\frac {\sqrt{1{+}\texttt{c}_{12}}}{\sqrt{2}\,\sqrt{1{+}\texttt{c}_{12}}+ \sqrt{1{+}c\,\texttt{c}_{23}}\,\sqrt{1{+}\texttt{c}_{31}}}] 
%\end{split}	
%\]
$\hspace*{5mm}\displaystyle [\frac {\sqrt{1+\mathfrak{c}\,\mathfrak{c}_{23}}}{\sqrt{2}\,\sqrt{1+\mathfrak{c}\,\mathfrak{c}_{23}}+\sqrt{1+\mathfrak{c}\,\mathfrak{c}_{31}}\,\sqrt{1+\mathfrak{c}\,\mathfrak{c}_{12}}}\,\;:\;\frac {\sqrt{1+\mathfrak{c}\,\mathfrak{c}_{31}}}{\sqrt{2}\,\sqrt{1+\mathfrak{c}\,\mathfrak{c}_{31}}+\sqrt{1+\mathfrak{c}\,\mathfrak{c}_{12}}\,\sqrt{1+\mathfrak{c}\,\mathfrak{c}_{23}}}$\vspace*{0.5mm}\\
	$\hspace*{4mm}: 
	\displaystyle\frac {\sqrt{1+\mathfrak{c}\,\mathfrak{c}_{12}}}{\sqrt{2}\,\sqrt{1+\mathfrak{c}\,\mathfrak{c}_{12}}+ \sqrt{1+\mathfrak{c}\,\mathfrak{c}_{23}}\,\sqrt{1+\mathfrak{c}\,\mathfrak{c}_{31}}}]_\Delta$	\\
	
$\displaystyle	\;=\; [\frac {\cosh(\frac{1}{2} \mathscr{a})}{\cosh(\frac{1}{2} \mathscr{a})+\cosh(\frac{1}{2} \mathscr{b})\cosh(\frac{1}{2} \mathscr{c})}$\\
$\hspace*{8mm}\displaystyle:\frac {\cosh(\frac{1}{2} \mathscr{b})}{\cosh(\frac{1}{2} \mathscr{b})+\cosh(\frac{1}{2} \mathscr{c})\cosh(\frac{1}{2} \mathscr{a})}
: \frac {\cosh(\frac{1}{2} \mathscr{c})}{\cosh(\frac{1}{2} \mathscr{c})+\cosh(\frac{1}{2} \mathscr{a})\cosh(\frac{1}{2} \mathscr{b})}]_\Delta.$\\\vspace*{1mm}

Another point on $O\vee H^{\!\star}$ is the center $\hat{N}$ of the Hart circle. Akopyan \cite{Ak} proved in 2011 that this Hart circle is a circumcircle of the traces of $\hat{G}$ and also a circumcircle of the traces of another point $\hat{H}$, which he called \textit{pseudo-orthocenter}. Moreover, he showed that the three points $\hat{G}, \hat{N}, \hat{H}$ are collinear on a line which passes also through the circumcenter $O$. It is easy to check by calculation that this line is $O\vee H^{\!\star}$. 
Vigara proposed to call this line \textit{Akopyan-Euler-line}. Akopyan used the name \textit{Euler line}.\vspace*{-2mm}\\

We introduce the \textit{Akopyan-measure} of an inner angle of a triangle as the sum of the other two inner angles diminished by a half of the triangle area \footnote{$^)$ See \cite{Ma}.}$^)$. For example, the Akopyan-measure of the angle $\angle_+ ACB$ of triangle $\Delta_0$ is $\hat{\mu}(\angle_+ ACB) = \frac{1}{2} (\alpha{\,+\,}\beta{\,-\,}\gamma{\,+\,}\pi i).$ By adding up the Akopyan-measures of the three inner angles of triangle $\Delta_0$ we get $\,\frac {1}{2}\textrm{area}(\Delta_0) + 2\pi i$. \vspace*{-2mm}

\noindent With the help of the Akopyan-measure we can formulate the following version of the Inscribed Angle Theorem, cf. \cite{Ak}: \vspace*{0.5mm}\\
$\hat{\mu}(\angle_+ ACB) = \frac {1}{2}\pi i \pm \mu(\angle_+ OAB)$, depending on whether the circumcenter $O$ lies outside the triangle $\Delta_0$ or not.
\vspace*{-2mm}\\

We like to shortly present Akopyan's explanation for $\hat{H}$ having properties similar to those of the orthocenter. For this we also make use of the Akopyan measure:\\
While the altitude $C{\vee}H$ meets the line $c$ at the point $C_{\!H}$ such that the angles $\angle_+ (AC_{\!H}C)$ and $\angle_+ (BC_{\!H}C)$ have both measures equal to $\frac{1}{2}\pi i$, the pseudo-altitude
$C{\vee}\hat{H}$ meets the line $c$ in $C_{\!\hat{H}}$ such that $\hat{\mu}(\angle_+ AC_{\!\hat{H}}C) = \hat{\mu}(\angle_+ BC_{\!\hat{H}}C) = \frac{1}{4} (\pi i{+}\alpha{+}\beta{+}\gamma).$ \vspace*{-1mm}\\
%((The same applies to the other two cevians of $\hat{H}$.))

The coordinates of $\hat{H}$ are $\hspace*{5mm}\displaystyle[\frac {\sqrt{1+\mathfrak{c}\,\mathfrak{c}_{23}}}{\sqrt{2}\,\sqrt{1+\mathfrak{c}\,\mathfrak{c}_{23}}-\sqrt{1+\mathfrak{c}\,\mathfrak{c}_{31}}\,\sqrt{1+\mathfrak{c}\,\mathfrak{c}_{12}}}{:}\cdots{:}\cdots]_\Delta $\\
 
$\;\hspace*{40.5mm}\;= \;\displaystyle[\frac {\cosh(\frac{1}{2} \mathscr{a})}{\cosh(\frac{1}{2} \mathscr{a})-\cosh(\frac{1}{2} \mathscr{b})\cosh(\frac{1}{2} \mathscr{c})}{:}\cdots{:}\cdots]_\Delta.$\vspace*{3.5mm}\\
%$\;\hspace*{10mm}\;\displaystyle\frac {\sqrt{1+\mathfrak{c}\,\mathfrak{c}_{31}}}{\sqrt{2}\,\sqrt{1+\mathfrak{c}\,\mathfrak{c}_{31}}-\sqrt{1+\mathfrak{c}\,\mathfrak{c}_{12}}\,\sqrt{1+\mathfrak{c}\,\mathfrak{c}_{23}}}
%	:\frac {\sqrt{1+\mathfrak{c}\,\mathfrak{c}_{12}}}{\sqrt{2}\,\sqrt{1+\mathfrak{c}\,\mathfrak{c}_{12}}-\sqrt{1+\mathfrak{c}\,\mathfrak{c}_{23}}\,\sqrt{1+\mathfrak{c}\,\mathfrak{c}_{31}}}]_\Delta.$

The points $\hat{H},\hat{N},\hat{G}, O$ form a harmonic range.\\
\noindent\textit{Proof}: We use the following abbreviations: $c_1\!:=\cosh(\mathscr{a}/2),\,c_2\!:=\cosh(\mathscr{b}/2),\,c_3\!:=\cosh(\mathscr{c}/2), 
\,s_1 :=\sinh(\mathscr{a}/2), s_2 :=\sinh(\mathscr{b}/2), s_3 :=\sinh(\mathscr{c}/2).$ Define vectors $\boldsymbol{g}, \boldsymbol{h}, \boldsymbol{n}, \boldsymbol{o}$
\noindent by\\
$\noindent\displaystyle\hspace*{8mm}\boldsymbol{g} = (g_1,g_2,g_3) ,\; g_j = \frac{c_j}{c_j+c_{j+1} c_{j+2}} ,\;\;\; \boldsymbol{h} = (h_1,h_2,h_3) ,\; h_j = \frac{c_j}{c_j-c_{j+1} c_{j+2}}$,\vspace*{0.5mm}\\

\noindent$\displaystyle\hspace*{8mm}\boldsymbol{n} = (n_1,n_2,n_3) ,\; n_j = c_j^{\;2}\big(c_j^{\;2}(2c_{j+1}^{\;2}c_{j+2}^{\;2}-c_{j+1}^{\;2}-c_{j+2}^{\;2})-c_{j+1}^{\;2}s_{j+1}^{\;2}-c_{j+2}^{\;2}s_{j+2}^{\;2}\big)$\vspace*{0.5mm}\\
$\noindent\displaystyle\hspace*{8mm}\boldsymbol{o} = (o_1,o_2,o_3) ,\; o_j = -s_j^{\;2}(1+c_j^{\;2}-c_{j+1}^{\;2}-c_{j+2}^{\;2})$,\vspace*{1.5mm}\\ 
\noindent and define real numbers $r, s, t$ by\\ 
$\noindent\displaystyle\hspace*{8mm}r = c_1 c_2 c_3$,\vspace*{0.5mm}\\
$\noindent\displaystyle\hspace*{8mm}s = (c_1 + c_2 c_3)(c_2 + c_3 c_1)(c_3 + c_1 c_2)(1+2r  - c_1^{\;2}-c_2^{\;2}-c_3^{\;2})$,\vspace*{0.5mm}\\
$\noindent\displaystyle\hspace*{8.5mm}t = (c_1 - c_2 c_3)(c_2 - c_3 c_1)(c_3 - c_1 c_2)(1-2r  - c_1^{\;2}-c_2^{\;2}-c_3^{\;2})$.\vspace*{1mm}\\ 
Then $\hat{G} = [g_1{:}g_2{:}g_3]_\Delta, \hat{H} = [h_1{:}h_2{:}h_3]_\Delta, \hat{N} = [n_1{:}n_2{:}n_3]_\Delta, O = [o_1{:}o_2{:}o_3]_\Delta\vspace*{1mm}$\\
and  $\hspace*{32.5mm}\;r\boldsymbol{o} + \boldsymbol{n} = s \boldsymbol{g}\;\;\;$,$\;\;\;r\boldsymbol{o} - \boldsymbol{n} = t \boldsymbol{h}.\hspace*{32mm}\Box$\vspace*{0.5mm}\\

The cevian line $C\vee \hat{G}$ meets the perpendicular bisector $C_G\vee C'$ at a point with coordinates
\[ 
\begin{split}	
[&\frac {\sqrt{1+\mathfrak{c}\,\mathfrak{c}_{23}}}{\sqrt{2}\,\sqrt{1+\mathfrak{c}\,\mathfrak{c}_{23}}+\sqrt{1+\mathfrak{c}\,\mathfrak{c}_{31}}\,\sqrt{1+\mathfrak{c}\,\mathfrak{c}_{12}}}:\frac {\sqrt{1+\mathfrak{c}\,\mathfrak{c}_{31}}}{\sqrt{2}\,\sqrt{1+\mathfrak{c}\,\mathfrak{c}_{31}}+\sqrt{1+\mathfrak{c}\,\mathfrak{c}_{31}}\,\sqrt{1+\mathfrak{c}\,\mathfrak{c}_{23}}}\\
	:& \frac {(1-\mathfrak{c} \mathfrak{c}_{12})\sqrt{1+\mathfrak{c}\,\mathfrak{c}_{12}}}{(\sqrt{2}\,\sqrt{1+\mathfrak{c}\,\mathfrak{c}_{23}}+\sqrt{1+\mathfrak{c}\,\mathfrak{c}_{31}}\,\sqrt{1+\mathfrak{c}\,\mathfrak{c}_{12}})(\sqrt{2}\,\sqrt{1+\mathfrak{c}\,\mathfrak{c}_{31}}+\sqrt{1+\mathfrak{c}\,\mathfrak{c}_{31}}\,\sqrt{1+\mathfrak{c}\,\mathfrak{c}_{23}})}]_\Delta.
\end{split}	
\]\vspace*{1.5mm}
This point is the center of a circle through the points $A, B, A_{\!\hat{H}}, B_{\!\hat{H}}$. \vspace*{-1.5mm}Thus, the line $A_{\!\hat{H}}\vee B_{\!\hat{H}}$ is antiparallel  to $A\vee B$ with respect to the angle $\angle_+ ACB$, and since $A_{\!\hat{H}}, B_{\!\hat{H}},A_{\!\hat{G}}, B_{\!\hat{G}}$ are points on the Hart circle, the line $A_{\!\hat{G}}\vee B_{\!\hat{G}}$ is parallel to $A\vee B\;$.\vspace*{-1mm}\\ 

We give a short hint of how to construct the cevians $\hat{G}\vee C$ and $\hat{H}\vee C$; see \cite{Ev} for a more detailed explanation.
Let $m$ be the sideline $A_G{\vee}B_G$ of the medial triangle. Define points $R_1 = \textrm{ped}(A,m)$,  $R_2 = \textrm{ped}(B,m)$, $S_1 = m \wedge \textrm{perp}(c,A)$, $S_2 = m \wedge \textrm{perp}(c,B)$. The bisectors of the angles $\angle_- R_1 A S_1$ and $\angle_- R_2 B S_2$ meet at a point on $\hat{G}\vee C$, the bisectors of the angles $\angle_+ R_1 A S_1$ and $\angle_+ R_2 B S_2$ meet at a point on $\hat{H}\vee C$.\\

The Hart circle meets each cevian line of $\hat{G}$ and $\hat{H}$  twice. One intersection point is the corresponding cevian point, the other one is also a significant point. The description of this second point generalizes as follows:\\
If $P$ and $Q$ are two points not on the sidelines of $\Delta$, then for each vertex $R\in \{A, B, C\}$, the cevian line $R{\vee}R_P$ meets the bicevian conic of  $P$ and $Q$ in a point which is the harmonic conjugate with respect to $\{R, P\}$ of the intersection point of $R\vee R_P$ with the tripolar of  $Q$.\\
\noindent\textit{Proof}: Without loss of generality we take $R = C$. For the matrix of the bicevian conic of $P = [p_1{:}p_2{:}p_3]_\Delta$ and $Q=[q_1{:}q_2{:}q_3]_\Delta$, see \ref{subsubsec:orthoaxis}. The second intersection of the cevian line $C\vee C_P$ with this conic is the point $S = [p_1 q_1 q_2:p_2 q_1 q_2:p_1 q_2 q_3{+}q_1 p_2 q_3{+}2 q_1 q_2 p_3]$, the intersection of $C\vee C_P$ with $Q^\tau$ is the point $T = [p_1 q_1 q_2:p_2 q_1 q_2:- q_3 (p_1 q_2 + p_2 q_1)]$.
It can be easily checked that $P,S,C,T$ form a harmonic range. $\Box$\vspace*{-2.5mm}\\

\noindent\textit{Addition}: The pole of $Q^\tau$ with respect to the bicevian conic of $P$ and $Q$, $[q_1 (2 p_1 q_2 q_3 + q_1 (p_2 q_3 + p_3 q_2)):\cdots:\cdots]_\Delta$, is a point on the line $P\vee Q$. If $P$ is a perspector of an inconic of $\Delta$, then its polar with respect to this conic agrees with its tripolar.\vspace*{12mm}\\

%p23 = (√\mathfrak{c}√(\mathfrak{c}(p12 + 1)^2 - 2 p12 w1^2) (w1^2 - w2^2) - \mathfrak{c}(p12 (w1^2 - w2^2) - w1^2 - w2^2) - w1^2 w2^2)/(p12 w1^2 (2 \mathfrak{c}- w1^2))
%p12 (√\mathfrak{c}(w1^2 - w2^2) √(\mathfrak{c}(p12^2 + 2 p12 + 1) - 2 p12 w1^2) + \mathfrak{c}(p12 (w1^2 - w2^2) - w1^2 - w2^2) + w1^2 w2^2)/(2 \mathfrak{c}(p12 (w1^2 - w2^2) - w2^2) + w2^4)
%2(w1^2 - w2^2) (√( (p12^2 + 2 p12 + 1) - 2 \mathfrak{c}p12 w1^2) +\mathfrak{c}p12) - ( -\mathfrak{c}w1^2 - \mathfrak{c}w2^2 + w1^2 w2^2)


\begin{thebibliography}{99}   



\bibitem {Ak} A.$\,$V. Akopyan, On some classical constructions extended to hyperbolic geometry,\\ arxiv:1105.2153, 2011.

\bibitem {DD} E. Danneels, N. Dergiades, A Theorem on Orthology Centers, \textit{Forum Geom.} 4 (2004), 135-141, 5-13.

\bibitem {DP} S. Dominte, T.D. Popescu, Concyclicities in Tucker-like configurations, available at\\ 
\texttt{https://www.awesomemath.org/wp-pdf-files/math-reflections/\\mr-2016-01/article1\_concyclicities\_tucker\_configurations.pdf}.

\bibitem {D} D. Douillet, Translation of the Kimberling's Glossary into barycentrics,
available at\\ \texttt{http://www.douillet.info/~douillet/triangle/Glossary.pdf}.


\bibitem {Ev} M. Evers, On centers and central lines of triangles in the elliptic plane, arXiv:1705.06187v3, 2018.

\bibitem {EE} L. Emelyanov, T. Emelyanova, A Note on the Schiffler Point, \textit{Forum Geom.} 3 (2003) 113-116.


\bibitem {GH} K. Ghalieh, M. Hajja, The Fermat Point of a Spherical Triangle, \textit{Math. Gaz.} 80 (1996), 561-564.

\bibitem {Gi} B. Gibert, Cubics in the triangle plane, available at\\ 
\texttt{http://bernard.gibert.pagesperso-orange.fr}.

\bibitem {GY} D. Grinberg, P. Yiu, The Apollonius Circle as a Tucker Circle, \textit{Forum Geom.} 2 (2002) 175-182.

\bibitem {Gr} D. Grinberg, Three properties of the symmedian point, 
available at\\ \texttt{http://www.cip.ifi.lmu.de/~grinberg/TPSymmedian.pdf}.

\bibitem {Ha} A.S. Hart, Extension of Terquem's theorem respecting the circle which bisects three sides of a triangle, \textit{Quarterly J. of Math.} 4 (1861) , 260-261.

\bibitem {HY} A.P. Hatzipolakis, P. Yiu, Reflections in Triangle Geometry, \textit{Forum Geom.} 9 (2009) 301–348.            

\bibitem {H1} \'A.$\,$G.$\,$Horv\'ath, On the hyperbolic triangle centers, arxiv:1410.6735v1, 2014.

\bibitem {H2} \'A.\,G.$\,$Horv\'ath, Hyperbolic plane-geometry revisited, arxiv:1405.1068v2, 2014.


\bibitem {Ki} C. Kimberling, Triangle centers as functions,  \textit{Rocky Mt. J. Math.} 23 (1993) 1269-1286.

\bibitem {ETC} C.Kimberling,  \textit{Encyclopedia of Triangle Centers} (ETC), available at\\ 
\texttt{http://faculty.evansville.edu/ck6/encyclopedia/ETC.html}.

\bibitem {KY} S.N. Kiss, P. Yiu, On the Tucker circles,  \textit{Forum Geom.} 17 (2017), 157-175.

\bibitem {L1} F. van Lamoen, Napoleon Triangles and Kiepert Perspectors, \textit{Forum Geom.} 3 (2003), 65-71.

\bibitem {L2} F. van Lamoen, Some Concurrencies from Tucker Hexagons, \textit{Forum Geom.} 2 (2002), 5-13.

\bibitem {Ma} P. Maraner, Fate of the Euler Line and the Nine-Point Circle on the Sphere, Contribution to the Wolfram Demonstrations Project, 2017, available at\\ 
\texttt{http://demonstrations.wolfram.com/ \newline
  FateOfTheEulerLineAndTheNinePointCircleOnTheSphere/}

\bibitem {Ru} R.A. Russell, Non-Euclidean Triangle Centers, arXiv:1608.08190v2, 2017.


\bibitem {Sa} G. Salmon, On the circle which touches the four circles which touch the sides of a given spherical triangle, \textit{Q. J. Math.} 6 (1864), 67-73.

\bibitem {St} J. Steiner, Verwandlung und Theilung sph\"arischer Figuren durch Construction,             
\textit{J. Reine Angew. Math.} 2 (1827), 45-63.

\bibitem {TL} I. Todhunter, J.G. Leathem, \textit{Spherical Trigonometry},  Macmillan \& Co. Ltd., 1914.

\bibitem {U1} A.$\,$A. Ungar, \textit{Hyperbolic Triangle Centers: The special relativistic approach}, Springer, New York, 2010.	

\bibitem {U2} A.$\,$A. Ungar, On the Study of Hyperbolic Triangles and Circles by Hyperbolic Barycentric Coordinates in Relativistic Hyperbolic Geometry,  	
arXiv:1305.4990, 2013.

\bibitem {Vi} R. Vigara, Non-euclidean shadows of classical projective theorems, arxiv:1412.7589, 2014.


\bibitem {W1} N.$\,$J. Wildberger, Universal Hyperbolic Geometry III: First Steps in Projective Triangle Geometry, \textit{KoG} 15(2011), 25-49.	

\bibitem {W2} N.$\,$J. Wildberger, A. Alkhaldi: Universal Hyperbolic Geometry IV: Sydpoints and Twin Circumcircles, \textit{KoG} 16 (2012), 43-62.


\bibitem {Y} P. Yiu, \textit{Introduction to the Geometry of the Triangle}, Florida Atlantic University Lecture Notes, 2001.


\bibitem {GG} GeoGebra, Ein Softwaresystem f\"ur dynamische Geometry und Algebra, invented by M. Hohenwarter and currently developed by IGI.\vspace*{3mm}		
				
\end{thebibliography}
\end{document}